\documentclass[12pt]{article}
\usepackage{amsfonts,amsmath,amsthm,amscd,amssymb,latexsym,amsbsy}
\usepackage[T1]{fontenc}
\usepackage{chngcntr,xcolor}
\definecolor{navyblue}{rgb}{0.0, 0.0, 0.5} 
\definecolor{my_red}{rgb}{1.0, 0.08, 0.0} 
\usepackage[colorlinks=true, citecolor = navyblue,  linkcolor=navyblue, filecolor=magenta, urlcolor=navyblue]{hyperref}
\usepackage[mathscr]{euscript}
\usepackage{enumitem}
\usepackage{geometry}
\usepackage{lipsum}        
\usepackage{bbm} 
\usepackage{graphicx}
\usepackage[backend=bibtex,style=authoryear, useprefix=false, maxcitenames=2, maxbibnames=30, dashed=true, natbib]{biblatex} 
\usepackage{tikz} 
\usetikzlibrary{positioning,fit,calc} 
\usepackage{float} 
\usepackage{subcaption} 

\DeclareFieldFormat{eprint:mrnumber}{%
  \ifhyperref
    {\href{http://www.ams.org/mathscinet-getitem?mr=MR#1}{MR#1}}
    {MR#1}%
  \addspace
  }
\DeclareFieldFormat{eprint:zbmath}{%
  zbMATH\addcolon\space
    \ifhyperref
      {\href{https://zbmath.org/#1}{\nolinkurl{#1}}}
      {\nolinkurl{#1}}%
  \addspace
  }

\renewbibmacro{begentry}{\midsentence}
\makeatletter
\AtBeginDocument{\toggletrue{blx@useprefix}}
\makeatother

\tikzset{block/.style={draw, rectangle,
           rounded corners, thick, text width=3.5cm, minimum height=1.5cm, align=center},   
line/.style={-latex}   
}  

\newtheorem{theo}{Theorem}[section]
\newtheorem{lem}{Lemma}[section]

\newtheorem{prop}{Proposition}[section]

\newtheorem{ex}{Example}[section]
\newtheorem{defi}{Definition}[section]
\newtheorem{rk}{Remark}[section]

\newcommand{\ra}{\rightarrow}
\newcommand{\p}[1]{{\mathbb{P}} \left( \, #1 \, \right) }

\newcommand{\e}{{\mathbb{E}} }

\newcommand{\skk}[1]{\left\{ #1 \right\}}

\renewcommand{\leq}{\leqslant}
\renewcommand{\geq}{\geqslant}
\renewcommand{\le}{\leqslant}
\renewcommand{\ge}{\geqslant}
\def\RV{\mathcal{RV}}
\def\ERV{E\mathcal{RV}}
\def\MRV{\mathcal{MRV}}
\def\oF{\overline F}
\def\la{\leftarrow}

\def\M{\mathbb{M}}
\def\R{\mathbb{R}}

\def\bB{\boldsymbol B}

\def\bTheta{\boldsymbol \Theta}

\def\bX{\boldsymbol X}

\def\bzero{\boldsymbol 0}
\def\a{\alpha}

\def\conv{\stackrel{v}{\to}}
\newcommand{\1}{\mbox{\rm 1\hspace{-0.3em}I}}

\usepackage[strict]{changepage}

\usepackage{framed}

\definecolor{formalshade}{RGB}{202, 226, 232}
\definecolor{formal}{RGB}{0, 53, 63}
\definecolor{mine}{RGB}{255, 191, 102}
\definecolor{questionshade}{RGB}{255, 83, 69}
\definecolor{question}{RGB}{0, 53, 63}

\title{Multi-Normex Distributions \\for the Sum of Random Vectors. Rates of Convergence}

\author{\Large  Marie Kratz$\dagger$ and Evgeny Prokopenko$\dagger$ $\ddagger$ \\[1ex]
\small $\dagger$ ESSEC Business School Paris, CREAR, France \\ \small $\ddagger$ Novosibirsk State University, Novosibirsk, Russia }

\begin{document}

\date{July 19, 2021}

\maketitle

\begin{abstract}
\noindent We build a sharp approximation of the whole distribution of the sum of iid heavy-tailed random vectors, combining mean and extreme behaviors. It extends the so-called 'normex' approach from a univariate to a multivariate framework. We propose two possible multi-normex distributions, named $d$-Normex and MRV-Normex. Both rely on the Gaussian distribution for describing the mean behavior, via the CLT, while the difference between the two versions comes from using the exact distribution or the EV theorem for the maximum. The main theorems provide the rate of convergence for each version of the multi-normex distributions towards the distribution of the sum, assuming second order regular variation property for the norm of the parent random vector when considering the MRV-normex case. Numerical illustrations and comparisons are proposed with various dependence structures on the parent random vector, using QQ-plots based on geometrical quantiles. 
\\[1ex]
\noindent {\emph 2010 AMS classification}: 
60E05; 
60F05; 
60G70; 
60F15; 
62G30; 
41A25 
\\[1ex]

\end{abstract}

\tableofcontents

\newpage


\section{Introduction}\label{sec:intro}

\paragraph{\sc Motivation.}
Looking for the most accurate possible evaluation of the distribution of the sum of random variables, or vectors, or  processes, with unknown distributions, has always been a classical problem in the probabilistic and statistical literature,  with various answers depending on the given framework and on the specific application in view. On one hand, (uni- or multivariate) Central Limit Theorems (CLT) or Functional ones  prove, under finite variance for the sum components or/and additional conditions, the asymptotic Gaussian behavior of the sum with some rate of convergence, focusing on the 'body' of the distribution. When considering heavy-tailed marginal distributions, Generalized CLT  with the convergence to stable distributions, handle the case of infinite variance (see e.g. \citet{samo:1994,Petrov:1995}, and references therein), while, in the case of finite variance, an alternative way is to consider trimmed sums,  removing extremes from the sample, to improve the rate of convergence; see e.g. \citet{Mori:1984, HahnAl:1991} and references therein.
\\[.7ex] When interested in tail distributions, CLTs may give poor results, especially when considering heavy tails. That is why different approaches have been developed,  among which large deviation theorems (see e.g. \citet{Petrov1975,Borovkov2020} for light tails and  \citet{Mikosch:Nagaev:1998,Foss2013,Lehtomaa:2017} for heavy tails, and references therein),
extreme value theorems (EVT) focusing on the tail only (see e.g. \citet{embrechts:kluppelberg:mikosch:1997, dehaan:ferreira:2006, resnickbook:2007}), and hybrid distributions combining (asymptotic) distributions for both the main and extreme behaviors when considering independent random variables (see e.g. \citet{csorgoAl:1988,zaliapinAl:2005,kratz:2014,mueller:2019}; we use the name given in \citet{kratz:2014} for this type of hybrid distribution/method/approach, namely {\it Normex} distribution/method/approach.   
Recall briefly the idea of Normex (for 'Norm(al)-Ex(tremes)') method. It consists of rewriting the sum of random variables as the sum of their ordered statistics, and splitting it into two main parts, a trimmed sum removing the extremes, and the extremes. Using that the trimmed sum of the first $n-k-1$ ordered statistics is conditionnally independent of the $k$ largest order statistics, given the $(n-k)$-th order statistics, we can express the distribution of the sum, integrating w.r.t. to the $(n-k)$-th order statistics and using a CLT for the conditional trimmed sum, and an EVT one for the $k$ largest order statistics.
Note that a benefit of Normex approach is that it does not require any condition on the existence of moments, as the CLT applies on truncated random variables.
\\[.7ex] The following example (as developed in \citet{kratz:2014}), simulating identically distributed and independent (iid) Pareto($\a$) (such that $\overline F(x) = x^{-\a}$ for $ x>1$) random variables (rv), with $\a=2.3$ (finite variance, but no third moment), illustrates perfectly the adding value of combining main and extreme behaviors.
\begin{figure}[h]
\centering\includegraphics[width=165mm]{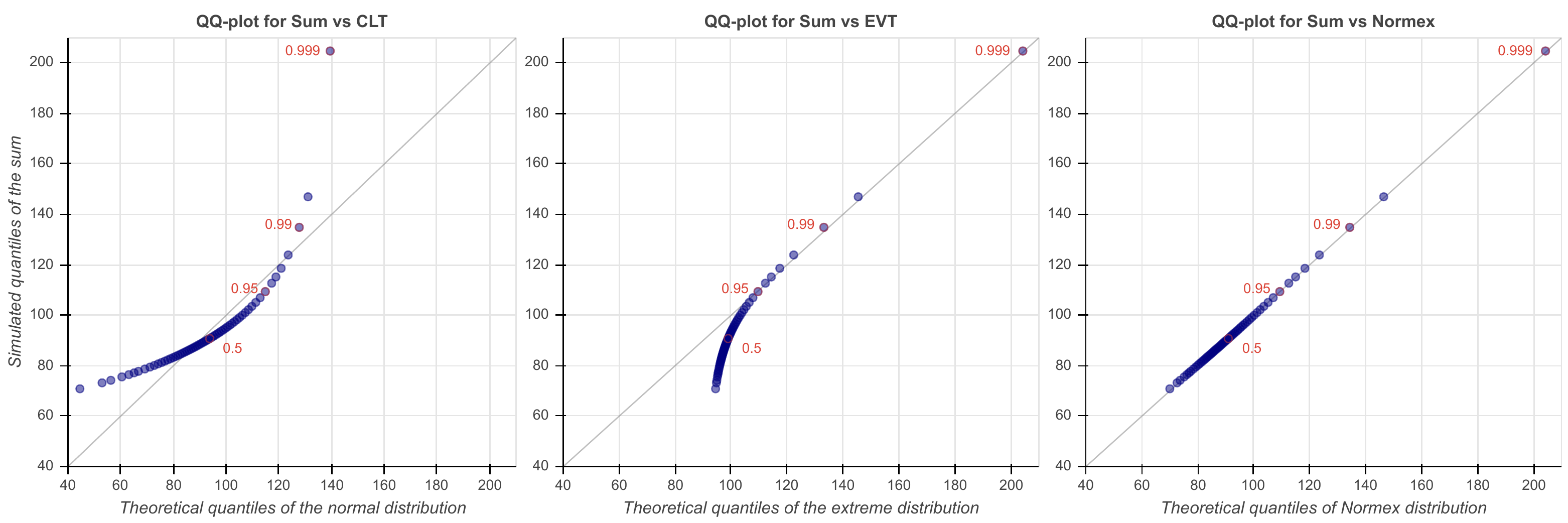}
\caption{\label{fig:dim1-compare3ways}\sf\small QQ-plots for the empirical distribution (sample size $= 10^7$) of the sum of 52 iid Pareto($\alpha = 2.3$) rv and three different approximations of the sum distribution using: CLT, EVT and Normex, respectively. The red circles and numbers denote the extreme quantiles and their levels.}
\end{figure}

Indeed, on the QQ-plots given in Figure~\ref{fig:dim1-compare3ways}, we can compare the fit of the distribution of the empirical sample with the following three distributions: the Gaussian one (CLT approach), the Fr\'echet one (EVT approach), the hybrid one ({\it i.e.} Normex, combining the CLT and the distribution of the maximum).
For the hybrid Normex distribution, we may consider either the exact distribution of the maximum, or its asymptotic approximation, the Fr\'echet distribution. Since both provide the same plot, we display it  in the third plot on the right.
Note that we choose a rather small number of components in the sum, $n=52$, also to illustrate the speed of convergence when using asymptotic theorems.
We observe that the CLT approach does not provide a sharp evaluation, even in the body of the distribution, due to this choice of $n$, which cannot compensate yet the fact that the distribution of the rv is asymmetric and skewed;  increasing $n$ will of course improve the fit in the body of the distribution. Given the fact that the Pareto distribution belongs to the Fr\'echet maximum domain of attraction, using the Fr\'echet distribution for the distribution of the sum of Pareto rv's gives a very sharp approximation in the tail (from the 93\% quantile), but not for the average behavior, as expected. Finally, a perfect match between empirical quantiles and Normex ones is observed for the whole distribution in the right plot, even for a small number of summands.

\paragraph{\sc Goal of the study.}
It is natural to extend the normex approach to a multivariate framework. With this goal of proposing a multi-normex method and distribution, we consider iid random vectors $\mathbf{X}_1, \dots, \mathbf{X}_n$, with parent random vector $\mathbf{X}$ having a heavy-tailed $d$-dimensional distribution $F_{\mathbf{X}}$ and density $f_{\mathbf{X}}$ (when existing).
Note that there are different ways to define multivariate extremes (see e.g. chap.8 in  \citet{beirlantAl:2004}). The chosen way in this paper is w.r.t. the norm $\|\cdot \|$ in $\mathbb{R}^d$, meaning that the ordered (w.r.t. the norm) vector of $(\mathbf{X}_1, \dots, \mathbf{X}_n)$, denoted by
$(\mathbf{X}_{(1)}, \dots, \mathbf{X}_{(n)})$, satisfies
\begin{equation}\label{ordered statistics}
  \|\mathbf{X}_{(1)}\| \leq  \|\mathbf{X}_{(2)}\|  \leq \dots \leq  \| \mathbf{X}_{(n)} \|.
\end{equation}
We propose two versions of multi-normex. The first one, named {\it $d$-Normex}, is a natural extension to any dimension $d$ of the univariate ($d=1$) normex method as developed in \citet{kratz:2014}: We approximate the distribution of the trimmed sum via the CLT and consider the distribution of the maximum $\mathbf{X}_{(n)}$. This latter distribution is approximated via the Extreme Value (EV) theorem in the second multi-normex version, named {\it MRV-Normex}.

Aiming at proving the benefit of using a multi-normex distribution for a better fit of the whole (unknown) distribution F,  
assuming $F$ heavy-tailed in the sense of $\|\mathbf{X}\| \in \RV_{-\alpha}$, {\it i.e.} regularly varying rv with $\alpha >0$ (which definition is recalled in Appendix~\ref{A:secRV}),
we focus analytically on the case $\alpha\in(2;3]$ (when $\|\mathbf{X}\|$ has a finite second moment, but no third moment), 
to compare the rates of convergence  when using the CLT and the multi-normex approach, respectively.
Note our focus on heavy tailed distributions (i.e. distributions belonging to the max domain of attraction of Fr\'echet), where the impact of using Normex distribution will be much stronger than in the light tail case (because of the one big jump principle), in particular for risk analysis and management.
 We prove that the normex approach leads, as expected, to a better speed of convergence for evaluating the distribution of the sum than the CLT does, for such type of heavy-tailed distributions.
 When varying the fatness of the tail measured by $\a>0$, we draw this comparison numerically, using geometrical multivariate quantiles (see e.g. \citet{chaudhuri:1996,DharAl:2014} or Appendix~\ref{A:secGeomQuant}).

\paragraph{\sc Structure of the paper.}~
In Section~\ref{sec:frame}, besides general notations, we recall the normex approach and the generalized Berry-Esseen inequality. Then we give two specific results on conditional distributions of order statistics, which will be needed for the construction of multi-normex distributions. The two next sections develop the two multi-normex versions, $d$-Normex in Section~\ref{sec:dNormex} and MRV-Normex in Section~\ref{sec:mrv}. These sections have the same structure: we first define the multi-normex distribution, then we study analytically its rate of convergence, before ending with some examples. In Section~\ref{sec:QQplots}, we consider those examples to study numerically the two versions of multi-normex distribution, comparing them with the empirical distribution of the sum (obtained via simulation) as well as, if relevant, with the Gaussian approximation when applying the CLT. Geometrical multivariate quantiles are computed to this aim and represented on QQ-plots.  Section~\ref{sec:concl} concludes. The appendix is divided in 3 parts. Appendix~\ref{A:secRV} recalls the various definitions of regular variation and its related concepts, as well as two properties of second order regular variation that are needed for the proofs of the analytical results. The proofs of all analytical results are developed in Appendix~\ref{A:proof}. An overview of geometrical multivariate quantiles is given in Appendix~\ref{A:secGeomQuant}.


\section{Framework}\label{sec:frame}

\subsection{Context and notations}\label{ssec:frame-intro}

\paragraph{\sc Normex idea -}

The Normex method clearly adapts to a multivariate framework. Using this approach, we split the maximum from the rest of the sum: 
\begin{equation}\label{eq:ideaNormex}
\mathbf{S}_n= \sum_{k=1}^{n-1} \mathbf{X}_{(k)}+\mathbf{X}_{(n)},
\end{equation}
taking into account the 'principle of one big jump', namely that the asymptotic tail behavior of the sum of heavy-tailed random vectors is led by that of the maximum. Indeed, this principle extends to dimension $d>1$ by borrowing the multivariate subexponential distribution definition of \citet{Gena:Sun:2016} (see Section 4 therein):
$$
 \p{\mathbf{S}_n>t\,\mathbf{A}} \;  \underset{t\to\infty} {\sim} \; \p{\mathbf{X}_{(n)}>t \,\mathbf{A}}, 
$$
where $\mathbf{A} \subset \R^d$ is open, increasing, such that $\mathbf{A}^c$ convex and $\mathbf{0}\notin \bar{\mathbf{A}}$.
Note also earlier works on that notion by \citet{cline:resnick:1992} and \citet{omey:2006}.
\\[0.7ex]
Combining the CLT for the trimmed sum given the maximum, 
and the distribution of the maximum or its asymptotic distribution, leads to a multivariate version of Normex distribution. We name this multi-normex distribution as $d$-Normex, when using the distribution of the maximum, and as MRV-Normex when considering its asymptotic distribution (when rescaled). 

\paragraph{\sc General notations -}
Before defining explicitly both versions of multi-normex distribution, let us introduce some general notations.
\\
 Let $f_{(i)}$ denote the $d$-dimensional density,  when existing, of the ordered vectors $\mathbf{X}_{(i)}$, for $i=1,\cdots,n$.
 The cumulative distribution function (cdf) of the norm $\|\mathbf{X}\|$ is denoted as $ F_{\|\mathbf{X} \|}(\cdot),$  which we assume to be absolutely continuous, and its probability density function (pdf) as $f_{\|\mathbf{X} \|}(\cdot) $. Note that we assume the absolutely continuous of $F_{\|\mathbf{X}\|}$ throughout the paper. Nevertheless, some of our results will be stated under the slightly stronger condition $F_{\mathbf{X}}$ absolutely continuous.

As we will work with truncated multidimensional distributions or vectors, let us introduce the following notions.

For any $y > 0$, we define the truncated (via the norm) multidimensional distribution $F_{\mathbf{X} \,|\, \| \mathbf{X}\|}(\cdot\,|\,y)$ of $\mathbf{X}$ on $\R^d$ as
\begin{equation}\label{eq:def_trunc_distr}
  F_{\mathbf{X}\,|\, \| \mathbf{X}\|}(\mathbf{B}\,|\, y) := \p{\mathbf{X} \in \mathbf{B} \,|\, \| \mathbf{X}\| \leq y}
\end{equation}
for any event $\bB$ of the Borel sigma-field  ${\cal B}(\mathbb{R}^d)$.
We denote by $f_{\mathbf{X} \, |\,\| \mathbf{X}\|}(\cdot) $ its pdf, when existing:
\begin{equation}\label{eq:def_trunc_pdf}
f_{\mathbf{X}\,|\,\| \mathbf{X}\|}(\mathbf{x}\,|\,y)= \frac{f(\mathbf{x})\,\1_{\left(\| \mathbf{x} \| \leq y\right)}}{ F_{\|\mathbf{X} \|}(y)}.
\end{equation}

Let $\overset{\circ}{\mathbf{X}}_y  \in \mathbb{R}^d$ denote the random vector with distribution $F_{\overset{\circ}{\mathbf{X}}_y }$ on  the $(d-1)$-sphere \\$\mathcal{S}_y = \skk{\mathbf{x} \in \mathbb{R}^d\,:\, \|\mathbf{x}\| = y}$ for $y > 0$, defined by
\begin{equation}\label{eq:def-Zcdf}
 F_{\overset{\circ}{\mathbf{X}}_y }(\mathbf{B}_y)= \p{ \overset{\circ}{\mathbf{X}}_y \in \mathbf{B}_y  } :=  \p{ \mathbf{X} \in \mathbf{B}_y \,|\, \|\mathbf{X}\| = y } \text{ for } \mathbf{B}_y \subseteq \mathcal{S}_y.
\end{equation}
With this definition, for any $\mathbf{A} \subseteq \mathbb{R}^d$, we can write
$$
\p{\mathbf{X} \in \mathbf{A}}  = \e\left[\e[\1_{\mathbf{X} \in \mathbf{A}} \,|\, \|\mathbf{X}\|] \right] 
= \int_0^{\infty} \e[\1_{\mathbf{X} \in \mathbf{A}}\,|\, \|\mathbf{X}\|=y]\, f_{\|\mathbf{X}\|}(y)dy,
$$
thus
\begin{equation}\label{eq: Z}
  \p{\mathbf{X} \in \mathbf{A}} = \int\limits_{\mathbb{R}_+} \p{\overset{\circ}{\mathbf{X}}_y \in \mathbf{A}_y} f_{\|\mathbf{X}\|}(y) \mathrm{d}y,\quad  \mbox{where}\quad \mathbf{A}_y:= \mathbf{A}\cap \mathcal{S}_y. 
\end{equation}

\paragraph{\sc Generalized Berry-Esseen inequality -}

As we are going to compare the rates of convergence when using, respectively, the CLT and the Normex method, let us recall the rate of convergence in the CLT provided by the Generalized Berry-Esseen inequality (see e.g. Corollary 18.3 of \citet{BhattacharyaRao:2010}), when assuming a 'moderate' heavy tail.
\begin{prop}[Generalized Berry-Esseen inequality]\label{thBEine}~

 Let $\mathbf{X}_{1}, \ldots, \mathbf{X}_{n}$ be i.i.d. centered random vectors with parent random vector $\mathbf{X}$ with values in $\R^{d}$, with positive-definite covariance matrix $\Sigma$.
If
\begin{equation*}
\e \left\|\mathbf{X}\right\|^{\alpha}<\infty, \quad \text{for some}\;\;\alpha\in [2,3],
\end{equation*}
then
\begin{equation}\label{eq:BEineq}
  \sup _{\mathbf{B} \in \mathcal{C}}\left|\p{\mathbf{S}_n \in \mathbf{B}}-\Phi_{\mathbf{0}, \Sigma}(\mathbf{B})\right| \leq  c \, \e \left(\left\|\Sigma^{-1/2} \mathbf{X}\right\|^{\alpha}\right) \,n^{-(\alpha-2) / 2},
\end{equation}
where $\mathcal{C}$ is the class of all Borel-measurable convex subsets of $\R^{d}$,  $c$ is a positive universal constant, and $\Phi_{\mathbf{0}, \Sigma}$ is the cdf of the centered normal multivariate distribution with covariance matrix $\Sigma$.
\end{prop}
Note that the (non-generalized) Berry-Esseen inequality, which holds for any  $\a \ge 3$,  corresponds to \eqref{eq:BEineq} when taking $\alpha = 3$.

\subsection{Preliminary results on order statistics}\label{ssec:orderStat}

We present two simple but elegant results on order statistics, needed for the proofs of the theorems, which are of interest in themselves, completing the vast literature on order statistics.
\begin{lem}\label{lem1}
The distribution of the first $n-1$ order statistics $ \mathbf{X}_{(1)}, \dots, \mathbf{X}_{(n - 1)}$, conditionally on the event $\| \mathbf{X}_{(n)}\| = y$, is the distribution of the $n-1$ ordered statistics from the truncated distribution $F_{\mathbf{X}\,|\,\| \mathbf{X}\|}(\cdot\,|\, y)$:
\begin{equation*}
  \mathcal{L}\left(  \mathbf{X}_{(1)}, \dots, \mathbf{X}_{(n -1)} \,\big|\, \| \mathbf{X}_{(n)}\| = y \right) =  \mathcal{L}\left(  \mathbf{Y}_{(1)}, \dots, \mathbf{Y}_{(n - 1)} \right),
\end{equation*}
where $\mathbf{Y}_1, \dots, \mathbf{Y}_{n- 1}$ are i.i.d. random vectors with multidimensional  distribution $F_{\mathbf{X}\,|\,\| \mathbf{X}\|}(\cdot\,|\, y)$ defined in \eqref{eq:def_trunc_distr}.
\end{lem}

Similarly, we  have the following result.
\begin{lem}\label{lem1-bis}
The distribution of the order statistics $ \mathbf{X}_{(1)}, \dots, \mathbf{X}_{(n - 1)}, \mathbf{X}_{(n )}$, conditionally on the event $\| \mathbf{X}_{(n)}\| = y$, is the distribution of the $n-1$ ordered statistics from the truncated distribution $F_{\mathbf{X}\,|\,\| \mathbf{X}\|}(\cdot\,|\, y)$, and of an independent random vector  $\overset{\circ}{\mathbf{X}}_y$ defined on the $(d-1)$-sphere $\mathcal{S}_y$, for $y > 0$:
\begin{equation*}
  \mathcal{L}\left(  \mathbf{X}_{(1)}, \dots, \mathbf{X}_{(n -1)}, \mathbf{X}_{(n)} \,\big|\, \| \mathbf{X}_{(n)}\| = y \right) =  \mathcal{L}\left(  \mathbf{Y}_{(1)}, \dots, \mathbf{Y}_{(n - 1)}\right) \times  \mathcal{L} \left(  \overset{\circ}{\mathbf{X}}_y \right),
\end{equation*}
where $\mathbf{Y}_1, \dots, \mathbf{Y}_{n- 1}$ are i.i.d. random vectors with multidimensional truncated  distribution $F_{\mathbf{X}\,|\,\| \mathbf{X}\|}(\cdot\,|\, y)$ defined in \eqref{eq:def_trunc_distr}, and the random vector $\overset{\circ}{\mathbf{X}}_y$ has the distribution $F_{\overset{\circ}{\mathbf{X}}_y}(\cdot)$ defined in \eqref{eq:def-Zcdf}.
\end{lem}

These lemmas are proved in Appendix~\ref{A:sec2}.


\section{A first multi-normex version: \texorpdfstring{$d$}{Lg}-Normex}\label{sec:dNormex}

We start building a first multi-normex version, using the Normex approach \eqref{eq:ideaNormex}, then approximating the distribution of the trimmed sum via the CLT and keeping the distribution of the maximum $\mathbf{X}_{(n)}$. It is a natural extension to any dimension $d$ of the univariate ($d=1$) Normex distribution as developed in \citet{kratz:2014}. Note that, when turning to data, the distribution of $\mathbf{X}_{(n)}$ may be approximated e.g. via simulations or, as will be done in the MRV-Normex, using another asymptotic theorem, the EV one.

\subsection{Definition}\label{ssec:def_dNormex}

\begin{defi}\label{def:Normex}
The so-called $\mathbf{d}${\bf -Normex} distribution function  is defined, for $\mathbf{B} \subset  \R^d$, as:
\begin{equation*}
   G_n(\mathbf{B})= \p{\mathbf{Z}+ \mathbf{X}_{(n)} \in \mathbf{B}}   
\end{equation*}
where $\mathbf{Z}$ is, conditionally on the event $(\|\mathbf{X}_{(n)}\| = y)$, a Gaussian random vector with mean $(n-1)\mathbf{\mu}(y)$ and covariance matrix $ (n-1) \Sigma(y)$, the functions $\mathbf{\mu}(\cdot)$ and $\Sigma(\cdot)$ are, respectively, the mean vector and covariance matrix of the truncated distribution $F_{\mathbf{X}\,|\, \| \mathbf{X}\|}$ defined in \eqref{eq:def_trunc_distr}.
\end{defi}
Another way to formulate Definition~\ref{def:Normex} is the following:
\begin{equation}\label{eq:normex-cdf}
 G_n(\mathbf{B})= \int_{\mathbb{R}^d} f_{(n)}(\mathbf{x}) \;\Phi_{(n-1)\mathbf{\mu}(\|\mathbf{x}\|),\, (n-1) \Sigma(\|\mathbf{x}\|) }\left(\mathbf{B} - \mathbf{x}\right)\,\mathrm{d}\mathbf{x},
\end{equation}
where $\Phi_{m, \Gamma}$ denotes the cdf of the Gaussian vector with mean $m$ and covariance matrix $\Gamma$.
  
\subsection{Rate of convergence}\label{ssec:MResult}

Let us turn to the evaluation of the Normex distribution for approximating the distribution of the sum of iid random vectors, studying its rate of convergence. 
We do it analytically. Although the multi-normex method  works for any $\a \ge 2$,
we state the result when assuming the same condition on moments as in the generalized Berry-Esseen inequality, namely $\a\in(2,3]$, to be able to compare the results and show explicitly the benefit of using Normex approximation. 
 Then we show numerically the general good fit of $d$-Normex  in an example (see Section~\ref{ssec:exs-dNormex}).

%
The analytical result given in Theorem~\ref{th:Result-max} shows that applying Normex method rather than the multivariate CLT improves, as expected, the accuracy of the evaluation of the (tail) distribution of the sum of heavy tailed vectors, with a better rate of convergence than the one of the CLT whenever the shape parameter $\alpha\in(2,3)$.

\begin{theo}\label{th:Result-max}
Let $\mathbf{X}_{1}, \ldots, \mathbf{X}_{n}$ be i.i.d. random vectors with parent random vector  $\mathbf{X}$ with values in $\R^{d}$ such that:
\begin{itemize}
\item[$(C1)$] For all $y >0$, the truncated (w.r.t. the norm) distribution $F_{\mathbf{X}\,|\, \| \mathbf{X}\|}(\cdot\,|\, y)$ defined in \eqref{eq:def_trunc_distr} is nondegenerate (i.e. for all $y > 0$ there is no  hyperplane $\mathcal{H} \subset \mathbb{R}^d$ such that $\p{\mathbf{X} \in \mathcal{H}\,|\, \|\mathbf{X}\| \leq y} = 1$).
\item[$(C2)$] The distribution of the rv $\|\mathbf{X}\|$ is absolutely continuous and regularly varying at infinity: $\|\mathbf{X}\| \in \RV_{-\alpha}$, with $\alpha >0$.
\end{itemize}
Then, for any $\alpha \in (2,3]$, there exists a slowly varying function $L(\cdot)$ such that
\begin{equation*}
   \sup _{\mathbf{B} \in \mathcal{C}}  | \p{\mathbf{S_n} \in \mathbf{B}} - G_n(\mathbf{B}) | \leq L(n) \, n^{- \frac{1}{2} + \frac{3-\alpha}{\alpha}},
\end{equation*}
where $\mathcal{C}$ is the class of all Borel-measurable convex subsets of $\R^{\mathbf{d}}.$
\end{theo}

Let us briefly indicate how to reach the upper bound of this main result; for more details, see the proof
 developed in Appendix~\ref{A:sec3}.
First, we use the law of total probability conditioning by $\mathbf{X}_{(n)}$. Second, we apply the (non-generalized) Berry-Esseen inequality for the truncated r.v. $\{\mathbf{Y}_i\}_{i \leq n-1}$ given the event $(\|\mathbf{X}_{(n)}\| = y)$. The right-hand side of the inequality is of the order of $\frac{1}{\sqrt{n}} \e \|\mathbf{Y}\|^{3}$, which is equivalent to  $\frac{1}{\sqrt{n}} \e \|\mathbf{X}_{(n)}\|^{3-\alpha}$ whenever $\alpha\leq 3$. Finally, to derive the upper bound of the main result, we use that $\|\mathbf{X}_{(n)}\|$ is of the order $n^{1/\alpha}$ under $(C2)$.

\begin{rk}\label{rk:main-theo}~\\[-4ex]
\begin{itemize}

\item[(i)] 
Note that we consider the case $\a\in (2,3]$ since it is the condition under which the Generalized Berry-Esseen holds. For $\a>3$,  the bound given in Theorem~\ref{th:Result-max} is the same as that of Berry-Esseen inequality, making the analytical comparison useless. Indeed, in such a case, the bound $\frac{1}{\sqrt{n}} \e \|\mathbf{Y}\|^{3}$ reduces simply to the order $\frac{1}{\sqrt{n}}$ (see \eqref{eq:proof_th3} in the proof), giving back the same rate as for the CLT. It means to look for an alternative way if we want to study analytically the Normex rate of convergence. We might use Edgeworth expansions, but 
it evolves too heavy computations (as we could experience for rv (the case $d=1$), conditioning on $X_{(n)}$). 
This is why we show numerically the benefit of using Normex distribution, as illustrated in Section~\ref{sec:QQplots}.
\item[(ii)] The rate of convergence given in Theorem~\ref{th:Result-max} is better than the one provided in the generalized Berry-Esseen inequality (Proposition~\ref{thBEine}),  whenever $\alpha \in (2,3)$ (and whatever $n$), as
\begin{equation*}
  \frac{\alpha - 2}{2} < \frac{1}{2} - \frac{3-\alpha}{\alpha}.
\end{equation*}
Note also that in the case $\alpha = 3$ and $\e\|\mathbf{X}\|^3 = \infty$, the inequality in Theorem~\ref{th:Result-max} is  slightly sharper than the inequality that can be obtained by the Berry-Esseen theorem (replacing  $n^{-1/2 + \varepsilon}$, $\varepsilon > 0$, by $L(n)n^{-1/2}$).
\item[(iii)]  One can apply the Normex method with any norm on $\R^d$, for instance the $L^1$ norm defined, for $\mathbf{x}=(x_1,\cdots,x_d) \in\R^d$, by $\|\mathbf{x}\|_1:= \sum_{i=1}^d |x_i|$.
\item[(iv)] Considering the $L^1$ norm, for positive random variables, Condition $(C2)$ translates into the assumption
$$
(C2^*) \qquad S_d:= \sum_{i=1}^d X^{(i)} \in \RV_{-\alpha},
$$
where $X^{(i)}$, for $i=1,\cdots, d$, denote the components of $\bX$.

We may want to relate this $\RV$ property on the sum, with conditions on the random vector itself. This topic has already been investigated in the literature; see e.g. \citet{basrak:davis:mikosch:2002b,barbe_etal:2006,embrechtsMainik:2013,cuberosAl:2015}.
For instance, assuming  $\bX$ multivariate regularly varying,  $\bX\in \MRV_{-\alpha}(b,\nu)$ (as defined in Definition~\ref{def:mrv}, Appendix~\ref{A:secRV}), implies
that the sum $S_d\in \RV_{-\alpha}(b)$. We will come back on the MRV notion in Section~\ref{sec:mrv}.
\end{itemize}
\end{rk}


\subsection{Example of the Multivariate Pareto-Lomax distribution}\label{ssec:exs-dNormex}

Let us consider a $d-$dimensional random vector $\mathbf{X} = (X^{(1)}, \dots, X^{(d)})$ having a multivariate Pareto-Lomax$(\a)$ distribution, with $\a>0$, {\it i.e.}  
with survival distribution function defined, for any non-negative real numbers $x_1, \dots, x_d$, by
\begin{equation*}
 \overline{F}_{\mathbf{X}}(x_1, \dots, x_d) :=  \p{ X^{(i)} > x_i,\, i=1,\cdots,d} =\left (1 + \sum_{i=1}^d x_i\right)^{-\alpha}.
\end{equation*}

First we compute the statistical characteristics, as (truncated) moments, needed for applying $d$-Normex. We consider the case $\alpha \in (2,3]$, as in Theorem~\ref{th:Result-max} (even if the method works for any $\alpha > 0$).

It is straightforward to compute the pdf, marginal distributions, expectation and covariance matrix of a multivariate Pareto-Lomax vector, namely:
\begin{equation}\label{density}
  f_{\mathbf{X}}(x_1, \dots, x_d) =  \frac{\alpha(\alpha+1)\dots(\alpha + d -1)}{(1 +  \sum_{i=1}^d x_i)^{\alpha + d}}, \quad  x_1, \dots, x_d \geq 0;
\end{equation}
\begin{equation}\label{def:one_dim_par}
  \p{X_j > x} = (1 + x)^{-\alpha}, \ \ x \geq 0, \, j\in\skk{1,\dots,d};
\end{equation}
\begin{equation*}
  \e [X_j ]= \frac{1}{\alpha -1},  \quad \mathrm{Var} (X_j)= \frac{\alpha}{(\alpha -1)^2(\alpha-2)}, \ \  j\in\skk{1,\dots,d};
\end{equation*}
\begin{equation*}
  \e[ X_j X_i ]= \frac{1}{(\alpha -1)(\alpha-2)},\ \ \mathrm{ Cov}(X_i, X_j) = \frac{1}{(\alpha -1)^2(\alpha-2)}, \ \ \  j\neq i\in\skk{1,\dots,d}.
\end{equation*}
Given the expression of the survival distribution of $\mathbf{X}$, and given that we can apply $d$-Normex for any norm on $\mathbb{R}^d$, we choose the example of the $L_1$ norm $\displaystyle \|\mathbf{x}\|_1:= \sum_{i=1}^d x_i$ to simplify the computations. We also take the example of $d=3$ for illustration.

We can express the cdf of the rv $\|\mathbf{X}\|$, for $y>0$, as
\begin{equation*}
  F_{\|\mathbf{X}\|}(y) = \p{  \|\mathbf{X}\| \leq y } = 1 - (1+y)^{-\alpha} - \alpha \,y (1 + y)^{-(\alpha+1)}   -\frac{\alpha(\alpha+1)}{2}\, y^2(1+y)^{-(\alpha + 2)}.    
\end{equation*}
Note that those expressions may also be obtained when applying the more general results \eqref{def:exampledensity} and \eqref{eq:cdf-norm} developed for any norm on $\mathbb{R}^d$ and for any $d$; see Appendix~\ref{A:sec4disc}.

Now, let us compute the moments for the $d$-dimensional truncated Pareto-Lomax  random vector, denoted by $\mathbf{Y}$, having cdf $F_{\mathbf{X}\,|\, \| \mathbf{X}\|}(\cdot\,|\,y)$ (see Definition~\ref{eq:def_trunc_distr}), expectation $\mathbf{\mu}(y)$ and covariance matrix $\Sigma(y)$.
We have, for any $j \in \skk{1,2}$, if $\alpha \neq 1$ (which is our case),
\begin{equation*}
\begin{split}
  \mathbf{\mu}_j(y) & := \e[Y_j] = \e \left[X_j\,\big|\, \|\mathbf{X}\| \leq y \right] \\
  & = \frac{1}{F_{\|\mathbf{X}\|}(y)} \dfrac{\left(y+1\right)^{-\alpha-2}\big(\left(y+1\right)^{\alpha+2}-\left(\alpha+2\right)y\left(1+y\left(\alpha+1\right)\left(\alpha y+3\right)/6\right)-1\big)}{\left(\alpha-1\right)},
\end{split}
\end{equation*}
and, if $\alpha \neq 1,2$ (also our case), 
\begin{equation*}
\begin{split}
  &\e \left[Y_1 Y_2\right] = \e \left[ X_1 X_2\,\big|\, \|\mathbf{X}\| \leq y  \right] 
 =  \frac{1}{F_{\|\mathbf{X}\|} (y)}  
\dfrac{1}{\left(\alpha-2\right)\left(\alpha-1\right)}\Bigl(1 - \left(y+1\right)^{-\alpha-2}\times\Bigr.\\
& \qquad \qquad\qquad\Bigl. \big(\left(\alpha+2\right)y\left(\left(\alpha+1\right)y\left(\a y \left(\left(\alpha-1\right)y/24+1/6\right)+1/2\right)+1\right)-1\big)\Bigr),
\end{split}
\end{equation*}
\begin{equation*}
\begin{split}
  &\e \left[Y_1^2\right] = \e \left[ X_1^2\,\big|\, \|\mathbf{X}\| \leq y  \right]  =\frac{1}{F_{\|\mathbf{X}\|} (y)}  \dfrac{2}{\left(\alpha-2\right)\left(\alpha-1\right)}\Bigg( 1 - \left(y+1\right)^{-\alpha-2}\times\\ 
 &  \qquad \qquad\qquad
{\bigg[\left(\alpha+2\right)y\Big(y\left(\alpha+1\right)\big(\alpha y\left(y\left(\alpha-1\right)/24+1/6\right)+1/2\big)+1\Big)+1\bigg]\Bigg)},
\end{split}
\end{equation*}
from which can be deduced the covariance matrix $\Sigma(y)=(\Sigma_{ij}(y))_{i,j}$.  

Therefore, the Gaussian cdf  $\displaystyle \Phi_{(n-1)\mathbf{\mu}(\|\mathbf{x}\|),\, (n-1) \Sigma(\|\mathbf{x}\|) } $ introduced in Definition~\ref{def:Normex} is explicitly determined, and so is the $d$-Normex distribution $G_n$ defined in \eqref{eq:normex-cdf}. \\

To illustrate the benefit of the Normex method (see Section \ref{sec:QQplots}), whatever the value $\alpha>0$, we draw QQ-plots (see Figure~\ref{fig:qqplot-23} to Figure~\ref{fig:qqplot 3,5}). For that, we simulate a sample from the multivariate Pareto-Lomax($\alpha$) distribution \eqref{density}. We proceed by induction:  the first coordinate $X^{(1)}$ has ($1$-dimensional) Pareto-Lomax distribution \eqref{def:one_dim_par} with parameter $\alpha >0;$ for the rest of the components, we can easily derive their conditional distribution: 
\begin{equation*}
   \mathcal{L} \left( X^{(k+1)} \,\big|\, \left(X^{(1)}, \cdots, X^{(k)}\right) = \left(x^{(1)}, \cdots, x^{(k)}\right) \right) =  \left(1+x^{(1)}+ \cdots+ x^{(k)}\right) \widehat{X},\quad k\ge 1,
\end{equation*}
where $\widehat{X}$ has the Pareto-Lomax distribution \eqref{def:one_dim_par} with parameter $\alpha + k$ and is independent of $X^{(1)}, \cdots, X^{(k)}.$


\section{MRV-Normex}\label{sec:mrv}

Here we investigate a more universal version of multi-normex, named MRV-Normex, using an asymptotic theorem for the maximum, namely the Extreme Value (EV) one. Given our focus on the sum of  iid heavy-tailed (w.r.t. the norm) random vectors, 
we consider the standard extreme value theory (EVT) framework of multivariate regularly varying (MRV), a natural extension of the regular variation in a multivariate framework. In fact, to obtain the rate of convergence of this multi-normex approximation, we assume slightly stronger assumption than MRV, asking for a {\it uniform} asymptotic independence of the polar coordinates of the random vector, as made explicit in Condition $(M_{\Theta})$ of Theorem~\ref{th:MRV result}  below.

In order to obtain the rate of convergence for the MRV-Normex approximation of the sum, we first need to discuss the rate of convergence in  EVT to control the difference between the norm of the maximum $\|\mathbf{X}_{(n)}\|$ and the limit Fr\'echet distribution. This is the object of the next subsection.

%
\subsection{Rate of convergence in the EV theorem. Discussion of its assumptions}

After recalling the Extreme Value Theorem, we discuss its rate of convergence, depending on the assumptions. 

\paragraph{Extreme Value (EV) Theorem -} 
{\it Let $\skk{X_n, n \geq 1}$ be i.i.d. random variables with c.d.f. $F_X$. Denote the maximum of this sequence as $\displaystyle \max_{1\le i\le n} X_i$. Assume that $F_X$ is in the maximum domain of attraction (MDA) of an extreme-value distribution $G_{\gamma}$, denoted by $F_{X}\in MDA(G_{\gamma})$ with $\gamma\in\R$, {\it i.e.}
 there exist normalizing constants $a_n > 0$ and $b_n \in \mathbb{R}$ such that
\begin{equation}\label{eq:EVT}
\p{\max_{1\le i\le n} X_i\leq a_{n} x+b_{n}} =F_X^n(a_{n} x+b_{n}) \,\underset{n\to\infty} {\longrightarrow} \, G_{\gamma}(x).
\end{equation}
}
Let us introduce the real function $g$ defined on $\R^+$ by :
\begin{equation}\label{eq:def-g}
  g:=\left(\frac{1}{-\log F_X}\right)^{\leftarrow}
\end{equation}
( $^{\leftarrow}$ denoting the left-continuous inverse function).

It is straightforward to show that the convergence \eqref{eq:EVT} is equivalent to  
\begin{equation}\label{eq:EVT2a}
  \lim _{t \rightarrow \infty} \frac{g(t x)-g(t)}{a(t)}=\frac{x^{\gamma}-1}{\gamma}, \quad \forall x > 0, 
\end{equation}
for some $\gamma\in\R$ and auxiliary positive function $a$ defined on $\R^+$,
{\it i.e.} $g$ is of extended regular variation (see Appendix~\ref{A:secRV}), denoted by $g\in \ERV_\gamma (a)$:
$$
 F_{X}\in MDA(G_{\gamma}) \quad \Leftrightarrow \quad g\in \ERV_\gamma (a).
$$
Note that this equivalence holds true when replacing $g$ by $\displaystyle U= \left(\frac{1}{1- F_X}\right)^{\leftarrow}$, for some auxiliary positive function $a$ (see e.g. \citet{dehaan:ferreira:2006}, Theorem 1.1.6).

If $F_X$ is differentiable, then the auxiliary function $a$ in \eqref{eq:EVT2a} can be chosen as $a(t):=t g^{\prime}(t)$, so that
\begin{equation}\label{eq:EVT2}
  \lim _{t \rightarrow \infty} \frac{g(t x)-g(t)}{t g^{\prime}(t)}=\frac{x^{\gamma}-1}{\gamma}, \quad \forall x > 0, \quad\text{with}\; \gamma\in\R.
\end{equation}

\paragraph{On the assumptions for the rate of convergence in the EV theorem -} 

To describe the  rate of convergence in the EV theorem, we first refer to two studies developed under slightly different assumptions, namely \citet{falk:marohn:1993} with a direct condition on the derivative of the distribution $F_X$ (see \eqref{Falks_cond}), and \citet{Haan:Resnick:1996} assuming a second-order von Mises condition on $g$ defined in \eqref{eq:def-g}. Note that the results obtained in both studies hold for $F_X$ belonging to any MDA. 
Then, focusing on the case of $F_X\in$MDA(Fr\'echet), we look for a condition on $F$ involving RV properties to replace the second-order von Mises condition on $g$ and to retrieve the exact rate of convergence described in \citet{Haan:Resnick:1996}. This is presented in Proposition~\ref{cvRate-Frechet}.
\\[1ex]
$\blacktriangleright$ {\sc De Haan and Resnick's result.}
In \citet{Haan:Resnick:1996}, the authors assume that the function $g$ defined in \eqref{eq:def-g}  satisfies the second-order von Mises condition, which we recall here.
\begin{defi}\label{def:vonMis}(second-order von Mises condition) A twice differentiable function $g:(0, \infty) \ra \mathbb{R}$ satisfies the {\it second-order von Mises condition } (shortly, $g\in2\!-\!von\,Mises(\gamma,-\rho)$) for $\gamma\in \mathbb{R}$ and $\rho>0$, if $g'$ is eventually positive and the function
  $$A(t):=\frac{t \, g^{\prime \prime}(t)}{g^{\prime}(t)}-\gamma+1$$
  has constant sign near infinity and is such that 
  $$
  A(t) \ra 0,   \text{ as } t \ra \infty, \quad \text{and} \quad |A| \in RV_{-\rho}.
  $$
An equivalent definition (see \citet{Haan:Resnick:1996}, Theorem 2.1) is that $g'$ has the representation:
\begin{equation}\label{eq: Karamata-deriv}
g'(t)=k\, t^{\gamma-1}\exp\left\{\int_1^t \frac{A(u)}{u} du\right\}, \quad\text{with}\;k\ne 0.
\end{equation}
\end{defi}
In particular, it implies that (see \citet{Haan:Resnick:1996}, Theorem 2.1)
$$ g'\in 2\RV_{\gamma-1, -\rho},$$
where the $2\RV$-class can be defined as \eqref{Proko_cond} below (see also Lemma~\ref{lem:representation_2RV} in Appendix \ref{A:secRV}).
De Haan and Resnick use the 2-von Mises condition representation \eqref{eq: Karamata-deriv}  to derive Potter bounds  for $g'$ and obtain the rate of convergence for \eqref{eq:EVT2}:
\begin{equation*}
  \frac{[g(t x)-g(t)] / (t g^{\prime}(t)) -\left(x^{\gamma}-1\right) / \gamma}{A(t)} \rightarrow \int_{1}^{x} u^{\gamma-1}\left(\frac{u^{\rho}-1}{\rho}\right) d u, \quad\text{as}\; t\to\infty,
\end{equation*}
and the following rate of convergence in \eqref{eq:EVT} for the total variation metric, given in Proposition~\ref{thResnick}.
\begin{prop}[Theorem 4.1 in \citet{Haan:Resnick:1996}]\label{thResnick}
~\\
Let $g$ be defined in \eqref{eq:def-g} such that  $g\in2\!-\!von\,Mises(\gamma,\text{-}\rho)$ with $\gamma \in \mathbb{R}$ and $\rho>0$. Then there exists a constant $C>0$ (that is defined explicitly) such that
\begin{equation}\label{eq:thResnick}
 \lim_{n \ra \infty} \frac{\sup\limits _{A \in B(\mathbb{R})}\left|\p{a_{n}^{-1}\big(\max\limits_{1\le i \le n} X_i - b_{n}\big) \in A}-G_{\gamma}(A)\right|}{|A(n)|} = C,
\end{equation}
where $a_n = n g'(n)$ and $b_n = g(n).$
\end{prop}
\begin{rk}\label{rk1}
From the definitions of  $2\!-\!von\,Mises$ and $\mathcal{RV_{-\rho}}$, we have
    \begin{equation*}
      A(t) = L(t) t^{-\rho},
    \end{equation*}
    where  $L(\cdot)$ is a slowly-varying function at infinity. From Proposition~\ref{thResnick}, it follows that
    \begin{equation*}
    \sup _{A \in B(\mathrm{R})}\left|\p{a_{n}^{-1}\big(\max_{1\le i \le n} X_i -b_{n}\big) \in A}-G_{\gamma}(A)\right| \,\underset{n\to\infty}{\sim} \, L(n) n^{-\rho} .
     \end{equation*}
\end{rk}

$\blacktriangleright$  {\sc Falk and Marohn's result.}  
When assuming that the distribution $F$ belongs to the Fr\'echet maximum domain of attraction (taking $\gamma=-\alpha<0$), which is our general assumption,  Falk and Marohn propose von-Mises type conditions to obtain an asymptotic bound for the total variation distance in the EV theorem (see Theorem 3.1  in \citet{falk:marohn:1993}). They show (see Theorems 2.2 \& 2.4 in \citet{falk:marohn:1993}) that these conditions are equivalent to 
\begin{equation}\label{Falks_cond}
f_X(t) = c\, t^{-\a - 1} (1 + h(t)), \quad \text{with} \quad h(t) = O(t^{-\rho}), \,\rho > 0,
\end{equation}
where $f_X(t)$ denotes the density function of the cdf $F_X$.
Note that it implies that $f_X\in RV_{-\a-1}$ taking the slowly varying function as a constant.
\\
Looking at the converse statement, the authors propose conditions on the remainder terms of the von-Mises type conditions (see Theorem 3.2 in \citet{falk:marohn:1993}) to ensure the slowly varying function to be a constant, so that, under those conditions, the total variation distance in the EV theorem implies \eqref{Falks_cond}.

Note that Condition~\eqref{Falks_cond}  is very close to the characterization of the $2\RV$ class given in \eqref{Proko_cond} (see Lemma~\ref{lem:representation_2RV} in Appendix~\ref{A:secRV}); nevertheless these conditions do not imply each other.
\\[.5ex]

$\blacktriangleright$  {\sc Our proposition.}  Assuming $F\in$MDA(Fr\'echet), we can state the following:
\begin{prop}\label{cvRate-Frechet}
Suppose $\bar F_X\in\RV_{-\a}$, with $\a>0$. The rate of convergence for the EV theorem given  in Proposition~\ref{thResnick} holds when replacing the condition $g\in2\!-\!von\,Mises(-\a,-\rho)$ ($g$ being defined in \eqref{eq:def-g}), where $\rho>0$, by the condition
\begin{equation}\label{Proko_cond}
f_X(t) = c t^{-\a - 1} (1 + h(t)), \quad \text{with} \quad h\in\RV_{-\beta}, \;\beta > 0. 
\end{equation}
In this case, the function $A(\cdot)$ in \eqref{eq:thResnick} belongs to the $\RV_{-\rho}$ class with $\rho:= - \min\skk{1, \frac{\beta}{\alpha}}$. 
\end{prop}
Indeed, since $F_X\in$MDA(Fr\'echet),  the Potter bounds \eqref{potterbounds} can be directly established from a $2\RV$ condition on $g'$ via Proposition 4 in \citet{hua:joe:2011b} (see Lemma~\ref{lemma:potterbound} in Appendix~\ref{A:secRV}), which in turn follows from a $2\RV$ condition on $f_X$ by Lemma~\ref{lemma:2RVrelation} below (which proof is provided in Appendix~\ref{A:sec4lem}). This latter condition is equivalent to
our assumption~\eqref{Proko_cond}, as stated in \citet{hua:joe:2011b}, Lemma 3 (recalled as Lemma~\ref{lem:representation_2RV} in Appendix~\ref{A:secRV}). Once we have the Potter bounds, we can replicate exactly the proof of Proposition~\ref{thResnick} as given in \citet{Haan:Resnick:1996}, obtaining the same rate of convergence.
\begin{lem}\label{lemma:2RVrelation}
If $\displaystyle f_X \in 2\RV_{-\alpha - 1, -\beta}$, with $\alpha >0$ and $\beta >0$, then the derivative  $g'$ of $g$ defined in \eqref{eq:def-g}, satisfies $g'\in 2\RV_{\frac1\a - 1,\, \rho}$ , where $\rho :=-\min\skk{1, \frac{\beta}{\alpha}} (<0)$.
\end{lem}

Now, let us comment in Remark~\ref{rk:Resnick_Haan} the choice of the normalizing sequences $(a_n)$  and $(b_n)$ in \eqref{eq:EVT} of the EV theorem. 
First, notice that the limit in \eqref{eq:EVT} remains unchanged when replacing $(a_n)$ and $(b_n)$ with $(\tilde a_n)$ and $(\tilde b_n)$, as long as 
\begin{equation}\label{eq_for_norm_seq}
    \lim_{n \ra \infty}\frac{\tilde a_n}{ a_n} = 1 \quad \text{and}\quad \lim_{n \ra \infty}\frac{b_n - \tilde b_n}{a_n}  = 0.
\end{equation}
\begin{rk}\label{rk:Resnick_Haan}~
\vspace{-1ex}
\begin{enumerate}
\item[(i)] {\it Remark p.117 in \citet{Haan:Resnick:1996}.} When $F_X$ has a pdf, de Haan and Resnick provide universal norming sequences $(a_n)$ and  $(b_n)$ in \eqref{eq:EVT} for the EV theorem, namely 
\begin{equation*}
    a_n = n g'(n) \quad \text{and}\quad  b_n = g(n).
\end{equation*}
They point out that this choice is optimal in the following sense: If one changes $(a_n)$ and $(b_n)$ into  $(\tilde a_n)$ and $(\tilde b_n)$ such that for some $c_a,\, c_b \in \mathbb{R}$, we have
\begin{equation*}
  \lim_{n \ra \infty}  \frac{1}{A(n)} \left(\frac{\tilde a_n}{a_n} - 1 \right) =  c_a  \quad\text{and}\quad \lim_{n \ra \infty}\frac{1}{A(n)}\left(\frac{b_n - \tilde b_n}{a_n} \right) = c_b, 
\end{equation*}
then the rate of convergence in \eqref{eq:thResnick} remains the same up to a constant. But if any of the above limits is infinite, then the supremum in \eqref{eq:thResnick} is no longer $O(A(n))$.
\item[(ii)] If the distribution $F_X$  of a rv $X$ belongs to the $2\RV_{-\alpha,-\beta}$ class, with $\alpha >0$ and $\beta >0$, then we can choose the normalizing sequences $(a_n)$ and $(b_n)$ in \eqref{eq:EVT} for the EV theorem as
$$
    a_n = c^{1/\alpha} \,n^{1/\alpha}\;\,\mbox{with}\; \,c:= \lim\limits_{x \ra\infty} x^{\alpha}\bar F_X(x),  \;\;\mbox{and}  \; \,b_n = 0. 
$$
\item[(iii)] 
In fact, the sequences $(a_n)$ and $(b_n)$ suggested in (i) and (ii) do not provide an optimal rate of convergence, as can be observed in the following simple example. 
\begin{ex}\label{ex_optimal_norm}
Let $X_1, \cdots, X_n $  be a sample from a Pareto-Lomax distribution with parameter $\alpha >0$ and survival function
\begin{equation*}
    \bar F_X(x) = (1+x)^{-\alpha},  \; x >0.
\end{equation*}
The Pareto distribution belonging to the maximum domain of attraction of the Fr\'echet distribution, there exist $a_n > 0, b_n \in \mathbb{R}$, such that
\begin{equation}\label{1}
    F_X^n(a_n x + b_n) \underset{n\to\infty}{\ra} e^{-x^{-\alpha}}, \; x >0.
\end{equation}
If we choose $a_n = n^{1/\alpha}$ and $b_n = 0$, then we have 
$$
    F_X^n(a_n x + b_n)  = \left(1 - \left( 1 + n^{1/ \alpha}x\right)^{-\alpha}\right)^n 
$$
and the difference in \eqref{1} is $O(n^{-1/\alpha})$.
Now, if we take $b_n = -1$ instead of $0$, then 
\begin{equation*}
    F_X^n(a_n x + b_n)  = \left(1 - \frac{x^{-\alpha}}{n}\right)^n 
\end{equation*}
and the difference in \eqref{1} is $O(1/n)$, which  is much smaller for $\alpha>1$. 
So playing with the sequence $(b_n)$ in this simple example provides a strong improvement in the rate of convergence in \eqref{1}. 
\end{ex}
\end{enumerate}
\end{rk}
%
\subsection{Rate of convergence for MRV-Normex }
%

\paragraph{Definitions -}
While the standard MRV definition is recalled in Appendix~\ref{A:secRV} (see Definition~\ref{def:mrv}), let us give 
an equivalent MRV definition using the pseudo-polar representation, to better understand why the introduction of this multi-normex version, namely MRV-Normex.
\begin{defi}\label{def:MRVpolar} 
The random vector $\bX\in \MRV_{-\alpha}$, with $\a>0$, if there exists a $d$-dimensional random vector $\bTheta$ with values in the unit sphere ${\cal S}_{1}$ in $\R^d$ w.r.t. the norm $\|\cdot\|$, such that, $\forall t >0$,
\begin{equation}\label{def:polarMRV}
\frac{\p{\|\bX\|>t\,u,\;\bX/\|\bX\|\in \cdot }}{\p{\|\bX\|>u}} \, \conv  \, t^{-\a}\,\p{\bTheta \in \cdot}\quad \text{as} \quad u\to \infty.
\end{equation}
\end{defi}
Using this latter definition of MRV, in particular the random vector $\bTheta$, and Remark~\ref{rk:Resnick_Haan}, we can define the {\it  MRV-Normex} distribution as follows.
\begin{defi}\label{def:NormexMRV}
The so-called {\bf  MRV-Normex} distribution function  is defined, for $\mathbf{B} \subset  \R^d$, by:
\begin{equation}\label{eq:normex-cdf-MRV}
 GM_n(\mathbf{B}):=  \p{H_{\a,n}\, \bTheta + Z \in \mathbf{B} }.
\end{equation}
 with 
$H_{\a,n}:=a_n H_\a + b_n$
where the random variable $H_\a$ (with $\a>0$) is Fr\'echet distributed ({\it i.e.} $\displaystyle  \p{H_\a \le x}=e^{-x^{-\a}}$, for $x>0$) and independent of the random vector $\bTheta$ introduced in \eqref{def:polarMRV}, the normalizing sequences satisfy the standard conditions of EV theorem, namely $a_n=c\,n^{1/\a}$ with $c^{\,\a} := \lim\limits_{y \ra\infty} y^{\alpha}\bar F_{\|X\|}(y) $, and $b_n=0$.
The $d$-dimensional random vector $Z$, also assumed to be  independent of $\bTheta$, is, conditionally to the event $(H_{\a,n} = y)$, with $y>0$, normally $\mathcal{N}_{(n-1)\mathbf{\mu}(y),\, (n-1) \Sigma(y) }$-distributed, the mean vector and covariance matrix being those of the truncated distribution $F_{\mathbf{X}\,|\, \| \mathbf{X}\|}(\cdot | y)$ defined in \eqref{eq:def_trunc_distr}.
\end{defi}
The MRV-Normex cdf can be rewritten as
\begin{equation}\label{eq:normex-cdf-MRV2}
   GM_n(\mathbf{B}):= \int\limits_{0}^{\infty} f_{ H_{\a,n}}(y) \,\p{ y \, \bTheta + Z_y \in \mathbf{B} } \mathrm{d}y,
\end{equation}
where $Z_y$ is $\mathcal{N}_{(n-1)\mathbf{\mu}(y),\, (n-1) \Sigma(y)}$-distributed and independent of $\mathbf{\Theta}$.\\[1ex]
%
Note that we may choose other normalizing sequences that satisfy \eqref{eq_for_norm_seq} in the definition of MRV-Normex distribution.  Moreover, this distribution may be defined mathematically for any MDA, using the Generalized Extreme Value distribution and the appropriate normalizing sequences $(a_n)$ and $(b_n)$, as the EV theorem  (see \eqref{eq:EVT}) holds for any MDA. Nevertheless, we focus on the MDA(Fr\'echet) for which a naive use of a light tail distribution for the sum would lead to a large error. 

\paragraph{Rate of convergence for the MRV-Normex approximation -}

We have the following result, which proof can be found in Appendix~\ref{A:sec4lem}.
\begin{theo}\label{th:MRV result}
Let $\mathbf{X}_{1}, \ldots, \mathbf{X}_{n}$ be i.i.d. random vectors with parent random vector  $\mathbf{X}$ with values in $\R^{d}$. Assume the following conditions: 
\begin{itemize}
\item[$(C1)$] given in Theorem~\ref{th:Result-max} (namely, $F_{\mathbf{X}\,|\, \| \mathbf{X}\|(\cdot\,|\,y)}$ non degenerate $\forall y > 0$);
\item[$(M_{\|\cdot\|})$] The distribution of the rv $\|\mathbf{X}\|$ is absolutely continuous 
and its pdf $f_{\|\mathbf{X}\|}\in2\RV_{-\alpha-1,-\beta}$  with $\alpha > 0$, $\beta > 0$;
\item[$(M_{\Theta})$] There exists a function $A$ such that $A(t) \ra 0$, $|A(t)| \in  \RV_{-\rho}$ with $\rho > 0$, and
\begin{equation*}
  \sup\limits_{\mathbf{B} \in \mathcal{S}_1} \left| \p{\frac{\bX}{\|\bX\|}\in \mathbf{B}\,\Big|\, \|\bX\|>t} - \p{\bTheta \in \mathbf{B} } \right| \,\underset{t\to\infty}{\sim} \, A(t),
\end{equation*}
where the supremum is taken over all measurable subsets of $\mathcal{S}_1.$
\end{itemize}
Then, for any $\alpha \in (2,3]$, $\beta>0$ and $\rho > 0$, there exists a slowly varying function $L(\cdot)$ such that
\begin{equation}\label{eq_main_MRV}
   \sup _{\mathbf{B} \in \mathcal{C}}  | \p{\mathbf{S_n} \in \mathbf{B}} - GM_n(\mathbf{B}) | \, \leq  \,\left( n^{- \frac{1}{2} + \frac{3-\alpha}{\alpha} } + n^{  - \frac{\rho}{\alpha}} +  n^{-\min(1,\frac{\beta}{\alpha})}\right) L(n),
\end{equation}
where $\mathcal{C}$ is the class of all Borel-measurable convex subsets of $\R^{\mathbf{d}}$ and $GM_n$ is defined in \eqref{eq:normex-cdf-MRV2}.\\[.7ex]
\end{theo}

\begin{rk}\label{rk:theoMRV}~ \\[-4ex]
	\begin{enumerate}
	\item If $\rho > \alpha$ and $\beta > \alpha$, then  \eqref{eq_main_MRV} rewrites as
$$
\sup _{\mathbf{B} \in \mathcal{C}}  | \p{\mathbf{S_n} \in \mathbf{B}} - GM_n(\mathbf{B}) | \, \leq  \,n^{- \frac{1}{2} + \frac{3-\alpha}{\alpha}} L(n),
$$
providing the same rate of convergence given in Theorem \ref{th:Result-max} for $d$-Normex.
	\item If the supremum considered in $(M_{\Theta})$ has a faster speed of convergence than that of  $\RV$, then we can set  $\rho = \infty$ and exclude the term $n^{  - \frac{\rho}{\alpha}}$ from inequality \eqref{eq_main_MRV}. 
	\item Discussion on Condition $(M_{\Theta})$:
	\begin{enumerate}
		\item First note that assuming $\|\mathbf{X}\|\in RV$ and  Condition $(M_{\Theta})$, which requires uniform convergence, is closely related to the MRV definition \eqref{def:polarMRV}. Nevertheless, replacing this technical condition $(M_{\Theta})$ with \eqref{def:polarMRV} might be investigated further.
		\item {\it Independent polar coordinates.} 
If the norm $\|\mathbf{X}\|$  and the direction $\mathbf{X}/ \|\mathbf{X}\|$ of $\mathbf{X}$ are independent, then $(M_{\Theta})$ is satisfied. Moreover, in the absolutely continuous case, $\|\mathbf{X}\|$  and  $\mathbf{X}/ \|\mathbf{X}\|$ are independent if and only if there exist non negative functions $h_1$ and $h_2$ such that the pdf  of $\mathbf{X}$ can be factorized as:
\begin{equation}\label{eq:factorization}
    f_{\mathbf{X}}(\mathbf{x}) = h_1(\|\mathbf{x}\|) \,h_2(\mathbf{x}/ \|\mathbf{x}\|), \ \mathbf{x} \in \mathbb{R}^d.
\end{equation}
		 \item  In some cases, the distribution of the random vector $\mathbf{\Theta}$ can be made explicit, as follows. For any $t > 0$ and $\mathbf{B} \in \mathcal{S}_1,$
\begin{equation*}
\p{\frac{\mathbf{X}}{\|\mathbf{X}\|} \in \mathbf{B} \,\Big|\, \|\mathbf{X}\| \geq t}  = \frac{1}{\overline{F}_{\|\mathbf{X}\|} (t)}\int\limits_{\frac{{\mathbf{u}}}{{\|\mathbf{u}\|}} \in \mathbf{B}, \|\mathbf{u}\| \geq 1} f(t \mathbf{u})\, t^d  \mathrm{d}\, \mathbf{u}. 
\end{equation*}
 If the latter integral has a limit, as $t \ra \infty$, of the form
$\displaystyle
 \int\limits_{\frac{{\mathbf{u}}}{{\|\mathbf{u}\|}} \in B, \|\mathbf{u}\| \geq 1} f_{\Theta}(\mathbf{u}) \,\mathrm{d}\mathbf{u}$, \\
then this limit defines the distribution of $\mathbf{\Theta}:$
\begin{equation*}
     \p{\mathbf{\Theta} \in \mathbf{B}}  =  \int\limits_{\frac{{\mathbf{u}}}{{\|\mathbf{u}\|} }\in \mathbf{B}, \|\mathbf{u}\| \geq 1} f_{\Theta}(\mathbf{u}) \,\mathrm{d}\mathbf{u},
\end{equation*}
from which we deduce that
\begin{equation*}
    \sup_{\mathbf{B}} \left| \p{\frac{\mathbf{X}}{\|\mathbf{X}\|} \in \mathbf{B} \,\Big|\, \|\mathbf{X}\| \geq t} -  \p{\mathbf{\Theta} \in \mathbf{B}} \right| = \frac{1}{2} \int\limits_{\|\mathbf{u}\| \geq 1}  \left| \frac{ f(t \mathbf{u}) \,t^d }{\overline{F}_{\|\mathbf{X}\|} (t)}  - f_{\Theta}(\mathbf{u})  \right| \mathrm{d}\mathbf{u}.
\end{equation*}
Further details are given in Appendix~\ref{A:sec4disc}.
		\end{enumerate}
	\end{enumerate}
\end{rk}

\subsection{Examples}\label{ssec:exs-MRV}

Let us develop two examples such that the marginal distributions are Pareto-Lomax, as in Example~\ref{ssec:exs-dNormex}, so that it allows  comparison with $d$-Normex.
We consider two cases for the parent random vector $\mathbf{X}$, when assuming its components to be, on one hand, independent, on the other hand, related with a survival Clayton copula. These are standard examples considered in the actuarial and risk literature (see e.g. \citet{das2020}), in particular in reinsurance context for the Clayton copula is (see e.g. \citet{dacoAL:2018} and references therein).
This second example includes itself two cases, when the polar coordinates of the considered vector  $\mathbf{X}$ are dependent (but asymptotically independent), and when they are independent. This latter case corresponds to Example~\ref{ssec:exs-dNormex}. 

We check that the conditions of Theorem~\ref{th:MRV result} are satisfied whenever $\a\in(2,3]$ (recall that this constraint on $\alpha$ appears only for the analytical comparison with the generalized Berry-Esseen inequality), but apply the MRV-Normex distribution for any positive $\alpha$, as the construction via asymptotic theorems remains valid whatever the value of this parameter.

\subsubsection{Independent Pareto-Lomax marginals}

Assume the components of the random vector $\mathbf{X}$ to be iid with Pareto-Lomax($\a$)  distribution \eqref{def:one_dim_par}. Then its (non truncated) moments remain the same as in Example~\ref{ssec:exs-dNormex} and its covariance matrix is diagonal. 
As we developed Example~\ref{ssec:exs-dNormex} with the $L_1$-norm, let us switch here to the  $L^\infty$-norm $\|\cdot\|_{\infty}$, more convenient in terms of computations in this framework.

The cdf of the norm of the vector $\mathbf{X}$ being, for $\alpha>0$,
\begin{equation*}
   F_{\|\mathbf{X}\|_\infty}(y) = \p{\max\skk{X^{(1)}, \cdots, X^{(d)}} \leq y} = (1  - (1+y)^{-\alpha})^d, \quad\text{for} \; y > 0, 
\end{equation*}
straightforward calculations give the following expressions for the truncated moments:
\begin{equation*}
\begin{split}
  \mu^{(1)}(y) & = \e\left( X^{(1)}\,|\, \|\mathbf{X}\|_\infty \leq y\right) = \e\left( X^{(1)}\,|\, X^{(1)} \leq y\right) \\
  & = \frac{1}{1 - (1+y)^{-\alpha}} \left(\frac{1}{\alpha-1}\left(1 - (1+y)^{-\alpha + 1} \right)- y(1+y)^{-\alpha}\right),
\end{split}
\end{equation*}
\begin{equation*}
\e\left( \left(X^{(1)}\right)^2\,\big|\, \|\mathbf{X}\|_\infty \leq y\right) 
  = \frac{1}{1 - (1+y)^{-\alpha}} \dfrac{2-(1+y)^{-\alpha}\left( \alpha y\left(\left(\alpha-1\right)y+2\right)+2\right)}{\left(\alpha-2\right)\left(\alpha-1\right)},
\end{equation*}
and zero for the truncated covariances. 

When looking for the distribution of $\mathbf{\Theta}$, notice that, for any $i\neq j$, for
any $\varepsilon_i > 0$ and $\varepsilon_j > 0$, we have 
$$
    \lim_{n \ra \infty} \p{ X^{(i)} > \varepsilon_i \,t, \,X^{(j)} > \varepsilon_j\, t\;\big|\, \|\mathbf{X}\|_\infty > t  } = 0.
$$
Therefore, the distribution of the random vector $\mathbf{\Theta}$ is discrete on the unit sphere $\mathcal{S}_1$, with values given by the basis vectors $\mathbf{e}_i = (0,\cdots,1,\cdots,0)$ (where $1$ is for the $i$-th component).
It is straightforward to verify that $F_{\|\mathbf{X}\|}\in 2\RV_{-\alpha,-\alpha}$, so that $(M_{\|\cdot\|})$ is satisfied, and that Condition $(M_{\Theta})$ holds with auxiliary function $A(\cdot) \in \RV_{-\alpha}$.

Finally, one may chose the normalizing sequences as $a_n = (d \,n)^{1/\alpha}$ (according to Remark ~\ref{rk:Resnick_Haan}(i) and $b_n = -1$ as in Example~\ref{ex_optimal_norm} of Remark ~\ref{rk:Resnick_Haan}(iii). 

The numerical implementation of this MRV-Normex approximation is then developed in Section~\ref{sec:QQplots}, along with that of the $d$-Normex one, for any positive $\alpha$, both multi-normex methods being compared to the Gaussian approximation whenever $\a\ge 2$. The QQ-plots are drawn in Figure~\ref{fig:qqplot_indep 2,3}.

\subsubsection{Pareto-Lomax marginal distribution with survival Clayton copula}\label{ex:ParetoClayton}

We introduce in this example some dependence among the components of $\mathbf{X}$, choosing a survival Clayton copula (so, with upper tail dependence). 
To lighten the expressions of the computed moments, we choose $d=2$.

We consider $\mathbf{X}=\left(X_{1}, X_{2}\right)$ with identical Pareto-Lomax $(\alpha, 1)$ marginal distributions, $\alpha>1$, {\it i.e.}
$$
\bar{F}_{1}(x)=\bar{F}_{2}(x)=(1+x)^{-\alpha}, \quad \forall x>0,
$$
and survival Clayton copula on $[0,1]^{2},$  with parameter $ \theta>0$, defined by
$$
\mathbb{P}\left(X_{1}>x_{1}, X_{2}>x_{2}\right)=\left[\left(1+x_{1}\right)^{\alpha \theta}+\left(1+x_{2}\right)^{\alpha \theta}-1\right]^{-1 / \theta},
$$
with pdf \;$\displaystyle
f\left(x_{1}, x_{2}\right)= \alpha^{2}(1+\theta)\left(1+x_{1}\right)^{\alpha \theta-1}\left(1+x_{2}\right)^{\alpha \theta-1}\left(\left(1+x_{1}\right)^{\alpha \theta}
+\left(1+x_{2}\right)^{\alpha \theta}-1\right)^{-\frac{1}{\theta}-2}$.

Considering the $L^\infty$-norm, $\|\cdot\|_{\infty}$, the survival cdf of the norm of the vector is:
\begin{equation*}
   \overline{F}_{\|\mathbf{X}\|_{\infty}}(t) = \p{\max\skk{X_1, X_2} > t} = 2(1+t)^{-\alpha} -(2(1+t)^{\alpha \theta} - 1)^{-1/\theta}. 
\end{equation*}

We computed the truncated moments (of order 1 and 2), using an integral calculator (based on  \textit{Maxima}, a computer algebra system developed by W. Schelter, MIT), providing explicit but long expressions, that is why we do not display them here 
(and give them in the online Appendix with the html version of all plots).

For positive $u_1$ and $u_2$ such that $\max\skk{u_1,u_2} \geq 1$,  we can write
$$
\frac{ f(t \mathbf{u})\, t^d }{ \overline{F}_{\|\mathbf{X}\|_\infty}(t) }  =    \frac{\alpha^2(1+ \theta) \left(u_1 + \frac{1}{t} \right)^{\alpha\theta - 1} \left(u_2 + \frac{1}{t} \right)^{\alpha\theta- 1}}{\left(2 - 2^{-1/\theta}\right)(1 + O(\frac{1}{t}))\left( \left(u_1 + \frac{1}{t} \right)^{\alpha\theta} + \left(u_2 + \frac{1}{t} \right)^{\alpha\theta}\right)^{1/\theta  +2 } }
$$
from which we deduce the limit as $t\to \infty$, namely
$$
    \lim\limits_{t \ra \infty}   \frac{ f(t \mathbf{u})\, t^d }{ \overline{F}_{\|\mathbf{X}\|_\infty}(t) }  = \frac{\alpha^2(1+ \theta) u_1^{\alpha\theta - 1} u_2^{\alpha\theta- 1}}{\left(2 - 2^{-1/\theta}\right)\left( u_1^{\alpha\theta} + u_2^{\alpha\theta}\right)^{1/\theta  +2 } } =: f_{\Theta}(\mathbf{u}). 
$$
Note that, for $ \|\cdot\| =  \|\cdot\|_\infty$, the function $\frac{ f(t \mathbf{u}) t^d }{ \overline{F}_{\|\mathbf{X}\|}(t) } $ depends on $t$, therefore the random variable $\|\mathbf{X}\|$ and random vector $\mathbf{X}/ \|\mathbf{X}\|$ are not independent. 

They will be independent when replacing the $L^\infty$-norm with the $L^1$-norm ($ \|\cdot\| =  \|\cdot\|_1 $) and choosing $\alpha\theta = 1$; this corresponds to the Pareto-Lomax Example~\ref{ssec:exs-dNormex}.  

Let us check the conditions of Theorem~\ref{th:MRV result}. It is straightforward to check that $F_{\|\mathbf{X}\|}\in 2\RV_{-\alpha, -\min(\alpha\theta, 1)}$, so that $(M_{\|\cdot\|})$ is satisfied.  
Some computations are required for Condition~$(M_{\Theta})$. 
We keep the maximum norm, {\it i.e.} $\|\cdot\| =  \|\cdot\|_\infty$, so that we exhibit an example with dependence between the polar coordinates (but with asymptotic independence), but consider the case $\alpha\theta = 1$ to simplify the computations. We obtain:
\begin{align}\label{eq:TVdifference1}
    & \sup_B \left| \p{\frac{\mathbf{X}}{\|\mathbf{X}\|} \in B \,|\, \|\mathbf{X}\| \geq t} -  \p{\Theta \in B} \right| = \frac{1}{2} \int\limits_{\|\mathbf{u}\| \geq 1}  \left| \frac{ f(t \mathbf{u}) t^d }{\overline{F}_{\|\mathbf{X}\|} (t)}  - f_{\Theta}(\mathbf{u})  \right| \mathrm{d}\mathbf{u} \nonumber\\
    & =\frac{\alpha(\alpha+1)}{2t} \int\limits_{\|\mathbf{u}\| \geq 1} t \, \left| \frac{t^2}{ \overline{F}(t)  (1 + t|\mathbf{u}|)^{\alpha + 2}} - \frac{1}{c_{\alpha} |\mathbf{u}|^{\alpha+2}}  \right| \mathrm{d}\mathbf{u}\nonumber \\
    & =\frac{\alpha(\alpha+1)}{2t} \int\limits_{\|\mathbf{u}\| \geq 1}  t\,\left| \frac{  c_{\alpha}|\mathbf{u}|^{\alpha+2} -  \overline{F}(t)   t^{\alpha} ( |\mathbf{u}| + \frac{1}{t})^{\alpha + 2}    }{  \overline{F}(t)  \, t^{\alpha}( |\mathbf{u}| + \frac{1}{t})^{\alpha + 2} \,c_{\alpha} |\mathbf{u}|^{\alpha+2}}\right| \mathrm{d}\mathbf{u},
\end{align}
where $c_{\alpha} = (2 - 2^{-\alpha})$ and $|\mathbf{u}| = u_1 + u_2$. \\
We can easily find a majorant of type $c/|\mathbf{u}|^{\alpha+2}$ for the integrand \eqref{eq:TVdifference1}, and, noticing that this integrand converges, as $t\ra \infty$, to
$\displaystyle  \frac{ \left| \hat{c}_{\alpha}|\mathbf{u}| + c_{\alpha}(\alpha+2) \right|   }{ c_{\alpha}^2 |\mathbf{u}|^{\alpha+3}}$, 
with $\hat{c}_{\alpha} := 2^{-\alpha - 1} - 2$,
we can conclude, via the dominated convergence theorem, that 
\begin{equation}\label{eq:TVdifference2} 
\int\limits_{\|\mathbf{u}\| \geq 1}  t\,\left| \frac{  c_{\alpha}|\mathbf{u}|^{\alpha+2} -  \overline{F}(t)   t^{\alpha} ( |\mathbf{u}| + \frac{1}{t})^{\alpha + 2}    }{  \overline{F}(t)   t^{\alpha}(( |\mathbf{u}| + \frac{1}{t})^{\alpha + 2}) c_{\alpha} |\mathbf{u}|^{\alpha+2}}\right| \mathrm{d}\mathbf{u} \underset{t\to\infty}{\ra}
\int\limits_{\|\mathbf{u}\| \geq 1}    \frac{ \left| \hat{c}_{\alpha}|\mathbf{u}| + c_{\alpha}(\alpha+2) \right|  }{ c_{\alpha}^2 |\mathbf{u}|^{\alpha+3}} \mathrm{d}\mathbf{u}  < \infty. 
\end{equation}
Combining \eqref{eq:TVdifference1} and \eqref{eq:TVdifference2} provides that Condition $(M_{\Theta})$ holds with $A(t)  = \frac{C_{\alpha}}{t} $ for some constant $C_{\alpha} \in (0,\infty)$. 

As in the previous example (case of independent components), one may chose the normalizing sequences $a_n = (c_\alpha n)^{1/\alpha}$ (see Remark~\ref{rk:Resnick_Haan}, (i))  and $b_n = -1$ (see Remark~\ref{rk:Resnick_Haan}, (iii)).
 
We refer to the next section for the numerical implementation of this example; see Figures~\ref{fig:qqplot_clayton 2,3}~\&~\ref{fig:qqplot_clayton1 2,3} for the QQ-plots. As we are in dimension 2, we also provide (ranked) scatter plots to compare the various approximations of the sum distribution, varying the parameters, in particular the dependence parameter $\theta$; see Figures~\ref{fig:scatterplot23}~\&~\ref{fig:ranks101}.

\section{QQ-plots of the various examples, illustrating both versions of the multi-normex method }\label{sec:QQplots}

\subsection{Construction of the QQ-plots}

Considering the examples given so far, we illustrate the benefit of the multi-normex method on $d-$dimensional QQ-plots based on geometrical quantiles.
We propose a brief overview of this latter notion in Appendix~\ref{A:secGeomQuant}, but refer mainly to \citet{DharAl:2014} for definitions and detailed explanations. 
Here are a few first key ideas of these objects, to help interpret the plots displayed in this section.
\\[.7ex]
The geometrical quantile, as given in Definition~\ref{def:geom quantile} (see Appendix~\ref{A:secGeomQuant}), is a generalization to higher dimension of the $1$-dimensional quantile that can be defined as the solution of some optimization problem (see $\eqref{def:quantile}$ in Appendix~\ref{A:secGeomQuant}). So, one can formulate the same optimization problem in the $d-$dimensional case, for which the solution will be named {\it geometrical quantile}. In this way, the geometrical quantile is a point of $\mathbb{R}^d$. While $1$-dimensional quantiles are parameterized by the interval $(0,1)$, multidimensional ones are parameterized by a  $d-$dimensional unit ball from $\mathbf{R}^d$. In this context, we refer to vectors of the unit ball as {\it levels}.
Although geometrical quantiles reflect the structure of a $d-$dimensional distribution, they are abstract objects and do not have a nice interpretation as $1$-dimensional quantiles do. Only the median has a geometrical sense: given a random vector, the median is a point of $\mathbb{R}^d$ such that the overall sum of the distances from this point to all values of the random vector is minimum (note that distances are multiplied by the 'probabilities' that the vector takes the considered values, respectively).
Moreover, if the vector has a finite second moment, then its extreme quantiles will share the same speed of convergence towards infinity as any vector having the same covariance matrix (see \citet{girard2015b,girard2015a}).
Nevertheless, we can construct QQ-plots with these geometrical quantiles, with the aim at comparing $d-$dimensional distributions by comparing the plots between each other.
\\[.5ex]
The construction of QQ-plots is similar as in the $1$-dimensional case. Considering two distributions on $\mathbb{R}^d$, we can solve the optimization problem for a fixed number of levels, say $N$, and obtain two sets of $N$ geometrical quantiles: $\skk{\mathbf{q}_1, \dots, \mathbf{q}_N}$ for the first distribution, and  $\skk{\mathbf{q}'_1, \dots, \mathbf{q}'_N}$ for the second.
Then we draw a  $2$-dimensional QQ-plot for each component of the $\R^d$-geometrical quantiles to get a visualization, obtaining $d$ QQ-plots: $\left\{\left(q_{i,j}, q'_{i,j}\right); j=1,\cdots,N\right\}$, for $i=1,\cdots,d$.
In \citet{DharAl:2014}, the authors proved that all the points (pairs of quantiles) lie close to a straight line with slope $1$ and intercept $0$ if and only if the two distributions are equal.

Turning to our previous examples, we consider multidimensional distributions with Pareto-Lomax($\a$) marginals defined in \eqref{def:one_dim_par}, varying their heaviness through the parameter $\alpha \in \skk{1.5, 2.3, 3.5}$, with the different structures of dependence given in Examples~\ref{ssec:exs-dNormex}~\&~\ref{ssec:exs-MRV}. We choose $n = 52$ and $d \in \skk{2,3}$. For each case, we evaluate both multi-normex distributions and the normal distribution (except for $\alpha=1.5$) obtained by the CLT; those distributions are evaluated empirically, via simulations, each sample being of size $10^7$. Then we proceed to their comparisons via the QQ-plots. \\
We construct the QQ-plots in 4 main steps:\\[-4ex]
\begin{enumerate}
\item  Simulate all the distributions to be compared: the distribution of the sum $\mathbf{S}_{n}$, the normal distribution from CLT, the d-Normex distribution $G_n$, and the MRV-Normex distribution $MG_n$. Namely,
	\begin{enumerate}
    	\item To obtain a simulated sample (of size $10^7$) for the sum $\mathbf{S}_{n}$, we simulate $n\times10^7$ random vectors $\mathbf{X}$ from the considered distribution. 
    \item To obtain a simulated sample from the $d$-Normex distribution defined in \eqref{eq:normex-cdf}: First, we build $n$ samples (of size $10^7$) from the $d$-dimensional multivariate Pareto distribution, from which we deduce a sample (of size $10^7$) for the ($d$-dimensional) maximum $\mathbf{X}_{(n)}$ (see \eqref{ordered statistics}). Second, for each element of the latter sample, we calculate its norm $y = \|\mathbf{X}_{(n)}\|$ and simulate a normal vector with mean $(n-1)\mu(y)$ and covariance matrix $(n-1)\Sigma(y)$ (described in Definition~\ref{def:Normex}), collecting then $10^7$ Gaussian vectors. Finally, we sum maximum and normal vectors to produce a sample (of size $= 10^7$)  from the $d$-Normex distribution.
    	\item In order to simulate the MRV-Normex distribution defined in \eqref{eq:normex-cdf-MRV}, we start simulating a sample (of size $10^7$) for the vector $\mathbf{\Theta}$ (representing the direction) and an independent sample for the Fr\'echet distributed rv $H_{\a,n}$ introduced in \eqref{eq:normex-cdf-MRV2}. Next, we collect the  $10^7$ normal vectors with mean  $(n-1)\mu(y)$ and covariance matrix $(n-1)\Sigma(y)$, where, now, $y =H_{\a,n}$. Finally, we aggregate all the constructed samples according to \eqref{eq:normex-cdf-MRV}.  
	\end{enumerate}
\item Fix the set of levels $ \mathcal{L} \subset \left\{ \mathbf{v} \in \mathbb{R}^{d},\|\mathbf{v}\|<1\right\}$ with different lengths and directions.  We choose $10$ lengths, $||\mathbf{v}|| \in \skk{0    , 0.2   , 0.4   , 0.6   , 0.8   , 0.9   , 0.9225, 0.945 , 0.9675, 0.99 }$ (half for the body of the distribution and half for its tail),  and all the directions with factor $\frac{\pi}{4}$. It represents a total of $235$ vectors for $d=3$ and of $145$ for $d = 2$. 
\item Calculate the geometrical quantiles for the simulated samples of all considered distributions. It means to solve the numerically optimization problem \eqref{def:emp-mlti-quantile} for each empirical distribution and  for all levels $\mathbf{v}$ from the set $\mathcal{L}$.
\item Draw a QQ-plot for the three pairs (or two, when the CLT cannot be applied): (sum, CLT), (sum, $d$-Normex) and (sum, MRV-Normex).
\end{enumerate}
Note that the numerical implementation has been performed with Python (SciPy library).  The scipy.optimize.minimize function based on the quasi-Newton method of Broyden, Fletcher, Goldfarb, and Shanno (see p.136 in \citet{Nocedal2006}) has helped solve the numerical optimization problem~\eqref{def:emp-mlti-quantile}. The gradient for the function in \eqref{def:emp-mlti-quantile} has been calculated analytically. The computation time of one geometrical quantile was, on average,  $53$ seconds (on the computer i7 2GHz, 16 GB RAM).

\subsection{Multivariate Pareto-Lomax distribution - Example~\ref{ssec:exs-dNormex}}
\label{ssec:plots_1}

We first show the QQ-plots for the multivariate Pareto distribution, developed in Section~\ref{ssec:exs-dNormex}, considering the dimension $d= 3$, the number of summands $n=52$, and varying $\alpha$.
We choose a small number of summands to highlight the performance of multi-normex methods, even in this case. The QQ-plots are given on a same row for each of the three components, and for each given approximation method (CLT, $d$-Normex and MRV-Normex). As expected, the plots are similar componentwise.
A zoom of the center of the graph is given in the upper left corner of each plot, since it is the only region where the quantiles are so concentrated, making it 
more difficult to judge if they form a straight line (Note that an html version of all plots, with the possibility to zoom anywhere, is available on the website \href{http://crear.essec.edu/members}{crear.essec.edu}). 
It should be noted that the center corresponds not only to levels with length less than or equal to $0.9$ (marked in blue), but also to extreme levels (in red) with length greater than $0.9$. 

\begin{figure}[H]
\centering
\begin{tabular}{c}
Multivariate Pareto ($\a=2.3$)\\
    \includegraphics[width=165mm]{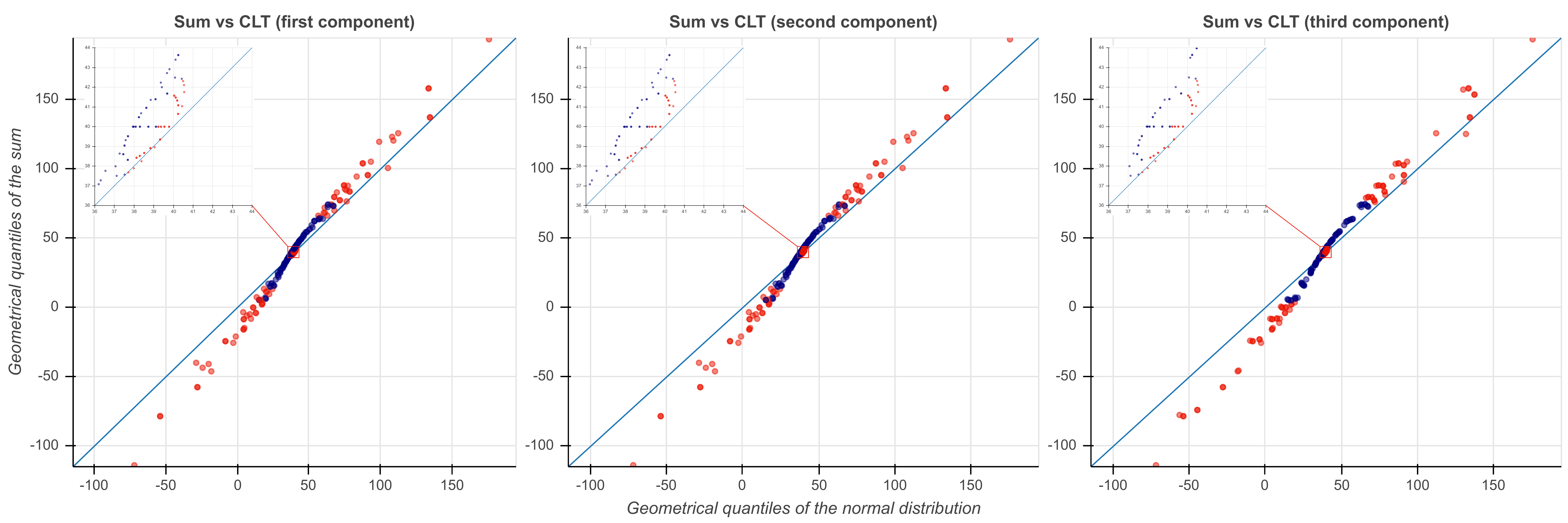} \\
    \includegraphics[width=165mm]{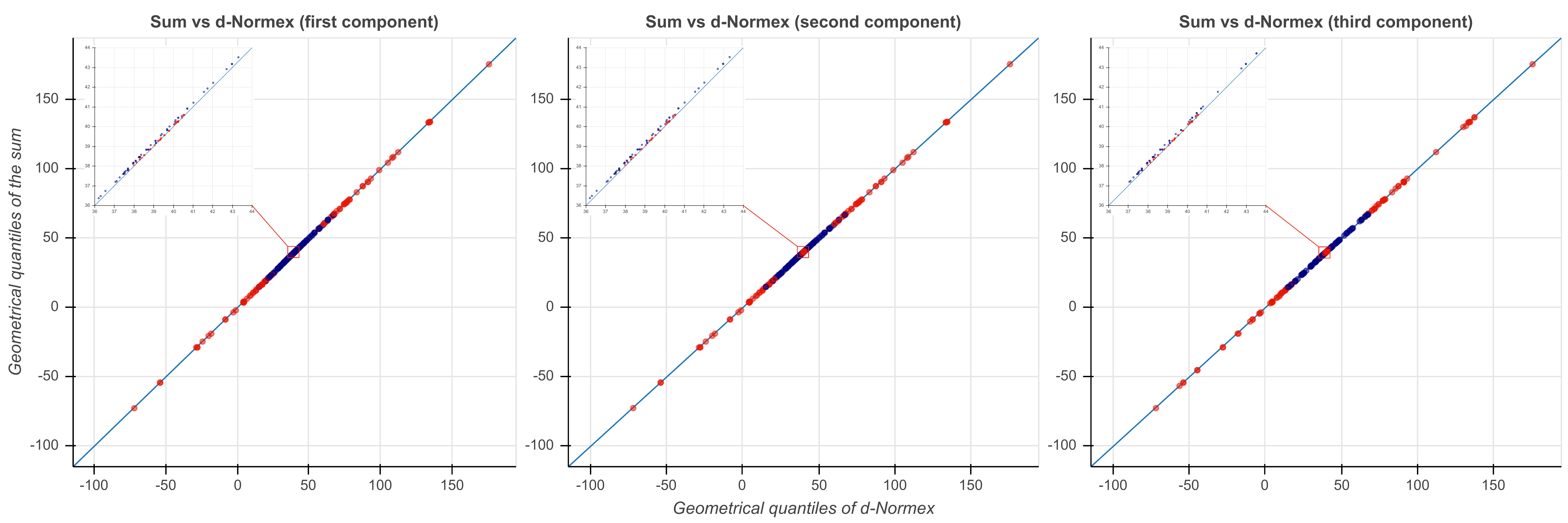} \\
    \includegraphics[width=165mm]{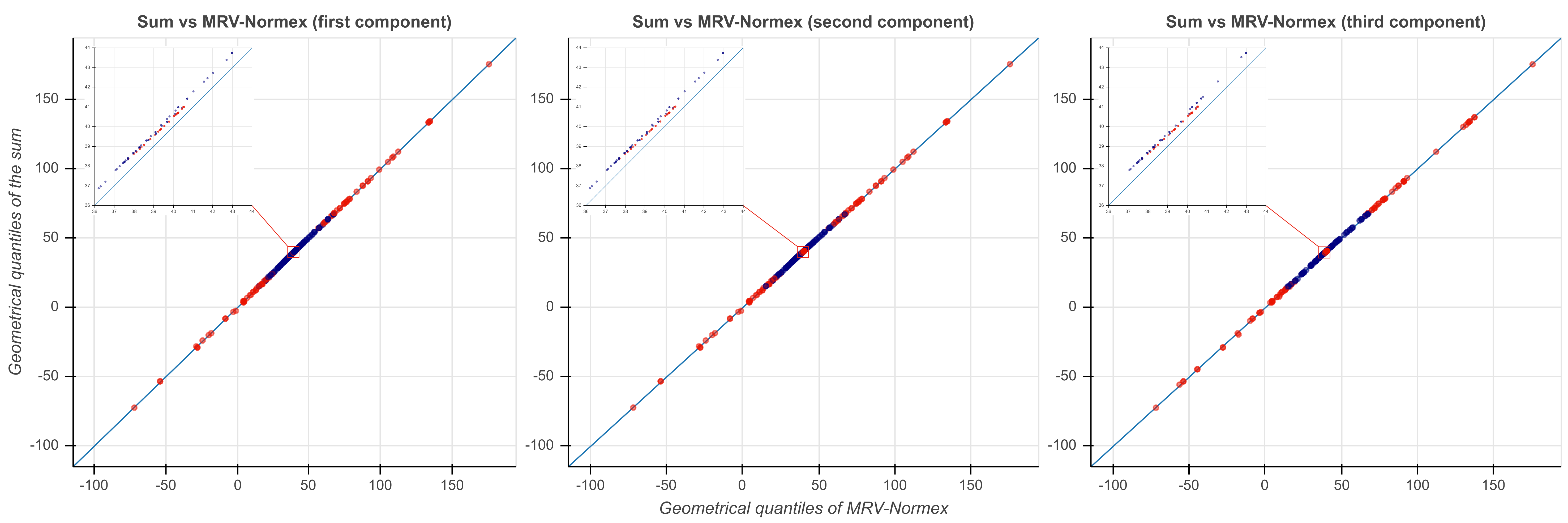} \\
     \end{tabular}
\caption{ \small  $3$-dimensional QQ-plots for the empirical distribution (sample size $= 10^7$) of the sum of $52$ iid trivariate Pareto ($\alpha = 2.3$) random vectors, with three different approximations of the sum distribution: CLT (first row), $d$-Normex (second row) and MRV-Normex (third row). Each column corresponds to a component (from the 1st to the 3rd). The red points on the plots correspond to extreme geometric quantiles (when the norm of the parameterized vectors is greater than $0.9$)}
\label{fig:qqplot-23}
\end{figure}

{\sc Case $\alpha \in (2,3]$.}
Figure~\ref{fig:qqplot-23} exhibits the QQ-plots in the case $\alpha = 2.3,$ i.e. in the framework of the multi-normex theorems, when using the Gaussian approximation (via the CLT) and both multi-normex approximations. We observe that all the points (pairs of quantiles) for $d$-Normex and MRV-Normex lie much closer to the line with slope $1$ and intercept $0$ than for the Gaussian distribution. Thus, we can see numerically that both multi-normex approximations better describe the distribution of the vectorial sum $\mathbf{S}_n$, as proved analytically in Theorems~\ref{th:Result-max} and~\ref{th:MRV result}. While the QQ-plots look quite similar for the two versions of multi-normex approximations, when zooming, the fit is slightly better when using the $d$-Normex distribution than the MRV-Normex one that uses the Fr\'echet distribution for  the distribution of the rescaled maximum. 

 {\sc Arbitrary $\alpha > 0$.}
Here we consider other examples of heavy tails than the case $\a\in (2,3]$, to illustrate the benefit of using Normex distributions rather than applying the generalized CLT with a Gaussian distribution (for finite variance) or a stable one (when the variance is infinite).  We do it numerically as the upper bound of the Berry-Esseen inequality does not allow us to compare analytically the rates of convergence in terms of the number of summands $n$. Nevertheless, the constant given in terms of $\a$ and the number $k$ of largest order statistics when trimming the sum, may make a difference, as noticed in the $1$-dimensional case (see \citet{kratz:2014}).
We give two examples, when $\a=1.5$ (Figure~\ref{fig:qqplot 1,5}), case where the summands have no variance, 
and when $\a=3.5$ (Figure~\ref{fig:qqplot 3,5}), {\it i.e.} beyond the frame of the generalized Berry-Esseen inequality. 
In Figure~\ref{fig:qqplot 1,5}, we only have two rows as the CLT does not apply and we did not build the QQ-plot for the stable distribution. The fit looks very good for both multi-normex methods, with barely no difference of fit between the two, as shown in the zoomed part.
Whenever $\a=3.5$, we again clearly observe in Figure~\ref{fig:qqplot 3,5} an overall fit that gives the advantage to the multi-normex distributions, but with more difference, when zooming,  between $d$-normex and MRV-normex than in the previous cases $\a=1.5$ and $\a=2.3$.
\begin{figure}[H]
\centering
	\begin{tabular}{c}
Multivariate Pareto ($\a=1.5$)\\
	\includegraphics[width=165mm]{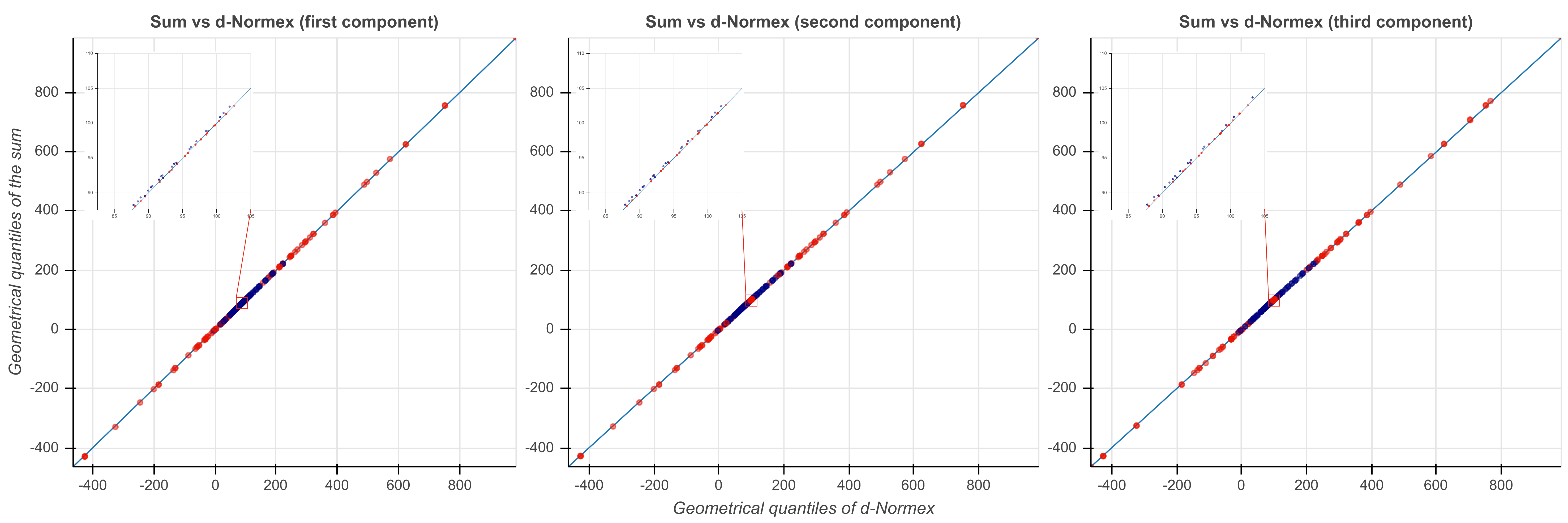}\\
	\includegraphics[width=165mm]{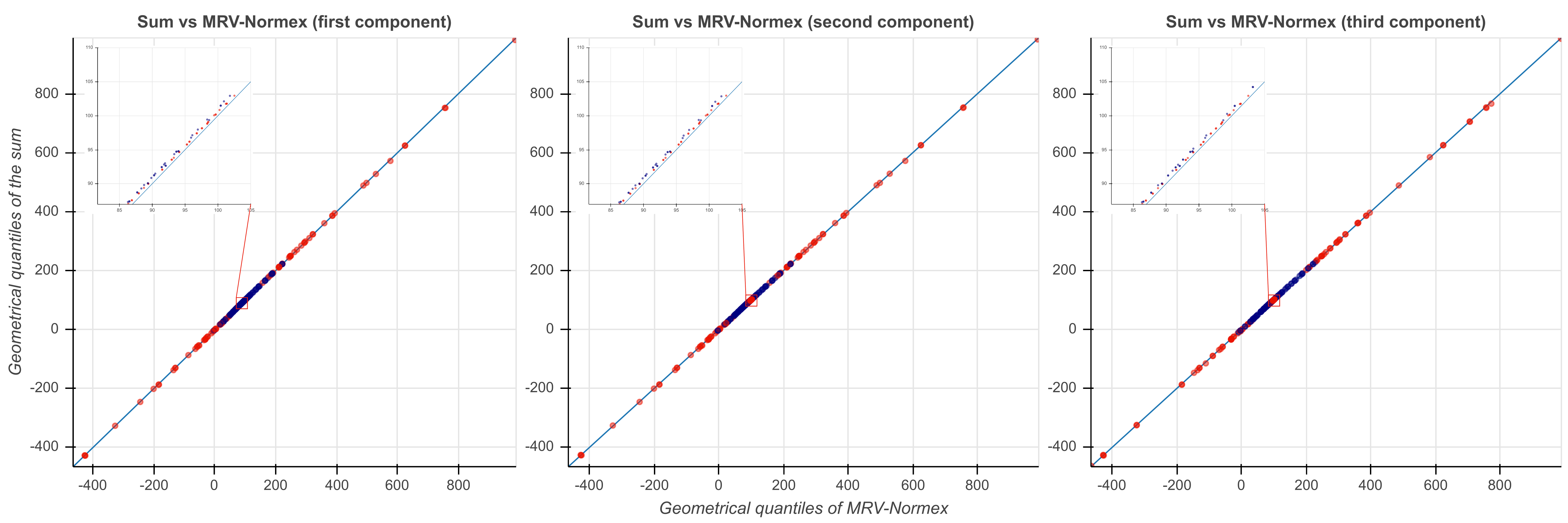}\\
     \end{tabular}
\caption{  \small  $3$-dimensional QQ-plots for the empirical distribution (sample size $= 10^7$) of the sum of $52$ iid trivariate Pareto ($\alpha = 1.5$) random vectors, with two different approximations of the sum distribution:  $d$-Normex (first row) and MRV-Normex (second row). Each column corresponds to a component (from the 1st to the 3rd). The red points on the plots correspond to extreme geometric quantiles (when the norm of the parameterized vectors is greater than $0.9$) }
\label{fig:qqplot 1,5}
\end{figure}
\begin{figure}[H]
\centering
	\begin{tabular}{c}
Multivariate Pareto ($\a=3.5$)\\
    	\includegraphics[width=165mm]{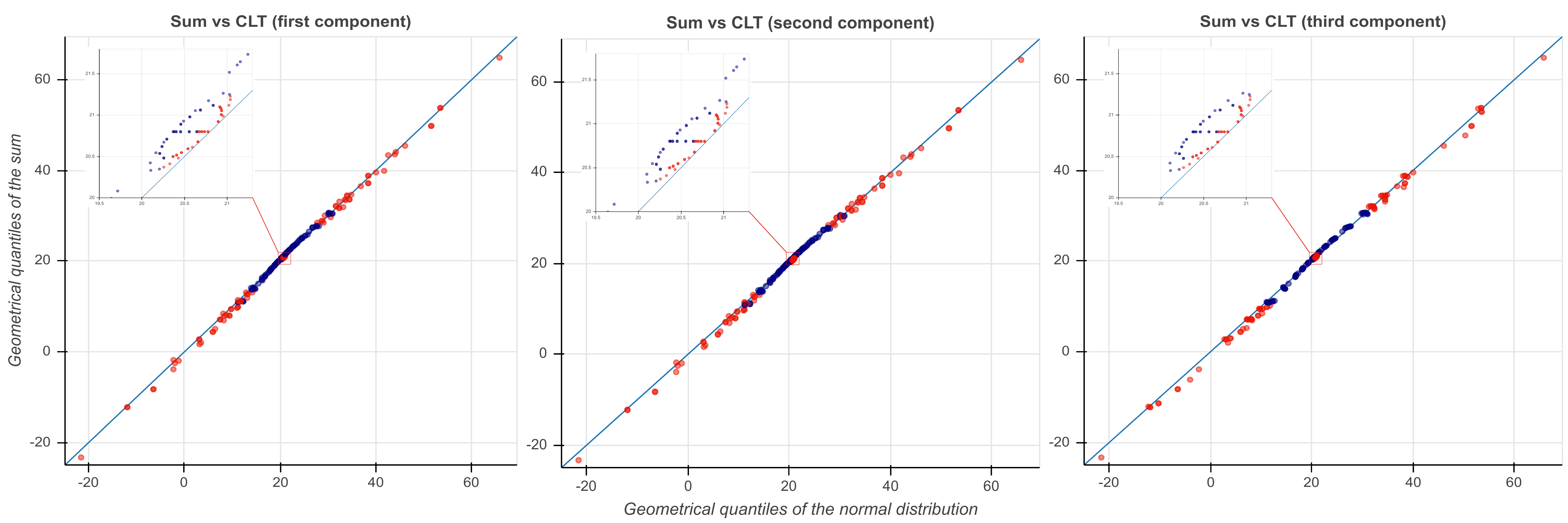} \\  	 		 
\includegraphics[width=165mm]{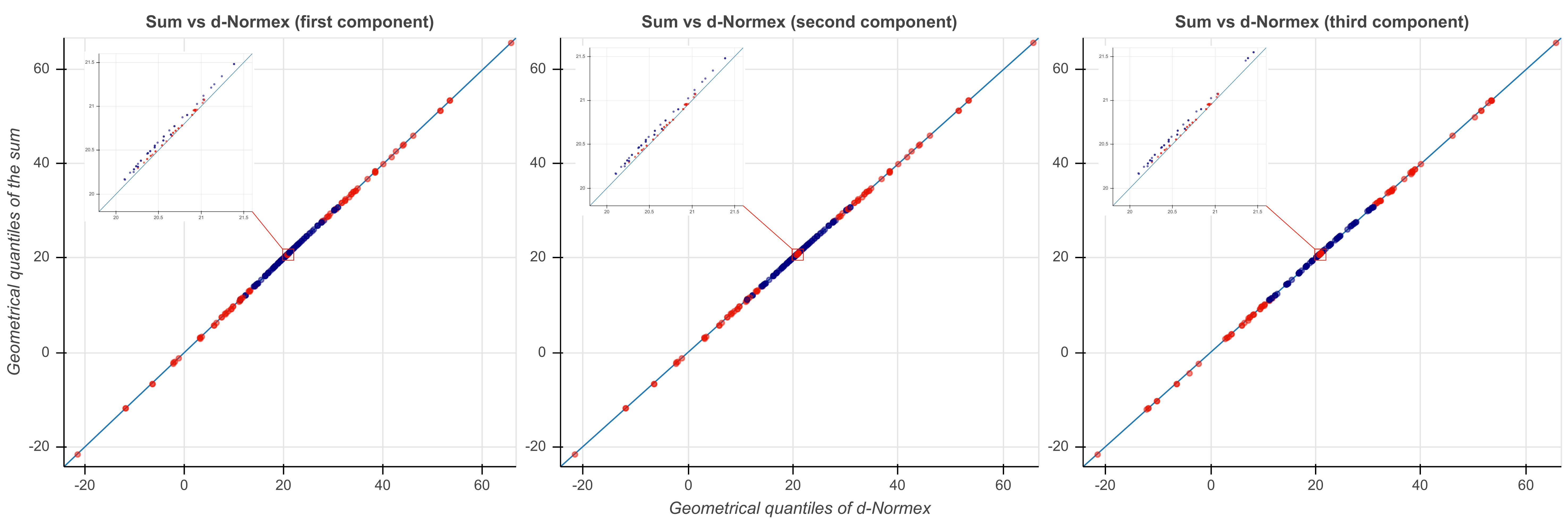} \\ 
    	\includegraphics[width=165mm]{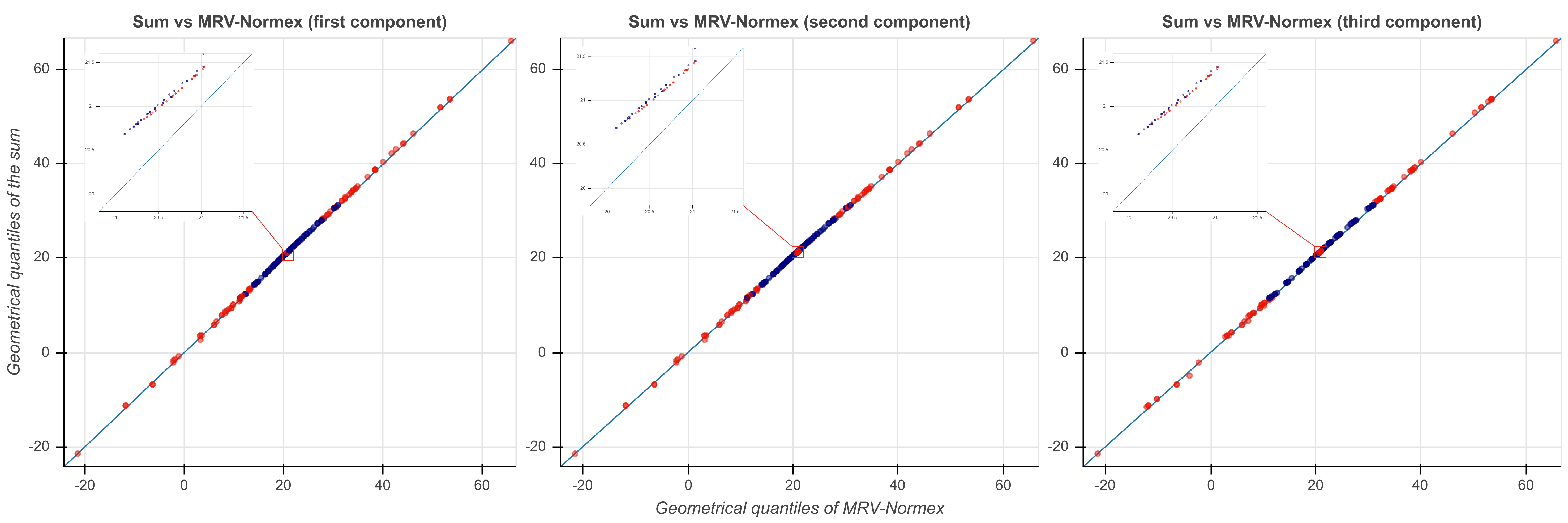} \\ 
  	\end{tabular}
\caption{ \small  $3$-dimensional QQ-plots for the empirical distribution (sample size $= 10^7$) of the sum of $52$ iid trivariate Pareto ($\alpha = 3.5$) random vectors, with three different approximations of the sum distribution: CLT (first row), $d$-Normex (second row) and MRV-Normex (third row). Each column corresponds to a component (from the 1st to the 3rd). The red points on the plots correspond to extreme geometric quantiles (when the norm of the parameterized vectors is greater than $0.9$). }
\label{fig:qqplot 3,5}
\end{figure}

\subsection{Pareto-Lomax marginal distributions with various dependence structures  - Examples~\ref{ssec:exs-MRV}}
\label{ssec:plots_dpdce}

Here, we consider the examples developed in Section~\ref{ssec:exs-MRV}, with Pareto-Lomax($\alpha$) marginal distributions, various dependence structures and the norm $\|\cdot\|_{\infty}$ instead of $\|\cdot\|_1$ as in the previous subsection. 
Figure~\ref{fig:qqplot_indep 2,3} displays QQ-plots for independent components of the vector $\mathbf{X}$, taking $\a=2.3$ to be in the half-closed interval $(2,3]$ considered in our theorems. 
Here also, the good fit of the multi-normex distributions appears clearly, in particular when comparing with the Gaussian approximation. Nevertheless, the difference between the two multi-normex distributions is more pronounced in the center, as can be observed in the zoomed part.

When the dependence structure is given via a survival Clayton copula with parameter $\theta$, we provide the QQ-plots in Figure~\ref{fig:qqplot_clayton 2,3}  assuming $\alpha\theta = 0.5$.
The previous observations hold too, with an increasing difference in the center between the two multi-normex distributions.
The case $\alpha\theta = 1$ is illustrated in Figure~\ref{fig:qqplot_clayton1 2,3}. Note that this choice of $\theta=1/\alpha$ is a standard example in the literature as it makes analytical computations much more tractable.
 
\begin{figure}[H]
\centering
\begin{tabular}{c}
Sample Scatter Plots ($\alpha = 2.3, \theta \approx 4.3$)\\
    \includegraphics[width=160mm]{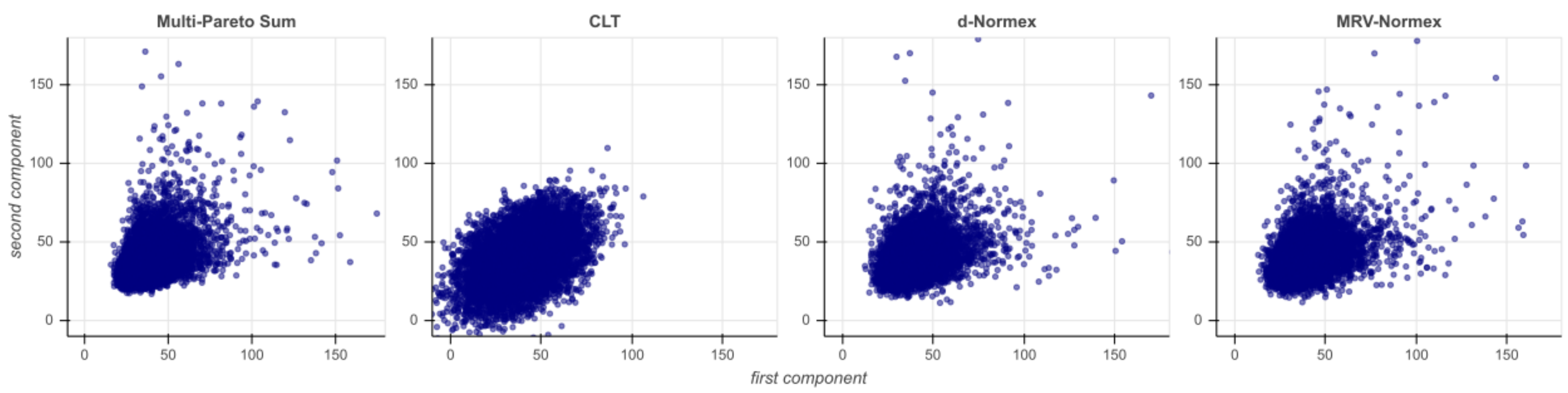} \\
    \includegraphics[width=160mm]{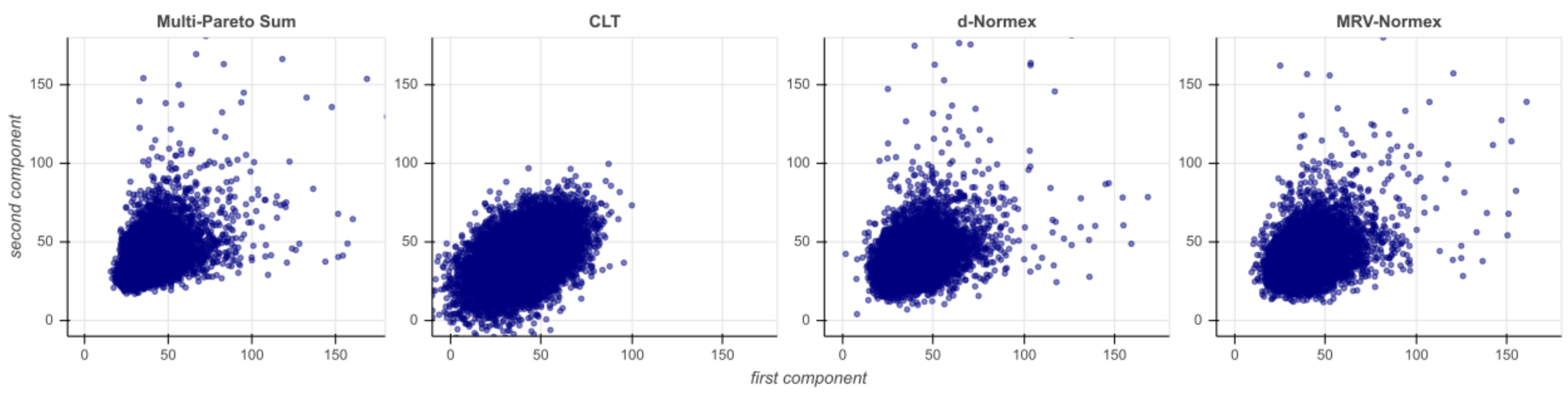} \\
     \end{tabular}
\caption{ \small Scatter plots for the Pareto-Lomax($\alpha=2.3$) marginal distributions with survival Clayton copula with parameter $\theta$ such that $\alpha \theta = 1$. The first column corresponds to the simulated sample (sample size = $10^4$) of the sum of 52 iid random vectors while the three next columns consider the different approximations of the sum distribution: CLT (2nd column), d-Normex (3rd column) and MRV-Normex (4th column). For the multi-normex approximations, two norms have been chosen, the $\|\cdot\|_{\infty}$-norm given in the 1st row, and $\|\cdot\|_{1}$-norm in the second one. }
\label{fig:scatterplot23}
\end{figure}

\begin{figure}[H]
\centering
\begin{tabular}{c}
Ranked Scatter Plots ($\alpha = 1.01, \theta \approx 0.99$)\\
    \includegraphics[width=120mm]{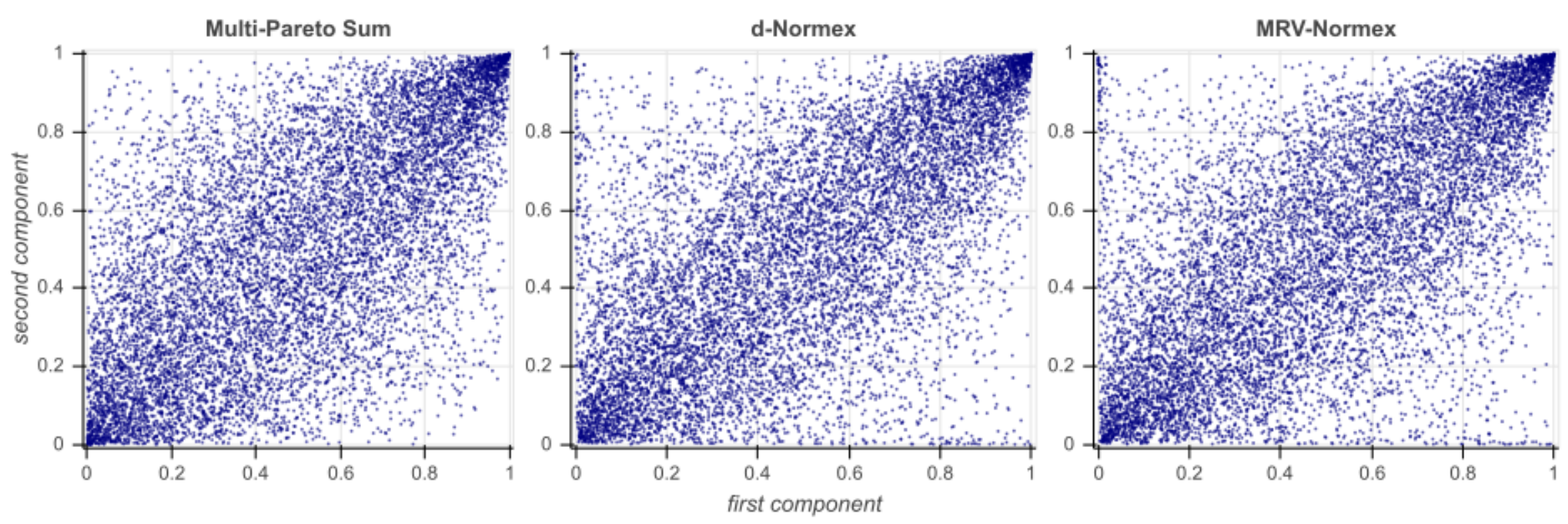} \\
    \includegraphics[width=120mm]{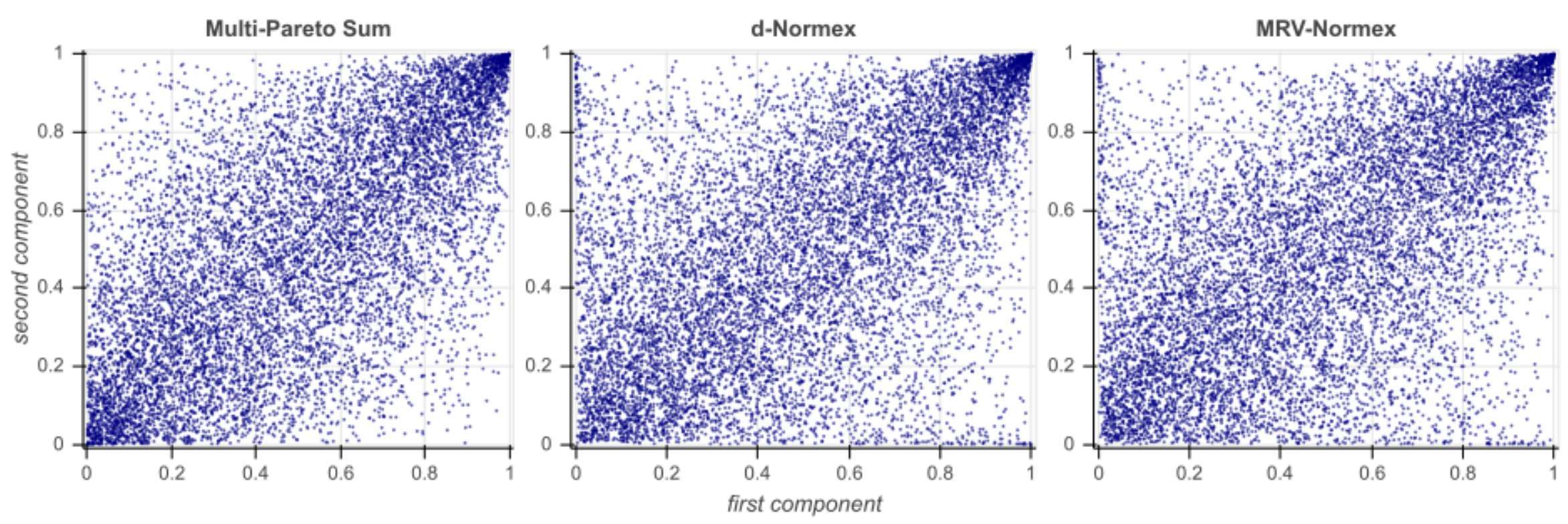} \\
     \end{tabular}
\caption{ \small  Ranked scatter plots for the Pareto-Lomax($\alpha=1.01$) marginal distributions with survival Clayton copula with parameter $\theta$ such that $\alpha \theta = 1$. The first column corresponds to the simulated sample (sample size = $10^7$) of the sum of 52 iid random vectors while the two next columns consider the multi-normex approximations: d-Normex (2nd column) and MRV-Normex (3rd column). Two norms have been chosen, the $\|\cdot\|_{\infty}$-norm given in the 1st row, and $\|\cdot\|_{1}$-norm in the second one.
}
\label{fig:ranks101}
\end{figure}

This case $\alpha\theta=1$ would also correspond to Example~\ref{ssec:exs-dNormex} if choosing the $L_1$-norm $\|\cdot\|_1$, allowing then a comparison of the results obtained respectively with the two norms. 
To easy such a numerical comparison based on norms, we reduce the dimension of the vectors to $d=2$ and provide (ranked) scatter plots, simulating $52\times10^4$ observations for $\mathbf{X}$ (hence a sample of size $10^4$ for the bivariate sum $\mathbf{S}_n$, with $n=52$). We display the (ranked) scatter plots on two rows, the first one for the case $\|\cdot\|_{\infty}$ and the second for $\|\cdot\|_{1}$.
Fixing $\a\theta=1$, we propose two choices for $\a$ to also observe the impact of the upper tail dependence in the survival Clayton copula: (i) $\a=2.3$, as previously, implying that $\theta \approx 0.43$ (weak upper tail dependence) and (ii) $\a=1.01$ to increase the value of $\theta$ (close to 1), hence the upper tail dependence, while keeping $\a>1$ (to have a finite expectation). Note that the CLT does not hold in the case (ii).
We exhibit the scatter plots for the case (i) in Figure~\ref{fig:scatterplot23}: It  is obvious that the Gaussian approximation given by the CLT does not reflect the distribution of $\mathbf{S}_n$, 
while both multi-normex scatter plots look more or less similar as that of $\mathbf{S}_n$, whatever the norm. 
It corresponds to the analytical results (Theorems~\ref{th:Result-max} and~\ref{th:MRV result}), which are independent of the norm.  
 Turning to (ii), as we have a stronger upper tail dependence, we display ranked scatter plots (rather than scatter plots). We observe in Figure~\ref{fig:ranks101} that the dependence structure looks alike among all the plots. Also for this range of $\a$, the norm does not seem to have any impact.
\begin{figure}[H]
\centering
	\begin{tabular}{c}
Independent Pareto-Lomax components  ($\a=2.3$)\\
    	\includegraphics[width=165mm]{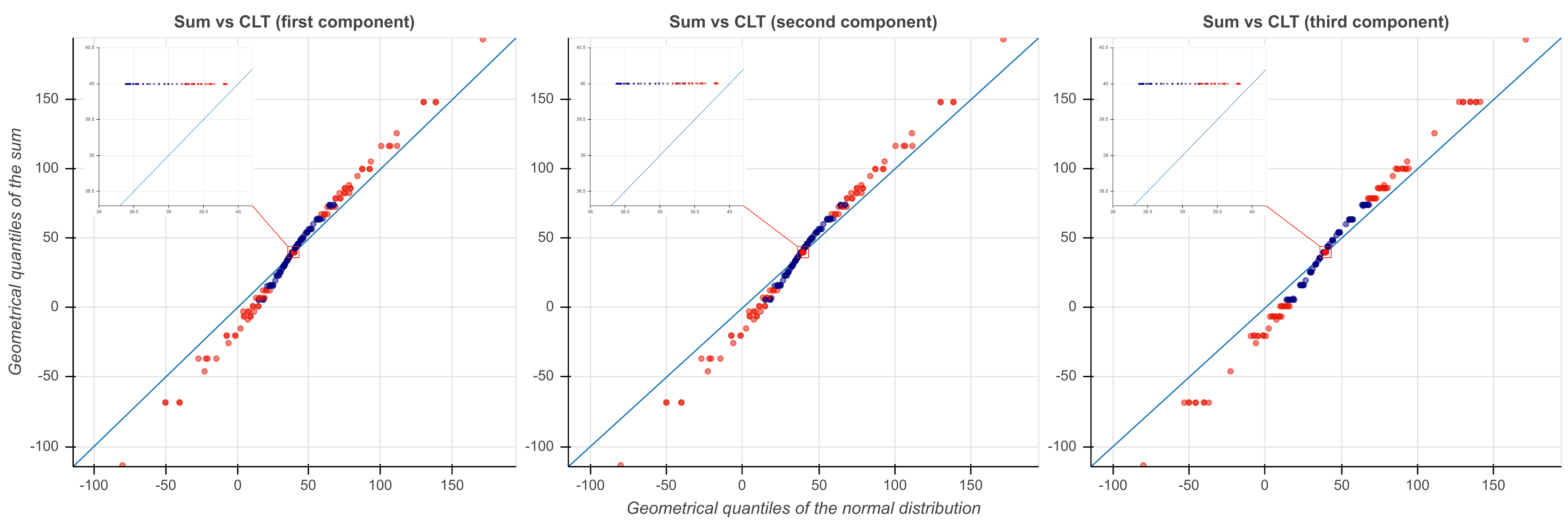} \\  	 		 
\includegraphics[width=165mm]{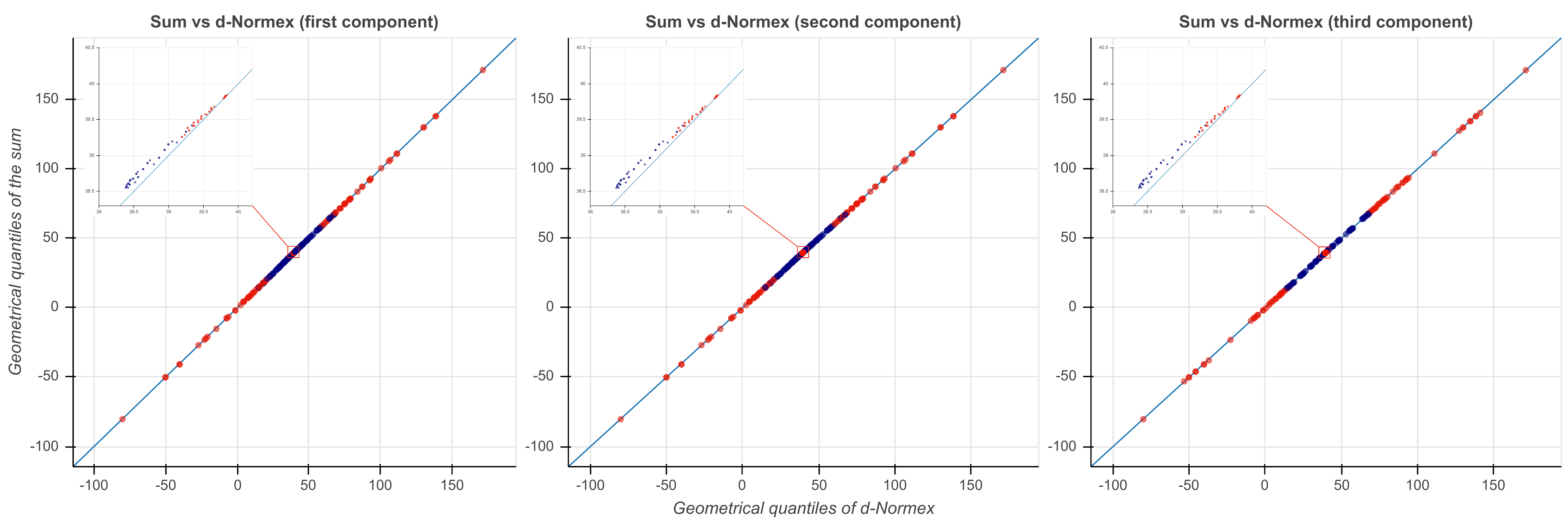} \\ 
    	\includegraphics[width=165mm]{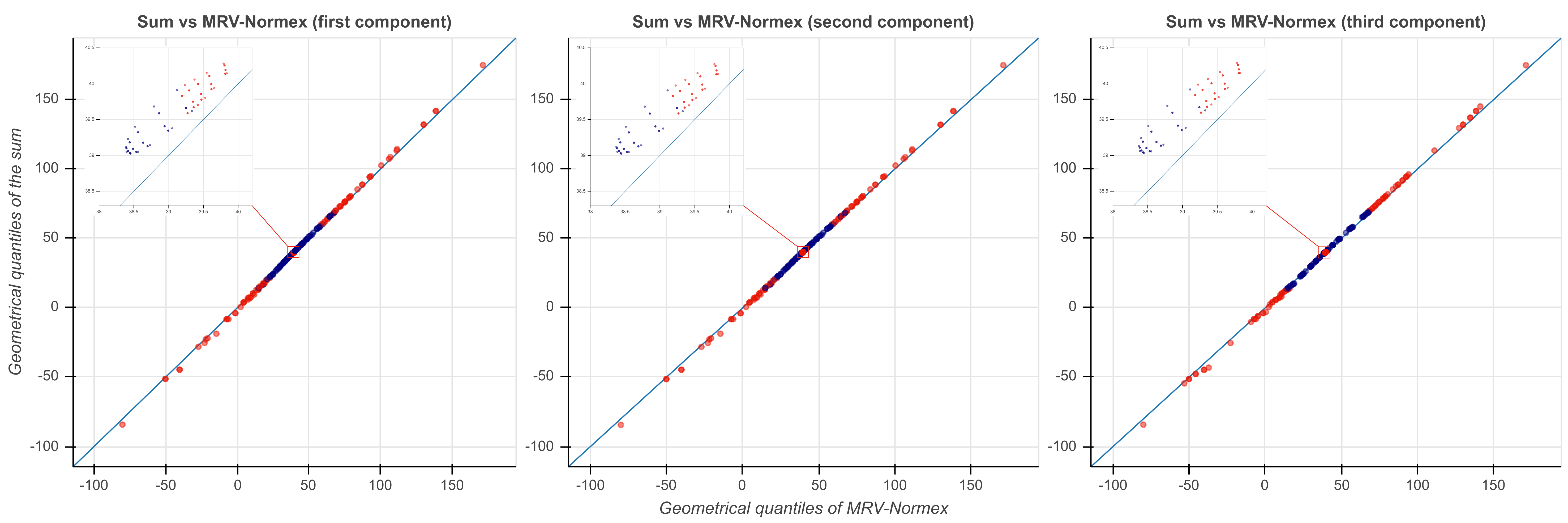} \\ 
  	\end{tabular}
\caption{ \small $3$-dimensional QQ-plots for the empirical distribution (sample size $= 10^7$) of the sum of $52$ iid random vectors with independent Pareto-Lomax components of the vectors ($\alpha = 2.3$), with three different approximations of the sum distribution:  CLT (first row), $d$-Normex (second row) and MRV-Normex (third row). Each column corresponds to a component (from the 1st to the 3rd). The red points on the plots correspond to extreme geometric quantiles (when the norm of the parameterized vectors is greater than $0.9$)  }
\label{fig:qqplot_indep 2,3}
\end{figure}

\begin{figure}[H]
    \centering
    	\begin{tabular}{c}
    	Pareto-Lomax with various Clayton copulas ($\a=2.3$) \\
        \begin{subfigure}[t]{0.5\textwidth}
            \centering
            \begin{tabular}{c}
            \includegraphics[width=75mm]{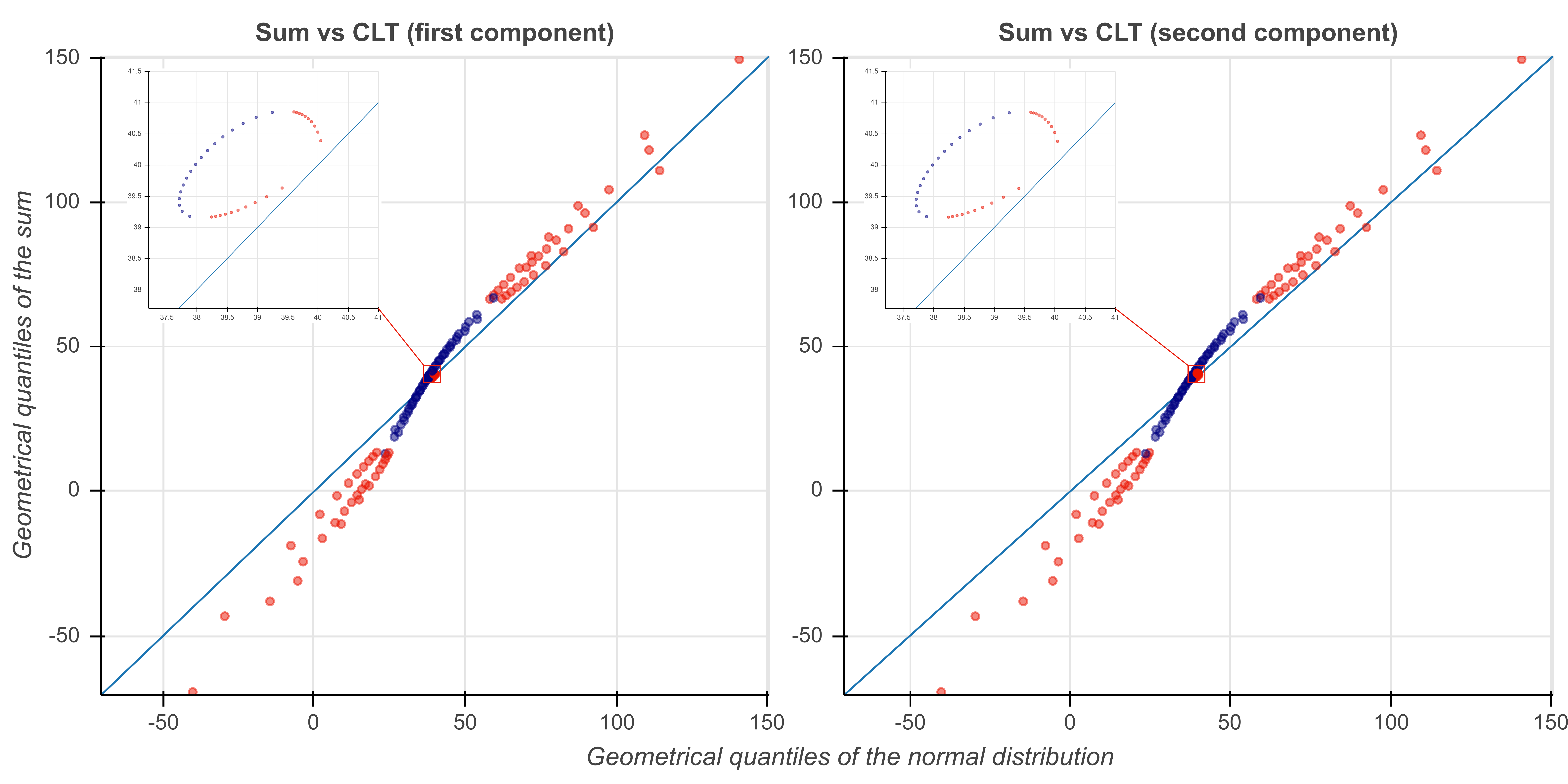} \\  		 
            \includegraphics[width=75mm]{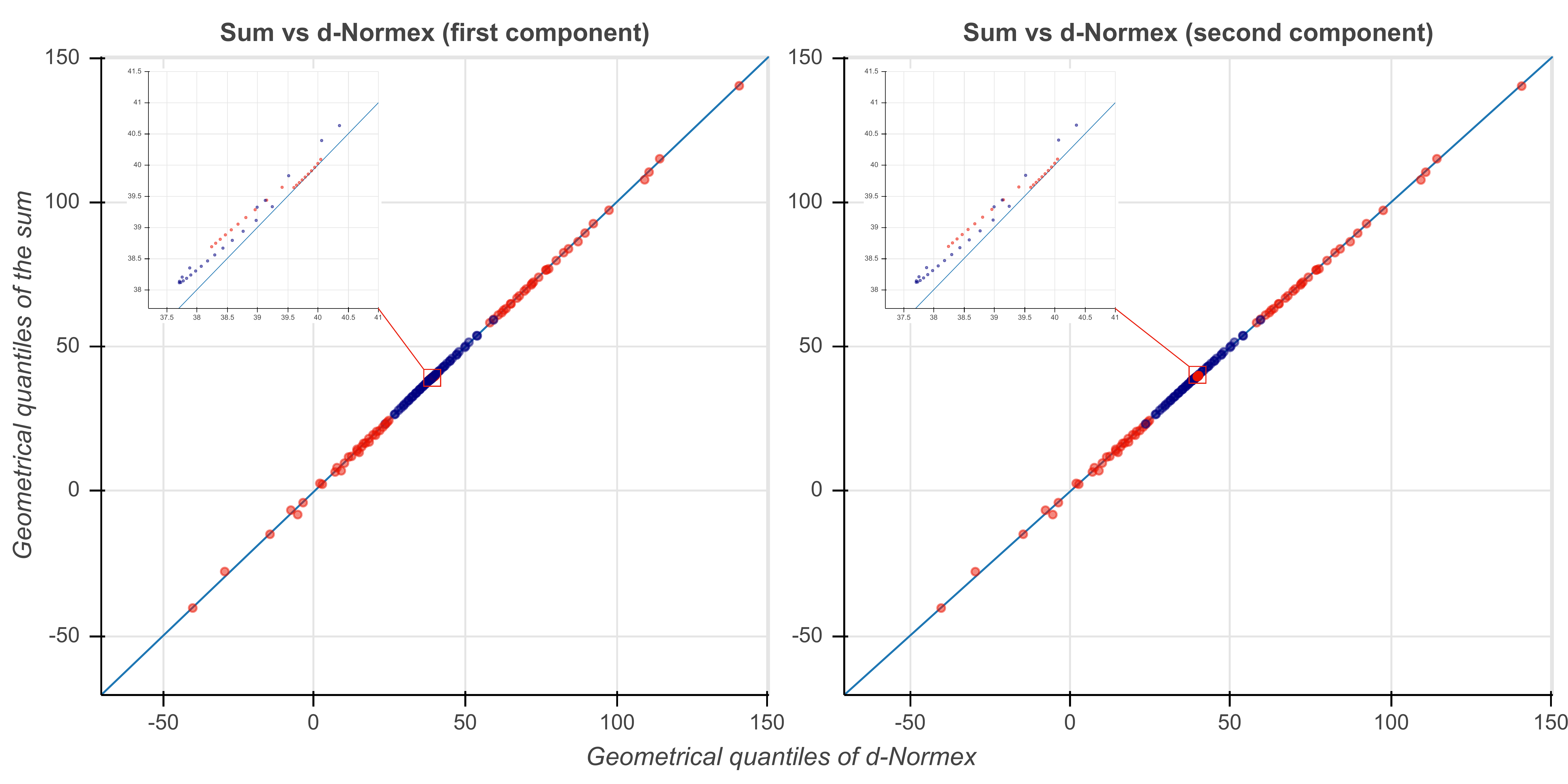} \\ 
        	\includegraphics[width=75mm]{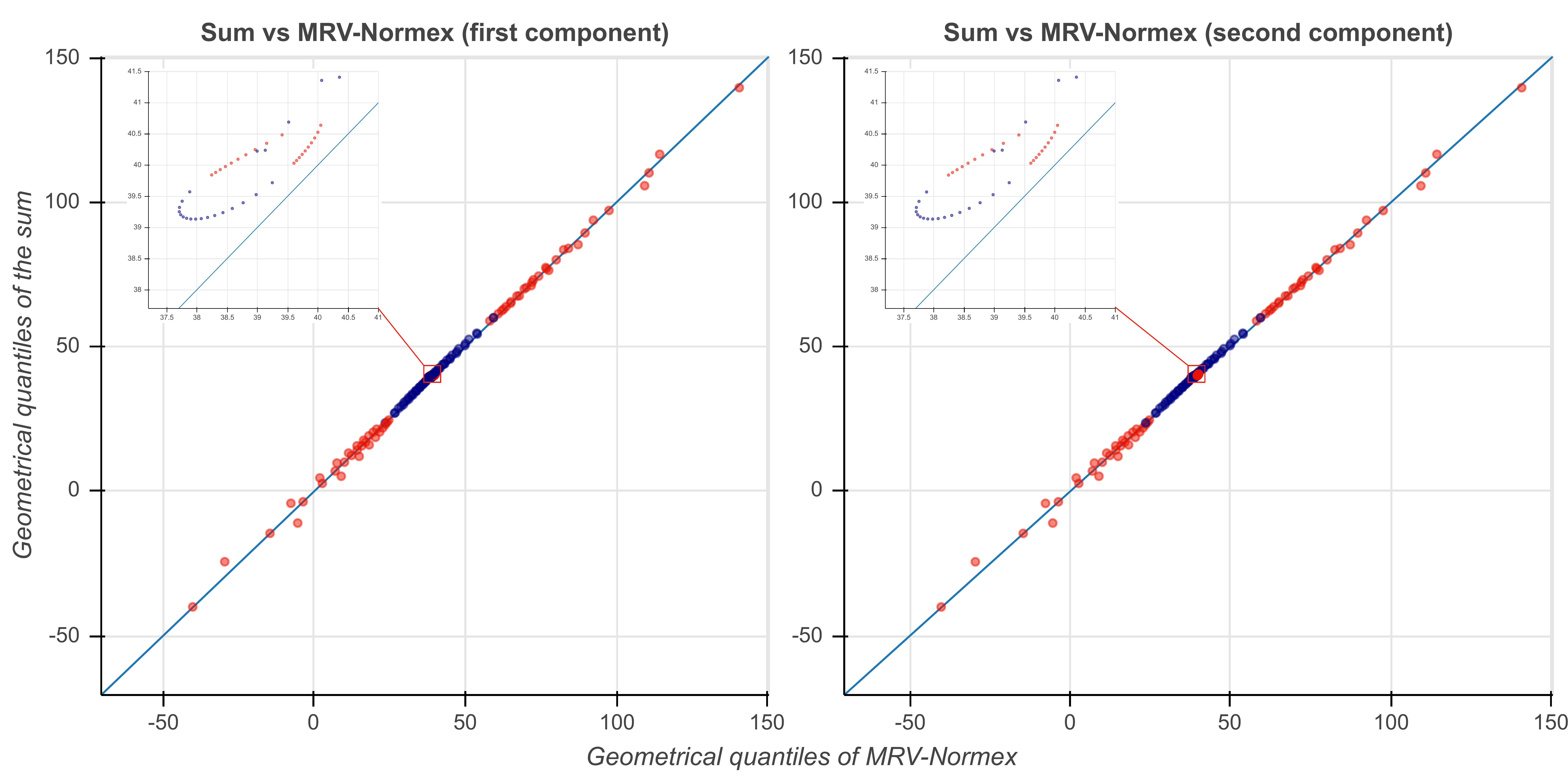} \\ 
        	\end{tabular}
        \caption{$\alpha\theta = 0.5, \theta \approx 0.22$}        
      \label{fig:qqplot_clayton 2,3}
            \end{subfigure}%

    ~ 
            \begin{subfigure}[t]{0.5\textwidth}
        \centering
            	\begin{tabular}{c}
    	\includegraphics[width=75mm]{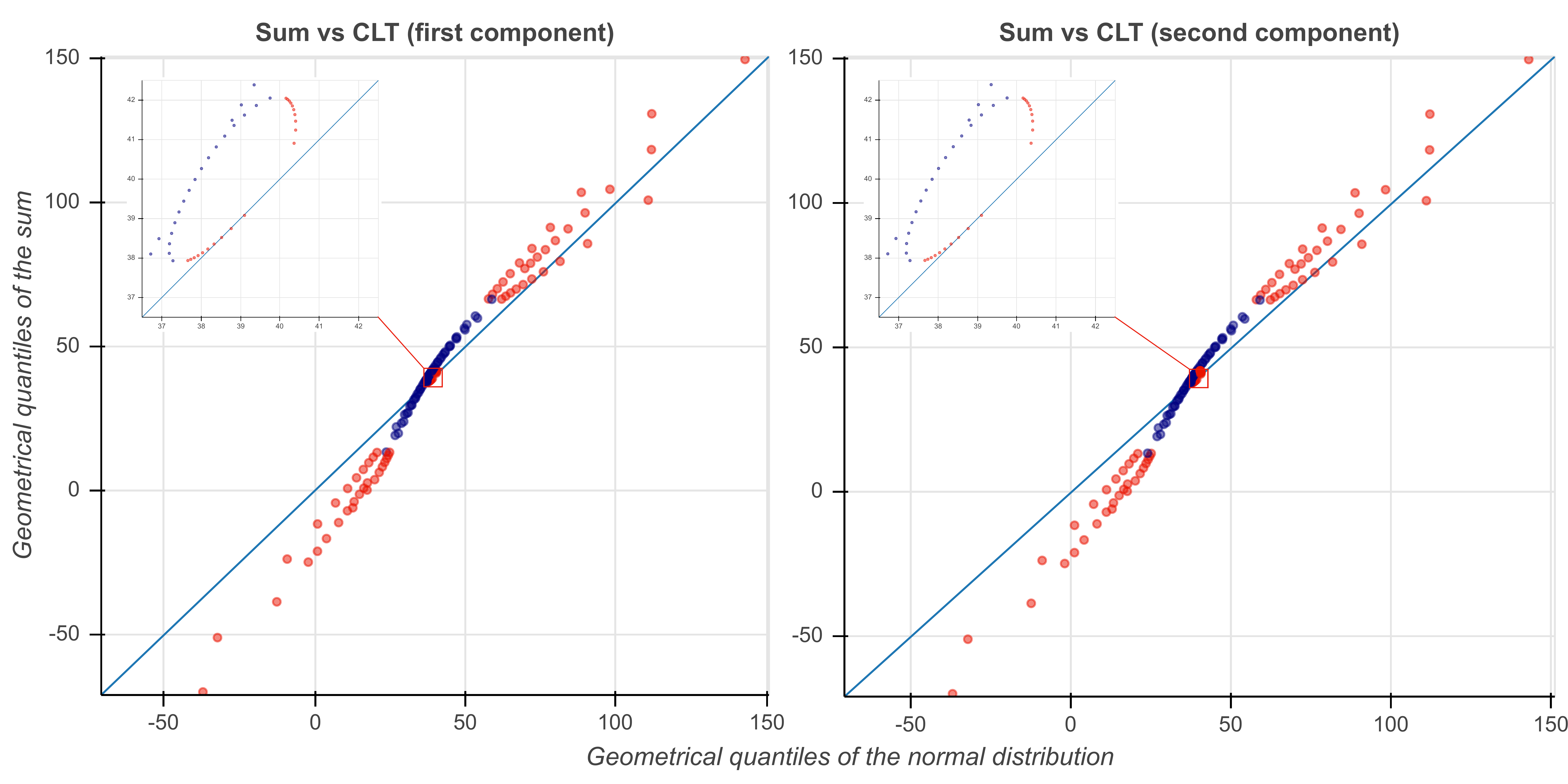} \\  	 		
        \includegraphics[width=75mm]{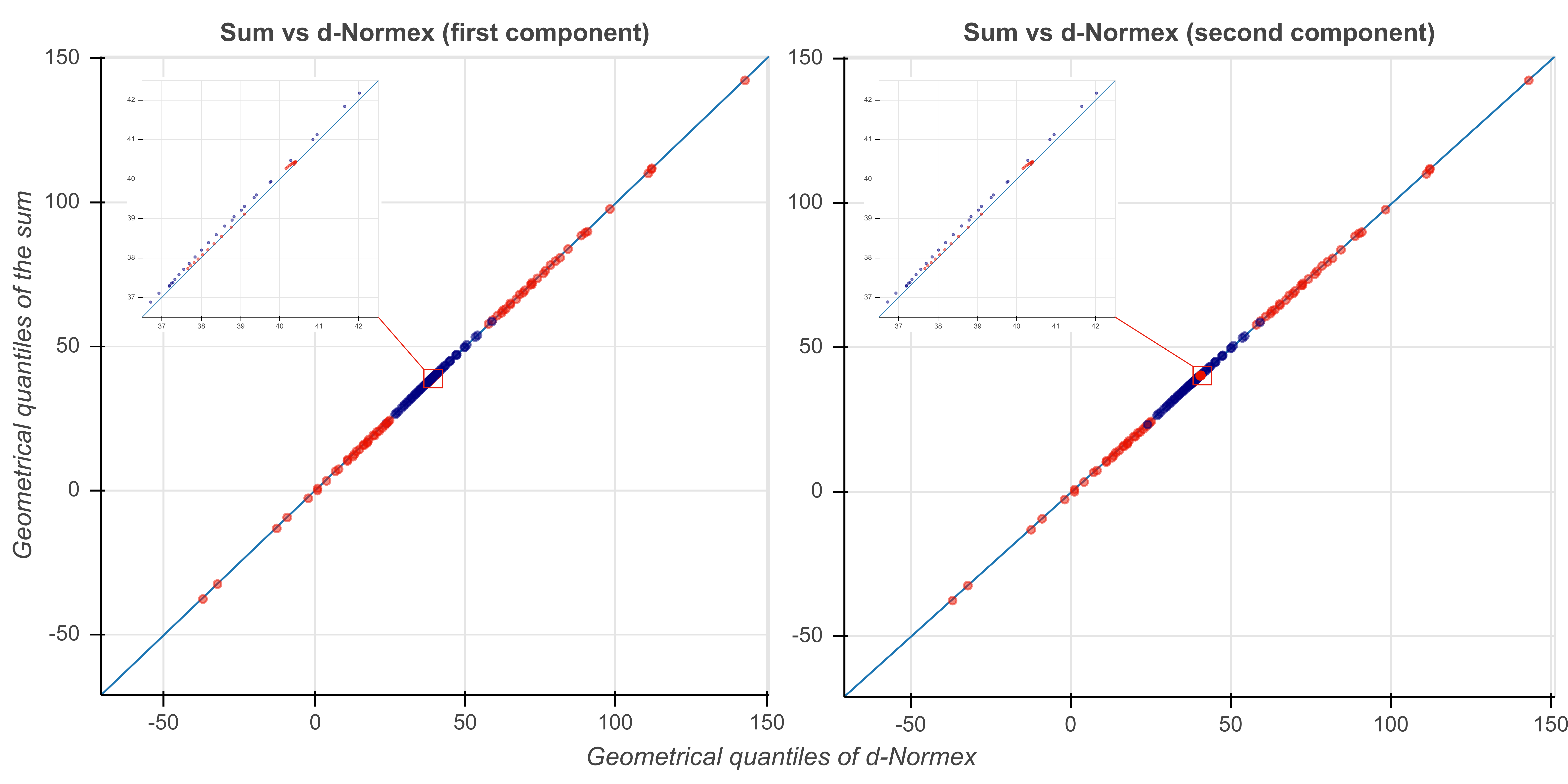} \\ 
    	\includegraphics[width=75mm]{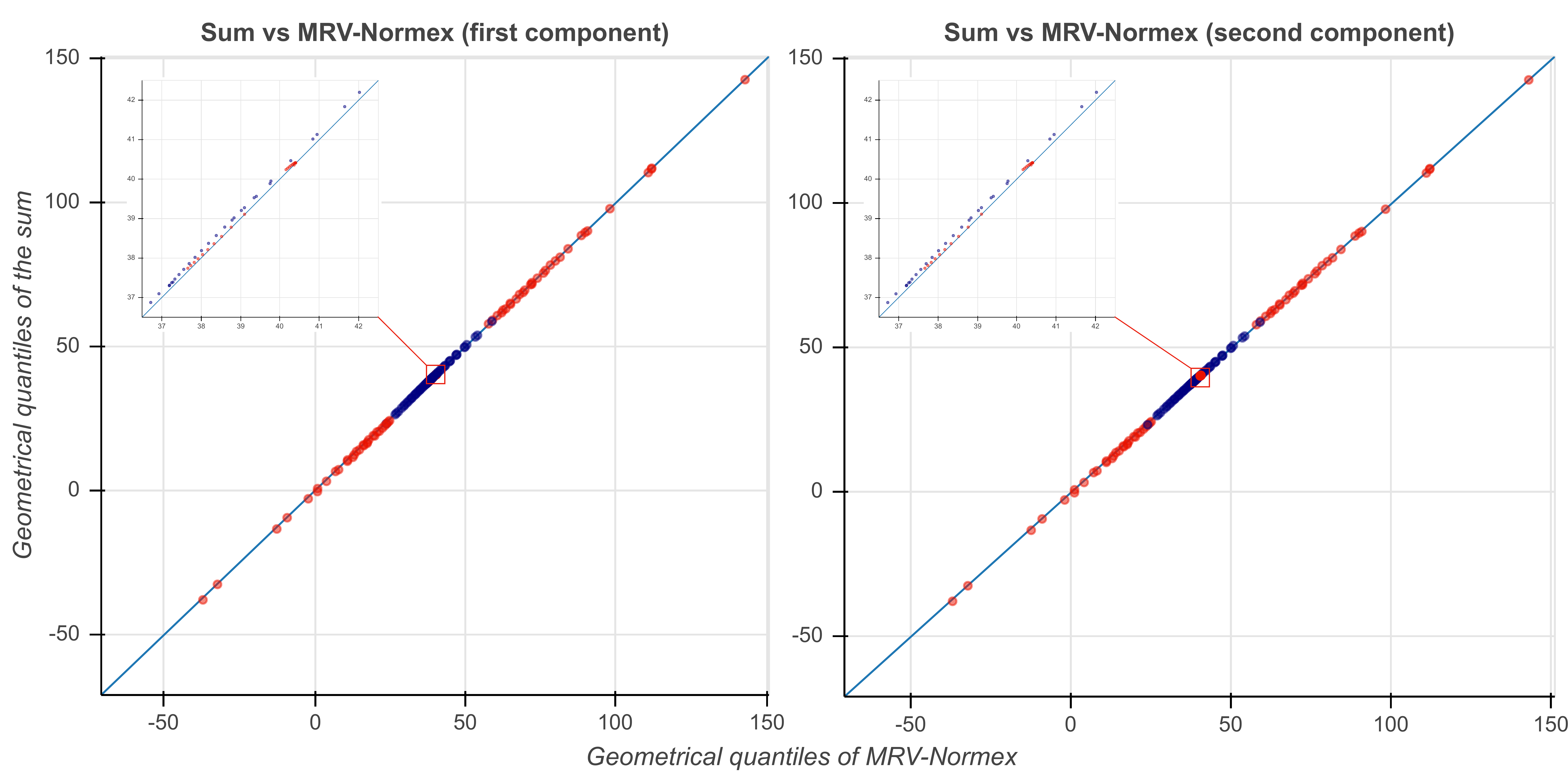} \\
        	       \end{tabular}
        \caption{$\alpha\theta = 1, \theta \approx 0.43$}      
      \label{fig:qqplot_clayton1 2,3}
            \end{subfigure}%
        \end{tabular}
    \caption{ \small  $2$-dimensional QQ-plots for the empirical distribution (sample size $= 10^7$) of the sum of $52$ iid random vectors with Pareto-Lomax marginal distributions ($\alpha = 2.3$) joint by a Clayton copula with parameter $\theta$ such that $\a\theta=0.5$ for the case (a), and $\a\theta=1$ for the  case (b). The three rows correspond to the approximations of the sum distribution: CLT (first row), $d$-Normex (second row) and MRV-Normex (third row). For each case (a) and (b), the two columns correspond to the 1st and the  2nd components. The red points on the plots correspond to extreme geometric quantiles (when the norm of the parameterized vectors is greater than $0.9$).}
\end{figure}

\section{Conclusion}\label{sec:concl}

The purpose of this study was to build a sharp approximation of the whole
distribution of the sum of iid random vectors under the presence of heavy tails. It has been done by extending the normex approach from a univariate to a multivariate framework, combining mean and extreme behaviors.  We proposed two possible multi-normex distributions, named $d$-Normex and MRV-Normex. Both rely on the Gaussian distribution for describing the mean behavior, via the CLT, while the difference between the two versions comes from using the EV theorem or the exact distribution for the maximum. The main theorems establish the rate of convergence of each version of the multi-normex distributions towards the distribution of the sum. This is done analytically whenever the shape parameter $\a$ of the tail of the marginal distribution belongs to the interval $(2,3]$, making the comparison with the generalized Berry-Esseen inequality relevant. For the MRV-Normex, second order regular variation conditions are needed to obtain the main theorem. These conditions are discussed w.r.t. the Extreme Value Theory literature.  Numerical comparisons are developed for any value of $\a$, for both multi-normex distributions, considering examples with different structure of dependence for the random vectors. Illustrations are made through multidimensional QQ-plots based on geometrical quantiles.

We focused on the case of heavy tailed random vectors, as it is of most interest in the risk literature.  Nevertheless, this method could be extended to light tails vectors (whenever $1/\a=0$)  as the rate of convergence for the EV theorem is also known in such a case. It would then require to introduce a specific metric to evaluate the error for the whole distribution taking into account the impact of extremes. 

Moreover, the MRV-Normex approach has been developed conditioning on the norm of the maximum. It could be done conditioning on the maximum itself (a vector). 
To widen the applicability of the multi-normex methods, simple approximations of truncated moments could also be suggested (as numerical ones or evaluating them using e.g. a Pareto approximation).
Finally, generalization of multi-normex distributions will be studied when introducing dependence between random vectors, then considering random processes. We intend to explore such topics in the near future.

\section*{Acknowledgments}

Evgeny Prokopenko acknowledges the support received from the National Research Agency of the French government through the program "Investment for the future" (ANR-16-IDEX-0008 CY Initiative) during his postdoctoral fellowship at ESSEC CREAR.
Evgeny Prokopoenko is  supported by the Mathematical Center in Akademgorodok under the agreement No. 075-15-2019-1675 with the Ministry of Science and Higher Education of the Russian Federation.

\printbibliography

\newpage

\begin{appendix}

\section{Regular Variation and related notions}\label{A:secRV}

Let us recall the definitions and properties of regular variation in both univariate and multivariate cases (see \citet{bingham:goldie:teugels:1989,geluk:dehaan:1987,dehaan:1970,dehaan:ferreira:2006,resnick:2002,resnickbook:2008,seneta:1976}), and point out some useful relations between each other.

\paragraph{Univariate distributions -}
\begin{defi}[Regular variation in one dimension; see \citet{bingham:goldie:teugels:1989}]\label{def:rv}
A rv $X$ with cdf $F$ has a regularly varying (right) tail with index $\alpha > 0$ if $\oF = 1-F \in \RV_{-\alpha}$ (or, by abuse of notation, $X \in \RV_{-\alpha}$), {\it i.e.}
$$
\lim_{t\to\infty}  \frac{\oF(tx)}{\oF(t)} = x^{-\a}, \quad x>0.
$$
Alternatively, we  say that $X$ has a regularly varying tail if there exists a function $b:\R_{+}\to\R_{+}$ with $b(t) \uparrow \infty$ as $t\to \infty$ such that
\begin{align*}
& \lim_{t\to\infty} t\,\p{X>b(t) x} = x^{-\a} .  
\end{align*}
We write $\oF\in \RV_{-\alpha}(b)$ or, by abuse of notation, $X \in \RV_{-\alpha}(b)$.
\end{defi}
As a consequence, $b(\cdot)\in \RV_{1/\a}$ and a natural choice is $b(t)=(1/\overline{F})^{\la}(t)$. \\[1ex]

\begin{defi}[Extended Regular variation; see e.g. \citet{dehaan:ferreira:2006}]\label{def:Erv}
A measurable function $f$ is of {\it extended regular variation}, denoted by $f\in \ERV(\gamma)$, if for some $\gamma\in\R$ and positive function $a$, 
$$
\lim_{t\to\infty} \frac{f(tx) - f(t)}{a(t)}= \frac{x^\gamma - 1}{\gamma}
$$
for all $x>0$.
\end{defi}

When looking at the speed of convergence to the limits given in the definition of (E)RV, comes the notion of "second-order (extended) regular variation".

\begin{defi}[Second order regular variation; see e.g. \citet{dehaan:resnick:1993}]
A rv $X$ with cdf $F$ such that $\bar{F} \in \mathcal{R} V_{-\alpha}(b)$ with $\alpha>0,$ possesses second order regular variation with parameter $\rho \leq 0$, denoted by $X\in 2\RV(-\alpha,\rho)$, 
if there exists a function $A(t) \underset{t \rightarrow \infty}{\rightarrow} 0$ that is ultimately of constant sign, $|A(t)| \in \mathcal{R} V_{\rho}$ with $\rho \leq 0$ and $c \neq 0$ such that
\begin{equation}\label{def:2RV}
\frac{t \bar{F}(b(t) x)-x^{-\alpha}}{A(b(t))} \,\underset{t \rightarrow \infty}{\longrightarrow}\, c x^{-\alpha} \frac{x^{\rho}-1}{\rho}, \quad x>0
\end{equation}
The right hand side of Eq. \eqref{def:2RV} is interpreted as the function $c x^{-\alpha} \log (x)$ when $\rho=0$.
\end{defi}

\begin{defi}[Second order extended regular variation; see e.g. \citet{dehaan:stadtmueller:1996}]
A measurable function $f$  is said to be of second-order extended regular variation (with second-order parameter $\rho\le 0$), denoted by $f\in2\ERV(\gamma,\rho)$, if there exist positive functions $a$ and $A$ with $\displaystyle \lim_{t\to\infty} A(t)=0$ such that 
$$
\frac{\frac{f(tx)-f(t)}{a(t)} -\frac{x^\gamma -1}{\gamma}}{A(t)} \,\underset{t \rightarrow \infty}{\longrightarrow}\, H(x)
$$
for all $x>0$, where the limit function $H$ is not a multiple of $(x^\gamma-1)/\gamma$. 
\end{defi}

{\sc An equivalent representation of 2RV functions and "potter bounds"}

It appears that a $2\RV$  function has a specific form.
\begin{lem}[Lemma 3 in \citet{hua:joe:2011b}]\label{lem:representation_2RV}
 Suppose $g \in 2 \RV_{\alpha, \rho}, \, \alpha\in \mathbb{R}, \rho<0,$ then we can write $g=k t^{\alpha} (1 + l(t)), t>0$ with
some constant $k \neq 0$ and regularly varying function $l(t) \in \mathrm{RV}_{\rho}.$
\end{lem}
In other words, Lemma~\ref{lem:representation_2RV} says that a slowly varying function $L(\cdot)$ can be a $2\RV_{0, \rho}$ function with parameters $0, \rho < 0$, only in the case $\lim_{t \ra \infty} L(t) \neq 0, \pm \infty$. In the case $\rho = 0$, the function $L$ belongs to so called $\Pi$-class and has the representation (see theorem B.2.12 in \citet{dehaan:ferreira:2006})
$$
L(t)=c_{1}+c_{2} a_{1}(t)+\int_{1}^{t} \frac{a_{2}(s)}{s} \mathrm{d}s,
$$
where $c_1, c_2 \in \mathbb{R}$ are constants and $a_1(\cdot), a_2(\cdot)$ are slowly varying functions such that \\
$\displaystyle a_1(t) \underset{t\to\infty}{\sim} a_2(t)$.
The following lemma will be useful to prove Proposition~\ref{cvRate-Frechet}.
\begin{lem}[Proposition 4 in  \citet{hua:joe:2011b}]\label{lemma:potterbound}
Suppose $g \in 2 \mathrm{RV}_{\alpha, \rho}$ with $\alpha \in \mathbb{R}$ and $\rho<0,$ then for any $\epsilon, \delta>0,$ there exists $t_{0}=t_{0}(\epsilon, \delta)$ such that for all $t, t x \geq t_{0}$ and $x>0$
\begin{equation}\label{potterbounds}
\left|\frac{g(t x) / g(t)-x^{\alpha}}{a(t)}-x^{\alpha} \frac{x^{\rho}-1}{\rho}\right| \leq \epsilon \, \max \left(x^{\alpha+\rho+\delta}, x^{\alpha+\rho-\delta}\right),
\end{equation}
where $a(t)=-\rho[\ell(t)] /(1+ \ell(t))$ with $g(t)=k t^{\alpha}(1+ \ell(t)), 0<k<+\infty$ and $l(\cdot) \in \mathrm{RV}_{\rho}.$
\end{lem}
In fact, Hua and Joe use Lemmas~\ref{lem:representation_2RV} and \ref{lemma:potterbound} for $g(t) = \Bar{F}(t)$,  and, for this reason, proved the lemmas in the case $\a < 0$. But their proof can be repeated line by line for an arbitrary $\alpha \in \mathbb{R}$.  
We believe that Lemma~\ref{lemma:potterbound} should also hold true when $\rho = 0$, but we will not  develop further on this question.

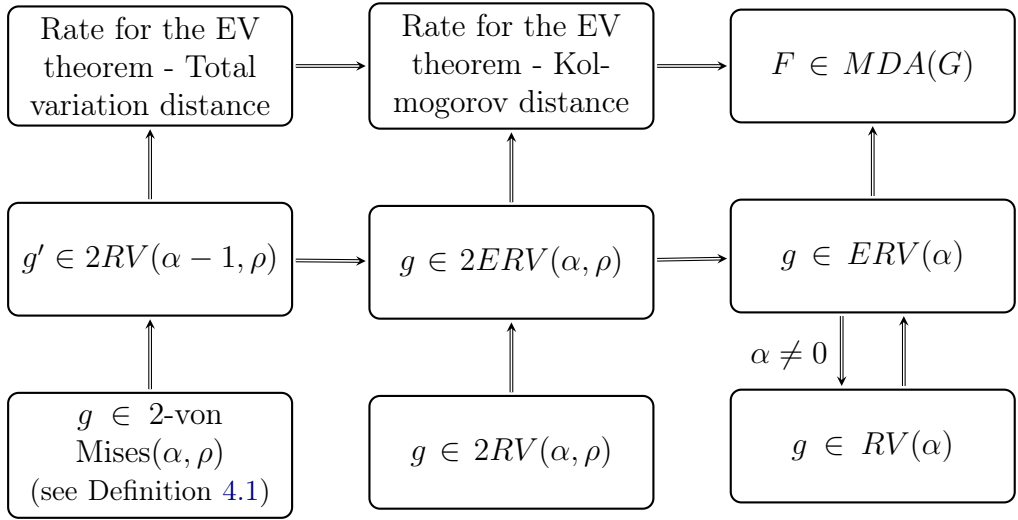
\begin{figure}
\centering
\begin{tikzpicture}  
  \node[block] (11) {Rate for the EV theorem - Total variation distance};  
  \node[block,right=of 11] (12) {Rate for the EV theorem - Kolmogorov distance};  
  \node[block,right=of 12] (13) {$F \in MDA(G)$};  
  \node[block,below=of 11] (21) {$g' \in 2RV(\a-1,\rho)$};
  \node[block,below=of 12] (22) {$g \in 2ERV(\a,\rho)$}; 
  \node[block,below=of 13] (23) {$g \in ERV(\a)$}; 
  \node[block,below=of 21] (31) {$g \in 2$-von Mises$(\a,\rho)$ {\small(see~Definition~\ref{def:vonMis})}};
  \node[block,below=of 22] (32) {$g \in 2RV(\a,\rho)$}; 
  \node[block,below=of 23] (33) {$g \in RV(\a)$}; 
  
  \draw[>=stealth, double,->, shorten >=1pt, shorten <=1pt]  (11)-- (12);  
  \draw[>=stealth, double,->, shorten >=1pt, shorten <=1pt] (12)-- (13);  
  \draw[>=stealth, double,->, shorten >=1pt, shorten <=1pt] (21)-- (22); 
  \draw[>=stealth, double,->, shorten >=1pt, shorten <=1pt] (22)-- (23); 
  \draw[>=stealth, double,->, shorten >=1pt, shorten <=1pt] (21)-- (11); 
  \draw[>=stealth, double,->, shorten >=1pt, shorten <=1pt] (22)-- (12);  
  \draw[>=stealth, double,->, shorten >=1pt, shorten <=1pt] (23)-- (13); 
  \draw[>=stealth, double,->, shorten >=1pt, shorten <=1pt] (31) -- (21); 
  \draw[>=stealth, double,->, shorten >=1pt, shorten <=1pt] (32) -- (22);  
 \draw[transform canvas={xshift=1em}, >=stealth, double,->, shorten >=1pt, shorten <=1pt] (33) -- (23);  
 \draw[transform canvas={xshift=-1em}, >=stealth, double,->, shorten >=1pt, shorten <=1pt]  (23) -- node [below of=(23), xshift=-0.7cm, yshift=0.95cm]{$\alpha \neq 0$}  (33);
\end{tikzpicture}  
\caption{\sf\small Relations between RV related notions, where $F$ is the cdf of a real rv $X$ and the function $g$ is defined by $g(x):=\left(-1/\log F(x)\right)^{\leftarrow}$}\label{fig:compareRVnotions}
\end{figure}

Finally, let us end this paragraph by representing  the relations between the various notions of (related) regular variation in Figure~\ref{fig:compareRVnotions}.

\paragraph{Multivariate Distributions -}
Multivariate regular variation appears as a natural extension of Definition \ref{def:rv}, when studying jointly heavy-tailed random variables. The notion of vague convergence of measures is used for convergence of measures on the non-negative Euclidean orthant $\R_{+}^{d}$ and its subsets; see \citet{resnickbook:2007} for further details.

\begin{defi}[Multivariate regular variation, \citet{resnickbook:2007}]\label{def:mrv}
 Suppose $\bX=(X_1,\ldots,X_d)$ is a random vector in  $[0,\infty)^d$. Then $\bX$ is {multivariate regularly varying}  if there exist $b(t) \uparrow \infty$ with $b(\cdot) \in \RV_{1/\alpha}$, $\a>0$, and a Radon measure $\nu\neq 0$ such that, as $t\to\infty$, 
$$
t\,\p{\frac{\bX}{b(t)} \in \;\cdot\; } \, \conv  \, \nu(\cdot) \quad \text{on } \;\; \M_{+}(E),
$$
where $\conv $ denotes vague convergence of measures on $\M_{+}(E)$, the class of Radon measures on Borel subsets of $E:=[0,\infty]^{d}\setminus \{\bzero\}$. We write  $\bX\in \MRV_{-\alpha}(b,\nu)$.
\end{defi}
The measure $\nu(\cdot)$ has a scaling property for relatively compact $A\subset E$, given by
$$
\nu(kA) = k^{-\alpha}\nu(A), \quad k>0\,.
$$
As already mentioned in the paper, an equivalent definition is through the pseudo-polar representation: $\bX\in \MRV_{-\alpha}$  if there exists a random vector $\bTheta$ with values in the unit sphere ${\cal S}_{d-1}$ in $\R^d$ w.r.t. the norm $\|\cdot\|$, such that, for all $t >0$,
$$
\frac{\p{\|\bX\|>t\,u,\;\bX/\|\bX\|\in \cdot }}{\p{\|\bX\|>u}} \, \conv  \, t^{-\a}\,\p{\bTheta \in \cdot}\quad \text{as} \;u\to \infty.
$$
Note that a characterization of MRV for $\bX$ is that any linear combination $\langle\boldsymbol{t},\bX\rangle$, $\boldsymbol{t}\in\R^d$, is regularly varying, $\langle\cdot,\cdot\rangle$ denoting the usual inner product in $\R^d$, namely (see \citet{basrak:davis:mikosch:2002b}):\\
There exist $\a>0$ and a slowly varying function $L$ such that, for all $\boldsymbol{x}$,
$$
\lim_{u\to\infty} \frac{\p{\langle\boldsymbol{x},\bX\rangle \geq u}}{u^{-\a}L(u)} = w(\boldsymbol{x}) \;\text{exists}\quad\text{and}\quad \exists\, \boldsymbol{x}_0\neq \bzero \;\text{s.t.}\; w(\boldsymbol{x}_0)>0.
$$

As for the dimension 1, we can define the notion of second order multivariate regular variation, namely:
\begin{defi}[Second order multivariate regular variation; see \citet{resnick:2002}]\label{def:2mrv} 
Suppose $\bX\in \MRV_{-\alpha}(b,\nu)$ and there exists $\displaystyle A(t)\underset{t\to\infty}{\to} 0$ that is ultimately of constant sign with $|A(\cdot)| \in \RV_{\rho},\, \rho \le 0$, such that
\begin{align}\label{eq:2mrv}
\frac{t\,\p{\frac{\bX}{b(t)} \in [\bzero,\boldsymbol{x}]^c } - \nu([\bzero,\boldsymbol{x}]^c)}{A(b(t))} \underset{t\to\infty}{\to} H(\boldsymbol{x})
\end{align}
locally uniformly in $\boldsymbol{x}\in (0,\infty]^d \setminus \{\boldsymbol{\infty}\}$, where $H$ is a function that is non-zero and  finite. Then $\bX$ is {\it{second order regularly varying}} with parameters $\alpha > 0$ and $\rho\le 0$. We write $\bX\in 2\MRV_{-\alpha,\rho}(b,A,\nu,H)$.
\end{defi}

\section{Proofs}\label{A:proof}
\subsection{Proofs of Lemmas of Section~\ref{sec:frame}}\label{A:sec2}

We can prove Lemma~\ref{lem1} and Lemma~\ref{lem1-bis} using pdf, when existing, or characteristic functions. We proceed for instance via pdf for the proof of the first lemma, as it is a very standard way (see e.g. \citet{david2004order}),
and via characteristic function for the second one to show that it can be done in full generality.

Starting with Lemma~\ref{lem1}, let us prove the following result.\\[.5ex]
{\it 
Let $f_{\mathbf{X}_{(1)},\cdots, \mathbf{X}_{(n-1)}\,|\, \|\mathbf{X}_{(n)}\|}(\mathbf{x}_{1}, \cdots, \mathbf{x}_{n-1} \,|\,y)$ be the conditional pdf of $\left(\mathbf{X}_{(1)}, \dots, \mathbf{X}_{(n - 1)}\right)$ given the event $\left(\| \mathbf{X}_{(n)}\| = y\right)$.
For any $\mathbf{x}_1, \dots, \mathbf{x}_{n-1} \in \mathbb{R}^d$, $y \geq 0$, such that $ \|\mathbf{x}_1\| \leq \dots\leq \|\mathbf{x}_{n-1}\| \leq y$, we have
\begin{equation}\label{eq:lem1-pdf}
  f_{\mathbf{X}_{(1)},\cdots, \mathbf{X}_{(n-1)}\,|\, \|\mathbf{X}_{(n)}\|}(\mathbf{x}_1, \dots, \mathbf{x}_{n-1}\,|\, y) = (n-1)! \,\prod_{i=1}^{n-1} f_{\mathbf{X}\,|\,\| \mathbf{X}\|}(\mathbf{x}_i\,|\, y),
\end{equation}
where  $f_{\mathbf{X} \,|\,\| \mathbf{X}\|}(\cdot\,|\, y)$ is defined in \eqref{eq:def_trunc_pdf}.
}\\[.5ex]
%
{\it Proof . } 
Let $[y, y + \delta y)$ be a `small' interval from $\mathbb{R}^+$ and, for $i = 1,\dots,n-1,$ let $[\mathbf{x}_i, \mathbf{x}_i + \delta \mathbf{x}_i) $ denote a `small' cube in $\mathbb{R}^d$ with initial vertex $\mathbf{x}_i$ and measure $\delta \mathbf{x}_i$. As $\mathbf{X}_{(1)}, \dots, \mathbf{X}_{(n )}$ are the ordered (by norm) statistics (see \eqref{ordered statistics}), there are $n!$ permutations of the vector $\left(\mathbf{X}_{1}, \dots, \mathbf{X}_{n}\right)$ to be $\left(\mathbf{X}_{(1)}, \dots, \mathbf{X}_{(n )}\right)$, hence we can write
\begin{equation*}
\begin{split}
  &\  \p{ \mathbf{X}_{(1)}\in [\mathbf{x}_1, \mathbf{x}_1 + \delta \mathbf{x}_1), \cdots, \mathbf{X}_{(n - 1)} \in [\mathbf{x}_{n-1}, \mathbf{x}_{n-1} + \delta \mathbf{x}_{n-1}),  \| \mathbf{X}_{(n )} \| \in [y, y + \delta y)  } \\
  &= n!\,\prod_{i=1}^{n-1}\Bigl( f(\mathbf{x}_i)  \delta \mathbf{x}_i\Bigr)\,
f_{\|\mathbf{X} \|} (y)  \delta y  + O((\delta y)^2).  
  \end{split}
\end{equation*}
Also, as $\|\mathbf{X}_{(n )} \|$ is the maximum of the iid rv
$\|\mathbf{X}_{1} \| , \dots, \|\mathbf{X}_{n} \|$, we have
\begin{equation}\label{def:density maximum}
   f_{\|\mathbf{X}_{(n)}\|}(y) = n F_{\|\mathbf{X}\|}^{n-1}(y) f_{\|\mathbf{X}\|}(y)
\end{equation}
and
\begin{equation*}
  \p{ \|\mathbf{X}_{(n)}\| \in [y, y + \delta y)  } = n \,F_{\|\mathbf{X} \|}^{n-1}(y) \,f_{\|\mathbf{X} \|}(y) \,\delta y + O\left((\delta y)^2\right).
\end{equation*}
Then
\begin{equation*}
  f_{\mathbf{X}_{(1)},\cdots, \mathbf{X}_{(n-1)}\,|\, \|\mathbf{X}_{(n)}\|}(\mathbf{x}_1, \dots, \mathbf{x}_{n-1}\,|\, y) = (n-1)! \,
\prod_{i=1}^{n-1} \frac{f(\mathbf{x}_i)}{ F_{\|\mathbf{X} \|}(y)},
\end{equation*}
hence the result \eqref{eq:lem1-pdf}. \hfill $\Box$ \\[1ex]

Turning to Lemma~\ref{lem1-bis}, it can be formalized as follows, using characteristic functions:\\[.5ex]
{\it 
For any vectors $\mathbf{x}_1, \cdots, \mathbf{x}_n \in \mathbb{R}^d$ and any real $y > 0$, we have
\begin{equation}\label{eq:lem1bis-chf}
\e \left[ e^{ \sum_{k=1}^n i\left\langle \mathbf{x}_k, \mathbf{X}_{(k)}\right\rangle} \Big| \|\mathbf{X}_{(n)}\|  = y\right] = \e  \left[e^{ \sum_{k=1}^{n-1} i\left\langle \mathbf{x}_k, \mathbf{Y}_{(k)}\right\rangle} \right] \e \left[ e^{ i\left\langle \mathbf{x}_n, \mathbf{X}_{n}\right\rangle} \Big| \|\mathbf{X}_n\| = y  \right].
\end{equation} 
}

\begin{proof} 
To simplify the notations, we denote by $R_1, \cdots, R_n$ the norms $ \|\mathbf{X}_1\|, \cdots, \|\mathbf{X}_n\|$. Since $\left(R_1, \cdots, R_n\right)$ are iid, for all $\mathbf{x}_1,\cdots,\mathbf{x}_n \in \mathbb{R}^d, y > 0$, we have
\begin{equation*}
\begin{split}
& \e \left[ e^{ \sum_{k=1}^n i\left\langle \mathbf{x}_k, \mathbf{X}_{(k)}\right\rangle} \Big| R_{(n)}  = y\right] = n!\, \e \left[ e^{ \sum_{k=1}^n i\left\langle \mathbf{x}_k, \mathbf{X}_{k}\right\rangle} \mathbbm{1}_{\skk{ R_1 \leq \cdots \leq R_{n-1} \leq R_n} }  \Big| R_{(n)}  = y\right]  \\
& = (n-1)! \, \e \left[ \left( e^{ \sum_{k=1}^{n-1} i\left\langle \mathbf{x}_k, \mathbf{X}_{k}\right\rangle} \mathbbm{1}_{\skk{ R_1 \leq \cdots \leq R_{n-1} \leq y} }\right) \left(n  e^{ i\left\langle \mathbf{x}_n, \mathbf{X}_{n}\right\rangle} \mathbbm{1}_{\skk{R_n =y}} \right)    \Big| R_{(n)}  = y\right]. 
\end{split}
\end{equation*} 
The random variables in the two big brackets being independent, since one depends on $\mathbf{X}_1, \cdots, \mathbf{X}_{n-1}$ and the other on $\mathbf{X}_n$, the last expression equals
\begin{equation}\label{eq:lastcond}
\begin{split}
\e \left[ (n-1)! \, e^{ \sum_{k=1}^{n-1} i\left\langle \mathbf{x}_k, \mathbf{X}_{k}\right\rangle} \mathbbm{1}_{\skk{ R_1 \leq \cdots \leq R_{n-1} \leq y} }   \Big| R_{(n)}  = y\right]  \e \left[n  e^{ i\left\langle \mathbf{x}_n, \mathbf{X}_{n}\right\rangle}  \mathbbm{1}_{\skk{R_n =y}}   \Big| R_{(n)}  = y\right].  
\end{split}
\end{equation} 
Note that, by Lemma \ref{lem1}, the first  expectation in \eqref{eq:lastcond} equals to
$\displaystyle
\e  \left[e^{ \sum_{k=1}^{n-1} i\left\langle \mathbf{x}_k, \mathbf{Y}_{(k)}\right\rangle} \right]
$.\\
Let us calculate the second expectation in \eqref{eq:lastcond}. Considering the $\sigma$-algebra $\sigma(R)$ generated by the iid rv $\left(R_1, \cdots, R_n\right)$ (with parent rv $R$), we can write, for any rv $\xi$, 
\\$\displaystyle 
    \e \left[\xi \,|\, R_{(n)}  \right] =    \e \left[\e \left[\xi \,|\,  \sigma(R) \right] \Big|  R_{(n)} \right] 
$
and
$\displaystyle
    \e \left[e^{ i\left\langle \mathbf{x}_n, \mathbf{X}_{n}\right\rangle}     \Big| \sigma(R)\right]  =      \e \left[e^{ i\left\langle \mathbf{x}_n, \mathbf{X}_{n}\right\rangle}     \Big|  R_n\right] =: h_{\mathbf{x}_n}(R_n)$.
\\We deduce that, for all $y>0$,
\begin{eqnarray*}
 \e \left[n \mathbbm{1}_{\skk{R_n =y}} e^{ i\left\langle \mathbf{x}_n, \mathbf{X}_{n}\right\rangle}    \Big| R_{(n)}  = y\right] 
&=&   \e \left[n \mathbbm{1}_{\skk{R_n =y}}  h_{\mathbf{x}_n}(R_n) \Big| R_{(n)}  = y\right] \\
&=& h_{\mathbf{x}_n}(y) \e \left[n \mathbbm{1}_{\skk{R_n =R_{(n)}}}   \Big| R_{(n)}  = y\right] = h_{\mathbf{x}_n}(y).
\end{eqnarray*}
Combining the results obtained for the two expectations in \eqref{eq:lastcond} provides \eqref{eq:lem1bis-chf}.
\end{proof}

\subsection{Proof of Theorem~\ref{th:Result-max}}\label{A:sec3}
%
\paragraph{\sc Preliminary result -}
To prove Theorem~\ref{th:Result-max}, we need to show that the inverse matrix of the covariance matrix $\Sigma(y)$ of the truncated random vector $\mathbf{Y}$ is bounded, namely: 
\begin{lem}\label{lem:covariance}
 Under the settings of Theorem \ref{th:Result-max}, for some $C > 0$, there exists $\delta >0$ such that 
$\p{\| \mathbf{X}\|\leq \delta} < 1$ and
  \begin{equation*}
    \| \Sigma(y)^{-1/2}\| \leq C \  \text{ for all } y \geq \delta.
  \end{equation*}
\end{lem}
\begin{proof}
 From the definition of the truncated distribution $F_{\mathbf{X}\,|\, \| \mathbf{X}\|}(\cdot | y)$ in \eqref{eq:def_trunc_distr} and continuity of  $F_{\|\mathbf{X}\|}(\cdot)$, one can show that,  for any sequence $(y_n)$ converging to $y_0\in\R^+$,
with $\displaystyle \p{\|\mathbf{X}\| \le  y_0} > 0$,  
we have
  \begin{equation*}
  \Sigma(y_n) \ra \Sigma(y_0) \text{ as } n \ra \infty,
  \end{equation*}
where, for convenience, we suppose $\Sigma(\infty):= \Sigma$, covariance matrix of the initial random vector $\mathbf{X}$ with c.d.f. $F$. Condition $(C2)$ in Theorem~\ref{th:Result-max} and the definition of $y_0$ imply that  $\Sigma(y_0)$ is positive definite.
Consequently, the square root of its inverse matrix exists, with finite norm:
\begin{equation}\label{eq:technical lem 1}
 \| \Sigma(y_0)^{-1/2}\| < \infty.
 \end{equation}
Now, assume that the statement of Lemma~\ref{lem:covariance} is false.
Then, for all $C >0$ and $\delta >0$ such that   $\p{\| \mathbf{X} \|\leq \delta } < 1$, there exists $y \geq \delta$ such that
$\displaystyle  \| \Sigma(y)^{-1/2}\| > C$.\\
  Therefore, there exists a sequence $(y_n)_n \geq \delta$ such that
  \begin{equation}\label{eq:technical lem 2}
  \| \Sigma(y_n)^{-1/2}\| \ra \infty \text{ as } n \ra \infty.
  \end{equation}
  But we can choose a subsequence $(y_{n_k})_{k \geq1}$ such that, as $k \ra \infty$, either (a) $y_{n_k} \ra \infty$, or (b) $y_{n_k} \ra y_0 \geq \delta$ for some point  $y_0 \geq \delta$.
 In both cases, \eqref{eq:technical lem 2} contradicts \eqref{eq:technical lem 1}. 
We conclude that the statement of Lemma~\ref{lem:covariance} is true.
\end{proof}

\paragraph{\sc Proof of Theorem~\ref{th:Result-max} -}

Let $C$ denote a positive constant that may vary from line to line, all along the proof.
Let $\mathbf{Y}_1, \dots, \mathbf{Y}_{n-1}$ be i.i.d. random vectors, with parent random vector $\mathbf{Y}$, having the truncated distribution $F_{\mathbf{X}\,|\, \| \mathbf{X}\|}(\cdot | y)$ defined in \eqref{eq:def_trunc_distr}.
We have, using conditional probabilities and applying Lemma~\ref{lem1}, 
\begin{equation*}
\begin{split}
\p{\mathbf{S}_n \in \mathbf{B}}   &= \p{  \mathbf{X}_{(1)}+  \dots  + \mathbf{X}_{(n)} \in \mathbf{B}} \\
&= \int_{\R^d} f_{(n)}(\mathbf{x}) \p{ \mathbf{X}_{(1)}+  \dots + \mathbf{X}_{(n-1)} \in \mathbf{B} - \mathbf{x}\,|\, \mathbf{X}_{(n)} = \mathbf{x} } \mathrm{d}\mathbf{x}\\
 &= \int_{\R^d} f_{(n)}(\mathbf{x}) \p{ \mathbf{T}_{n-1} \in \mathbf{B} - \mathbf{x} } \mathrm{d}\mathbf{x}, \quad \text{where} \quad T_{n-1}:= \sum_{k=1}^{n-1} \mathbf{Y}_k.
\end{split}
\end{equation*}
Using Definition~\ref{def:Normex} of the $d$-Normex distribution,  we can write, for any $\mathbf{B} \in \mathcal{C}$,
\begin{equation}\label{eq:proof_th1}
  \left| \p{\mathbf{S}_n \in \mathbf{B}} - G_n(\mathbf{B}) \right|
=  \int_{\R^d} f_{(n)}(\mathbf{x}) \,\left|  \p{\mathbf{T}_{n-1} \in \mathbf{B} - \mathbf{x}} - \Phi_{(n-1)\mathbf{\mu}(\|\mathbf{x}\|), (n-1) \Sigma(\|\mathbf{x}\|)} \left(\mathbf{B} - \mathbf{x}\right) \right| \mathrm{d}\mathbf{x}.
\end{equation}
Further, 
under $(C2)$,  the pdf of the rv  $\|\mathbf{X}\|$ is regularly varying,
$\displaystyle f_{\|X\|}(\cdot)\in\RV_{-\a-1}$,  {\it i.e.} (see e.g., Bingham et al. \citet{bingham:goldie:teugels:1989}) there exists a slowly varying function $L(\cdot)$ such that \begin{equation}\label{eq:fnormRV}
  f_{\|X\|}(y)  = L(y) y^{-\alpha-1}.
\end{equation}
Let us choose $\delta>0$ that satisfies Lemma~\ref{lem:covariance} and such that 
the function of $y$ defined by $\displaystyle \sup_{t \in [\delta, y]} L(t)$ is slowly varying (see e.g. ex.4 p.58 in \citet{bingham:goldie:teugels:1989}). \\
Splitting the integration domain of the integral of \eqref{eq:proof_th1} into two disjoint sets  $\skk{ \mathbf{x}: \| \mathbf{x}\| < \delta}$ and $ \skk{ \mathbf{x}: \| \mathbf{x}\| \geq \delta}$, we have,  
on one hand, 
\begin{equation}\label{eq:proof_th2}
\begin{split}
  &\int_{ \| \mathbf{x}\| < \delta}  f_{(n)}(\mathbf{x}) \,\left|  \p{\mathbf{T}_{n-1} \in \mathbf{B} - \mathbf{x}} - \Phi_{(n-1)\mathbf{\mu}(\|\mathbf{x}\|), (n-1) \Sigma(\|\mathbf{x}\|)} \left(\mathbf{B} - \mathbf{x}\right) \right| \mathrm{d}\mathbf{x} \\  &\leq  \int_{ \| \mathbf{x}\| < \delta} f_{(n)}(\mathbf{x}) \mathrm{d}\mathbf{x} = \p{\|\mathbf{X}_{(n)}\| < \delta} =\left(\p{\|\mathbf{X}\| < \delta}\right)^n.
  \end{split}
\end{equation}
On the other hand, since $\mathbf{Y}_1, \dots, \mathbf{Y}_{n-1}$  are i.i.d bounded random vectors with a non degenerate distribution (via $(C2)$),  we can use the non-generalized Berry-Esseen inequality recalled in Proposition \ref{thBEine}, and obtain
\begin{equation}\label{eq:proof_th3}
\begin{split}
&\int_{ \| \mathbf{x}\| \geq \delta} f_{(n)}(\mathbf{x}) \,\left|  \p{\mathbf{T}_{n-1} \in \mathbf{B} - \mathbf{x}} - \Phi_{(n-1)\mathbf{\mu}(\|\mathbf{x}\|), (n-1) \Sigma(\|\mathbf{x}\|)} \left(\mathbf{B} - \mathbf{x}\right) \right| \mathrm{d}\mathbf{x} \\
 &\leq \frac{C }{\sqrt{n}} \int_{\| \mathbf{x}\| \geq \delta} f_{(n)}(\mathbf{x})\, \e \left(\left\| \left(\Sigma(\|\mathbf{x}\|)\right)^{-1/2} \left( \mathbf{Y} - \mathbf{\mu}(\|\mathbf{x}\|) \right) \right\|^3 \right) \mathrm{d}\mathbf{x},
 \end{split}
 \end{equation}
$\mathbf{Y}$ denoting the parent random vector of $(\mathbf{Y}_1, \dots, \mathbf{Y}_{n-1})$.
%
It is straightforward to see that the last integral term in \eqref{eq:proof_th3} can be written only in terms of $\|\mathbf{x}\|$, which we denote by $y$, and of the pdf of the rv $\|\mathbf{X}_{(n)}\|$:
\begin{equation}\label{eq:proof_th31}
\int_{\| \mathbf{x}\| \geq \delta}\!\!\! f_{(n)}(\mathbf{x})\, \e \left(\left\| \left(\Sigma(\|\mathbf{x}\|)\right)^{-1/2} \!\! \left( \mathbf{Y} - \mathbf{\mu}(\|\mathbf{x}\|) \right) \right\|^3 \right) \!\mathrm{d}\mathbf{x}
=  \int_{y \geq \delta} \!\!\!f_{\|\mathbf{X}_{(n)}\|}(y)\, \e \left(\left\| \left(\Sigma(y)\right)^{-1/2} \!\!\left( \mathbf{Y} - \mathbf{\mu}(y) \right) \right\|^3 \right) \!\mathrm{d}y.
\end{equation}
For  $ y \ge \delta $, using  Lemma~\ref{lem:covariance}, it is straightforward to see that 
\begin{equation*}
  \e \left(\left\| \Sigma^{-1/2}(y) \left( \mathbf{Y} - \mathbf{\mu}(y) \right) \right\|^3 \right) 
 \leq C\,  \e  \, \left\| \mathbf{Y} - \mathbf{\mu}(y)  \right\|^3 
\leq C\left( \e   \left\| \mathbf{Y}  \right\|^3 + \|\mathbf{\mu}(y) \|^3\right) \leq C\, \e   \left\| \mathbf{Y}  \right\|^3. 
\end{equation*}
Splitting the latter expectation into 2 parts according to $(\left\|  \mathbf{X} \right\| \le \delta)$ or $(\left\|  \mathbf{X} \right\| \in [\delta, y])$ (and recalling that $\mathbf{Y}$ has the truncated cdf $F_{\mathbf{X}\,|\, \| \mathbf{X}\|}(\cdot | y)$),   we can write
\begin{equation}\label{eq:proof_th4}
\begin{split}
  \e \left(\left\| \Sigma^{-1/2}(y) \left( \mathbf{Y} - \mathbf{\mu}(y) \right) \right\|^3 \right) 
& \leq C \left( \delta^{3} +  \, \e\left( \left\|  \mathbf{X} \right\|^3\, \mathbbm{1}_{\skk{\left\|  \mathbf{X} \right\| \in [\delta, y] }}\right) \right) 
\\  & \leq \,C \left( \delta^{3} + \int_{\delta}^{y} t^{3} L(t) t^{-\alpha - 1} \mathrm{d}t \right) \leq  C L(y)\, y^{3-\alpha},
  \end{split}
\end{equation}
  where the last inequality comes from the Karamata's integral theorem (\citet{Karamata:1933}; see e.g. Proposition 1.5.8 p.26 in \citet{bingham:goldie:teugels:1989}).
So, from \eqref{eq:proof_th1}--\eqref{eq:proof_th4} we obtain
\begin{equation}\label{eq1}
\!\! \!\left| \p{\mathbf{S}_n \in \mathbf{B}} - G_n(\mathbf{B}) \right|   \leq
\frac{C }{\sqrt{n}} \int_{y \geq \delta} f_{\|\mathbf{X}_{(n)}\|}(y) L(y)  y^{3-\alpha} \mathrm{d}y
   = \frac{C }{\sqrt{n}} \, \e \left( L(\|\mathbf{X}_{(n)}\|)     \|\mathbf{X}_{(n)}\|^{3-\alpha} \mathbbm{1}_{\skk{  \|\mathbf{X}_{(n)}\| \geq \delta }}\right).
 \end{equation}
 For $n$ large enough such that $n^{1/\alpha} > \delta$, using \eqref{def:density maximum} and \eqref{eq:fnormRV} gives
  \begin{align*}
     \e \left( L(\|\mathbf{X}_{(n)}\|)     \|\mathbf{X}_{(n)}\|^{3-\alpha}  \mathbbm{1}_{\skk{  \|X_{(n)}\| > n^{1/\alpha} }} \right)&= \int_{n^{1/\alpha}}^{\infty} L(y)  y^{3-\alpha}  n F_{\|\mathbf{X}\|}^{n-1}(y) f_{\|\mathbf{X}\|}(y) \, \mathrm{d} y \\
     &\leq  n \int_{n^{1/\alpha}}^{\infty} L(y) y^{2-2\alpha} \,\mathrm{d}y  \underset{n\to\infty}{\sim}
L(n^{1/\alpha}) n^{\frac{3-\alpha}{\alpha}}.
  \end{align*}
  Consequently, we have, $L$ denoting a slowly varying function at infinity that may vary from line to line,
  \begin{eqnarray}\label{eq2}
&& \e \left( L(\|\mathbf{X}_{(n)}\|)     \|\mathbf{X}_{(n)}\|^{3-\alpha} \mathbbm{1}_{\skk{  \|\mathbf{X}_{(n)}\| \geq \delta }}\right) \nonumber\\
&&\quad\le  \e \left( L(\|\mathbf{X}_{(n)}\|)     \|\mathbf{X}_{(n)}\|^{3-\alpha} \mathbbm{1}_{ \skk{ \|\mathbf{X}_{(n)}\| \in [\delta, \,n^{1/\alpha}] }}\right)  +  \e \left( L(\|\mathbf{X}_{(n)}\|)     \|\mathbf{X}_{(n)}\|^{3-\alpha} \mathbbm{1}_{ \skk{  \|\mathbf{X}_{(n)}\| > n^{1/\alpha} }}\right) \nonumber\\
&& \quad\leq \sup_{y \in [\delta, \,n^{1/\alpha}]} L(y) n^{\frac{3-\alpha}{\alpha}} + L(n^{1/\alpha}) n^{\frac{3-\alpha}{\alpha}} \,= \,L(n) \,n^{\frac{3-\alpha}{\alpha}}.
  \end{eqnarray}
Combining inequalities \eqref{eq1} and \eqref{eq2} provides the result. \hfill $\Box$

\subsection{Proof of the results given in Section~\ref{sec:mrv}}\label{A:sec4}

\subsubsection{Proof of Lemma~\ref{lemma:2RVrelation} and Theorem~\ref{th:MRV result}}\label{A:sec4lem}

First let us recall Lemma~\ref{lemma:2RVrelation} for convenience:\\
{\it 
If $f_X \in 2\RV_{-\alpha - 1, -\beta}$, with $\alpha >0$ and $\beta >0$, then the derivative  $g'$ of $g$ defined in \eqref{eq:def-g}, satisfies $g'\in 2\RV_{\frac1\a - 1,\, \rho}$ , where $\rho :=-\min\skk{1, \frac{\beta}{\alpha}}$.
}

\begin{proof} 
First we prove that $g$ belongs to $2\RV_{\gamma - 1, \rho}$, with $\gamma=1/\a$.\\
It is straightforward to check that: 
\begin{equation}\label{eq:rep_for_g}
g(t) = F_X^{\leftarrow}(e^{-\frac{1}{t}}) = \Bar{F}_X^{\leftarrow}(1 - e^{-\frac{1}{t}}).
\end{equation}
Notice that the function $h(t) := 1 - e^{-\frac{1}{t}}$ belongs to $2\RV_{-1,-1}$.  
Now, from the $2RV$ condition on $f_X$, we obtain that $\Bar{F}_X \in 2\RV_{-\alpha, - \beta}$ with parameters $ \alpha >0, \beta >0$; see e.g.  Proposition 6  in \citet{hua:joe:2011b}. 
Then, applying Proposition 2.6 in \citet{lv:mao:hu:2012a} and denoting by $2\RV_* (0+)$ the class of regularly varying functions at $0+$, 
we have 
\begin{equation*}
    \Bar{F}_X^{\leftarrow} \in 2\RV_{-\frac{1}{\alpha},  -\frac{\beta}{\alpha}} (0+),
\end{equation*}
from which we deduce, using Proposition 2.8 in \citet{lv:mao:hu:2012a} (for a composition of functions) and  \eqref{eq:rep_for_g}, 
\begin{equation*}
    g \in 2\RV_{\frac{1}{\alpha},  -\min\skk{1, \frac{\beta}{\alpha}}}.
\end{equation*}
Now let us look on the derivative $g'$. From \eqref{eq:rep_for_g}, we have
$\displaystyle g'(t) = \frac{e^{-\frac{1}{t}} \,t^{-2}}{f_X\left( g(t)\right)} $.

Again, by Proposition 2.8 and Proposition 2.5 in \citet{lv:mao:hu:2012a}, 
we have
\begin{equation*}
    \frac{1}{f_X\left( g(t)\right)}  \in 2\RV_{\frac{1}{\alpha} + 1,  -\min\skk{1, \frac{\beta}{\alpha}}} \quad
\text{and} \quad 
   e^{-\frac{1}{t}}\, t^{-2} \in  2\RV_{-2, -1}. 
\end{equation*}
Hence, multiplying the $2\RV$ functions, we obtain
$\displaystyle g' \in 2\RV_{\frac{1}{\alpha} -1,  -\min\skk{1, \frac{\beta}{\alpha}}}$.
\end{proof}

\begin{proof}[Proof of Theorem~\ref{th:MRV result}]

Recall that $\mathbf{Y}_1, \dots, \mathbf{Y}_{n-1}$ denote i.i.d. random vectors having the truncated distribution  $F_\mathbf{X}(\cdot\,|\, y)$ defined in \eqref{eq:def_trunc_distr}, while $\overset{\circ}{\mathbf{X}}_y$ is a family of random vectors with distribution \eqref{eq:def-Zcdf}, independent of $\skk{\mathbf{Y}_k}.$
Using conditional probabilities and applying Lemma~\ref{lem1-bis}, we have
\begin{equation*}
\begin{split}
\p{\mathbf{S}_n \in \mathbf{B}}   
&= \int\limits_{0}^{\infty} f_{\|\mathbf{X}_{(n)}\|}(y) \p{ \mathbf{X}_{(1)}+  \dots + \mathbf{X}_{(n-1)} + \mathbf{X}_{(n)} \in \mathbf{B} \,|\, \|\mathbf{X}_{(n)}\| = y } \mathrm{d}y
\\ &= \int\limits_{0}^{\infty} f_{\|\mathbf{X}_{(n)}\|}(y) \p{ \mathbf{T}_{n-1}  + \overset{\circ}{\mathbf{X}}_y \in \mathbf{B}  } \mathrm{d}y, \quad \text{where} \quad \mathbf{T}_{n-1}:= \sum_{k=1}^{n-1} \mathbf{Y}_k.
\end{split}
\end{equation*}
Using the definition \eqref{eq:normex-cdf-MRV2} of the MRV-Normex cdf, then the triangle inequality, we have, for any $\mathbf{B}\in \cal{C}$,
\begin{equation}\label{local1}
\begin{split}
   \left| \p{S_n \in \mathbf{B}} - MG_n(\mathbf{B}) \right| &= \left|  \int\limits_{0}^{\infty} f_{\|\mathbf{X}_{(n)}\|}(y)  \p{\overset{\circ}{\mathbf{X}}_y + \mathbf{T}_{n-1} \in \mathbf{B}} \mathrm{d} y - \int\limits_{0}^{\infty} f_{H_{\a,n}}(y) \p{ y \cdot \bTheta + Z_y \in \mathbf{B} } \mathrm{d}y\right|  \\
     & \leq  \int\limits_{0}^{\infty} \left|  f_{\|\mathbf{X}_{(n)}\|}(y)  - f_{H_{\a,n}}(y)  \right| \mathrm{d} y      \\
     & + \int\limits_{0}^{\infty} f_{\|\mathbf{X}_{(n)}\|}(y)  \Bigl| \p{ \overset{\circ}{\mathbf{X}}_y + \mathbf{T}_{n-1} \in \mathbf{B}  } - \p{\overset{\circ}{\mathbf{X}}_y + Z_y \in \mathbf{B} }  \Bigr| \mathrm{d} y \\
     & + \int\limits_{0}^{\infty} f_{\|\mathbf{X}_{(n)}\|}(y)  \Bigl|  \p{\overset{\circ}{\mathbf{X}}_y + Z_y \in \mathbf{B}} - \p{ y\cdot\bTheta +   Z_y } \Bigr|\mathrm{d} y.
\end{split}
\end{equation}
Recall that, via the Scheff\'{e}  theorem (see e.g. \citet{Bill1968}, p. 224), if the cdf $F$ and $G$ are absolutely continuous, then
$\displaystyle
       \sup _{A \in B(\mathrm{R})}\left| F(A) - G(A) \right| = \frac{1}{2} \int_{-\infty}^{\infty} \left|F'(x) - G'(x) \right| \mathrm{d}x$.

Therefore, for the first integral in the right-hand side of the inequality \eqref{local1}, using Proposition~\ref{cvRate-Frechet} and this recall, there exists a slowly-varying function at infinity $L$ such that, for all $n \geq 1$,
\begin{equation}\label{local1_1}
 \int\limits_{0}^{\infty} \left|  f_{\|\mathbf{X}_{(n)}\|}(y)  - f_{H_{\a,n}}(y)  \right| \mathrm{d} y  \leq L(n) n^{-{\min\skk{1,\frac{\beta}{\alpha}}}}.
\end{equation}
%
For the second integral appearing in \eqref{local1}, following the same approach as for the proof of Theorem~\ref{th:Result-max}, we can prove that there exists a slowly varying function at infinity, $L$, such that, for all $n \geq 1$,
  \begin{equation}\label{eq:lem6:1}
 \int\limits_{0}^{\infty} f_{\|\mathbf{X}_{(n)}\|}(y) | \p{ \overset{\circ}{\mathbf{X}}_y + \mathbf{T}_{n-1} \in \mathbf{B}  } - \p{ \overset{\circ}{\mathbf{X}}_y+ Z_y\in B } | \mathrm{d}y \leq L(n) n^{-\frac{1}{2} + \frac{3-\alpha}{\alpha}}.
  \end{equation}
Indeed, choosing the constant $\delta >0$ as in the proof of Theorem~\ref{th:Result-max}  and splitting the integration domain of the integral in \eqref{eq:lem6:1} into two disjoint sets  $\skk{ y: y < \delta}$ and $ \skk{ y: y \geq \delta}$, we can write
\begin{equation}\label{eq:proof_lem6:1}
  \int\limits_{0}^{\delta}  f_{\|\mathbf{X}_{(n)}\|}(y) \left| \p{ \mathbf{T}_{n-1}  + \overset{\circ}{\mathbf{X}}_y \in \mathbf{B}  } - \p{\overset{\circ}{\mathbf{X}}_y+  Z_y \in B } \right| \mathrm{d}y 
\leq  \int\limits_{0}^{\delta}  f_{\|\mathbf{X}_{(n)}\|}(y) | \mathrm{d}y 
= (F_{\|\mathbf{X}\|}(\delta))^n.
\end{equation}
Let  $C$ denote a positive constant and $L$ a slowly varying function at infinity, which may vary from line to line.
Since $\mathbf{Y}_1, \dots, \mathbf{Y}_{n-1}$  are i.i.d bounded random vectors ($\mathbf{Y}$ denoting their parent random vector) with non degenerate distribution (condition $(C2)$), we can use once again the non-generalized Berry-Esseen inequality, and obtain
\begin{eqnarray}\label{eq:proof_lem6:2}
&&\int\limits_{\delta}^{\infty}  f_{\|\mathbf{X}_{(n)}\|}(y) \left| \p{ \mathbf{T}_{n-1}  + \overset{\circ}{\mathbf{X}}_y \in \mathbf{B}  } - \p{ \overset{\circ}{\mathbf{X}}_y +  Z_y \in \mathbf{B} } \right| \mathrm{d}y \nonumber\\
&& \leq \int\limits_{\delta}^{\infty}  f_{\|\mathbf{X}_{(n)}\|}(y)  
\sup_{\mathbf{B} \in \mathcal{C}} | \p{ \mathbf{T}_{n-1}  \in \mathbf{B}  } - \p{  Z_y \in \mathbf{B} } |\,
\mathrm{d}y \nonumber \\
 &&\leq \frac{C }{\sqrt{n}}\int\limits_{\delta}^{\infty}  f_{\|\mathbf{X}_{(n)}\|}(y) \, \e \left(\left\| \left(\Sigma(y)\right)^{-1/2} \left( \mathbf{Y} - \mathbf{\mu}(y) \right) \right\|^3_2 \right) \mathrm{d}y.
 \end{eqnarray}
Replicating the arguments in the proof of Theorem~\ref{th:Result-max}, we obtain \eqref{eq:proof_th4}, {\it i.e.}
$$\displaystyle \e \left(\left\| \Sigma^{-1/2}(y) \left( \mathbf{Y} - \mathbf{\mu}(y) \right) \right\|^3_2 \right) \le C L(y)\, y^{3-\alpha}.$$
We deduce that
\begin{equation}\label{eq:proof_lem6:4}
\begin{split}
&\int\limits_{0}^{\infty} f_{\|\mathbf{X}_{(n)}\|}(y) \left| \p{ \mathbf{T}_{n-1}  + \overset{\circ}{\mathbf{X}}_y \in \mathbf{B}  } - \p{\overset{\circ}{\mathbf{X}}_y +  Z_y \in B } \right| \mathrm{d}y \\
& \leq
\frac{C }{\sqrt{n}} \int_{ y \geq \delta} f_{\|\mathbf{X}_{(n)}\|}(y) L(y)  y^{3-\alpha} \mathrm{d}y
   = \frac{C }{\sqrt{n}} \, \e \left( L(\|\mathbf{X}_{(n)}\|)     \|\mathbf{X}_{(n)}\|^{3-\alpha}  \mathbbm{1}_{\skk{\|\mathbf{X}_{(n)}\| \geq \delta}}\right).
   \end{split}
 \end{equation}
We are back to the upper bound given in \eqref{eq1}. Combining it with \eqref{eq2} provides the claimed result \eqref {eq:lem6:1}.\\

The last integral in \eqref{local1} is taken care of, in the following technical lemma.
\begin{lem}\label{lem6.3}
There exists a slowly varying function at infinity, $L$, such that, for all $n \geq 1$, for any $\mathbf{B}\in \cal{C}$,
 \begin{equation*}
\int\limits_{0}^{\infty} f_{\|\mathbf{X}_{(n)}\|}(y)  \Bigl|  \p{\overset{\circ}{\mathbf{X}}_y + Z_y \in \mathbf{B}} - \p{ y\cdot\bTheta +   Z_y \in \mathbf{B}} \Bigr|\mathrm{d} y \leq  L(n) n^{-\frac{\rho}{\alpha}}.
\end{equation*}
\end{lem} 

Combining inequalities  \eqref{local1}, \eqref{local1_1}, \eqref{eq:lem6:1} and Lemma~\ref{lem6.3} provides the statement of Theorem~\ref{th:MRV result}.
\end{proof}

Let us turn to the proof of Lemma~\ref{lem6.3}.

\begin{proof}[Proof of Lemma \ref{lem6.3}]
 Define the sequence $(a_n)_{n\ge 1}$  with $\displaystyle a_n:= F^{\leftarrow}_{\|\mathbf{X}\|}\left(n^{-\frac{\ln n}{n}}\right) \underset{n\to\infty}{\ra} \infty$. We have
\begin{equation}\label{local2}
 \p{\|\mathbf{X}_{(n)}\| \leq a_n} = ( F_{\|\mathbf{X}\|}(a_n) )^n= n^{-\ln n}.
\end{equation}
Note that $a_n=\bar{F}^{\leftarrow}_{\|\mathbf{X}\|}\left(1-n^{-\frac{\ln n}{n}}\right)$.
Since the function $u(n):=1-n^{-\frac{\ln n}{n}} \in RV_{-1}$ (at infinity) and $\|\mathbf{X}\|\in\RV_{-\alpha}$, so $\bar{F}_{\|\mathbf{X}\|}^{\leftarrow}\in \RV_{-1/\alpha}$ at $0^+$, therefore 
\begin{equation}\label{local20}
a_n=\bar{F}^{\leftarrow}_{\|\mathbf{X}\|}\circ u\, (n)\, \in \RV_{1/\alpha} \quad\text{ (at infinity)}
\end{equation}
 (for a reminder on properties of $RV$, see e.g. \citet{lv:mao:hu:2012a}). 

As $\displaystyle \bar{F}_{\|\mathbf{X}\|}(a_n) = n \left(1 - n^{-\frac{\ln n}{n}}\right)$, there exists also a slowly varying function at infinity, $L$, such that
\begin{equation}\label{local22}
 n \bar{F}_{\|\mathbf{X}\|}(a_n) = L(n).
\end{equation}

Splitting the integral of Lemma~\ref{lem6.3} into two sets $\skk{y \leq  a_n}$ and $\skk{y > a_n}$ and using the pdf for the maximum \eqref{def:density maximum}, we can write
\begin{equation}\label{local3}
\begin{split}
  & \int\limits_{0}^{\infty} f_{\|\mathbf{X}_{(n)}\|}(y)  \Bigl|  \p{\overset{\circ}{\mathbf{X}}_y + Z_y \in \mathbf{B}} - \p{ y\cdot\bTheta +   Z_y \in \mathbf{B}} \Bigr|\mathrm{d} y
\\&\leq  \p{\|\mathbf{X}_{(n)}\| \leq a_n}
 + \int\limits_{a_n}^{\infty} n f_{\|\mathbf{X}\|}(y) F_{\|\mathbf{X}\|}^{n-1}(y) \Bigl|  \p{\overset{\circ}{\mathbf{X}}_y + Z_y \in \mathbf{B}} - \p{ y\cdot\bTheta +   Z_y \in \mathbf{B} } \Bigr|\mathrm{d} y 
\\& \le n^{-\ln n} + \int\limits_{a_n}^{\infty} n f_{\|\mathbf{X}\|}(y) \Bigl|  \p{\overset{\circ}{\mathbf{X}}_y + Z_y \in \mathbf{B}} - \p{ y\cdot\bTheta +   Z_y \in \mathbf{B}} \Bigr|\mathrm{d} y,
\quad\text{using}\,\eqref{local2}.
\end{split}
\end{equation}
To estimate the integral given in \eqref{local3}, we use a coupling technique, namely: For any random vector $(\xi, \eta)$ in some complete separable metric space, and for any measurable set $B$, we have 
\begin{equation*}
   \left| \p{\xi \in B} - \p{\eta \in B} \right| \leq \p{\xi \neq \eta}. 
\end{equation*}
So, for any joint distribution of $\overset{\circ}{\mathbf{X}}_y$ and $\bTheta$ (or, more generally, for any joint distribution of $\mathbf{X}$ and $\bTheta$), we have 
\begin{equation*}
    \Bigl|  \p{\overset{\circ}{\mathbf{X}}_y + Z_y \in \mathbf{B}} - \p{ y\cdot\bTheta +   Z_y \in\mathbf{B}} \Bigr| \leq \p{\overset{\circ}{\mathbf{X}}_y\neq y\cdot \bTheta },
\end{equation*}
which gives the following upper bound for the integral in \eqref{local3}, using \eqref{eq:def-Zcdf}:
\begin{equation*}
\int\limits_{a_n}^{\infty} n f_{\|\mathbf{X}\|}(y) \p{\frac{\overset{\circ}{\mathbf{X}}_y}{y}\neq \bTheta}\mathrm{d} y = n\,\p{  \frac{\mathbf{X}}{\|\mathbf{X}\|} \neq \bTheta, \|\mathbf{X}\| > a_n }.  
\end{equation*}
Therefore, we have
\begin{equation}\label{eq:local3bis}
\eqref{local3} \,\le \, n^{-\ln n} + n\,\p{\|\mathbf{X}\| > a_n }\p{  \frac{\mathbf{X}}{\|\mathbf{X}\|} \neq \bTheta   \, \Big|\, \|\mathbf{X}\| > a_n }.
\end{equation}
Now, at given $n$, we define the joint distribution of $\mathbf{X}$ and $\mathbf{\Theta}$ such that:\\
(i) $\mathbf{\Theta}$ is independent of the event $\skk{ \|\mathbf{X}\| > a_n }$,
\\(ii) ($\mathbf{X}\,/\,\|\mathbf{X}\|,\,\mathbf{\Theta}$) has a joint distribution defined by Dobrushin's theorem (see \citet{Dobrushin1970}). 
\\Hence we can write, using Condition $(M_{\Theta})$ to get the asymptotic behavior,
\begin{equation*}
    \p{  \frac{\mathbf{X}}{\|\mathbf{X}\|} \neq \bTheta   \, \Big|\, \|\mathbf{X}\| > a_n } = \sup_{\mathbf{B} \subseteq \mathcal{S}_1}\left| \p{  \frac{\mathbf{X}}{\|\mathbf{X}\|}  \in \mathbf{B}\, \Big|\, \|\mathbf{X}\| > a_n } - \p{ \bTheta   \in \mathbf{B}} \right| \;\underset{n\to\infty}{\sim} \,A(a_n),
\end{equation*}
where $A(a_n)\in \RV_{-\rho/\alpha}$ by combining $(M_{\Theta})$ and \eqref{local20}.\\
Reporting this last result in \eqref{eq:local3bis} and using \eqref{local22}, we obtain
$$
\eqref{local3} \,\le \, n^{-\ln n} + n \bar{F}_{\|\mathbf{X}\|}(a_n) \, A(a_n) = n^{-\ln n} + L(n) \, A(a_n), 
$$
from which the result of Lemma~\ref{lem6.3}  follows.
\end{proof}

\subsubsection{Discussion on Condition \texorpdfstring{$(M_\Theta)$}{Lg}}\label{A:sec4disc}

\paragraph{Proof of 3(b) given in Remark \ref{rk:theoMRV} -}
Let us prove that Condition~\eqref{eq:factorization} given on the form of the pdf of a random vector is a necessary and sufficient condition to have the independence between the norm and direction of this vector. 

\begin{proof}
Notice that the unit ball $\skk{\mathbf{x}\,:\, \|\mathbf{x}\| \leq 1}$ is a convex set, 
which boundary $\mathcal{S}_1$ is a differentiable manifold.  Further we want to integrate on $\mathcal{S}_1$, and assume w.l.o.g. that the atlas for $\mathcal{S}_1$ consists of one chart only. (Indeed, to provide the integration in the case of many charts, one would just need to use a partition and sum the integrals over all charts.) 

If a map $\varphi:\, \R^{d-1} \ra \mathcal{S}_1$ parameterizes a set $B$ on the unit sphere $\mathcal{S}_1$,  one can introduce the map
\begin{align*}
	\widehat{\varphi} \,:\, &\R^{d-1}   \ra \mathcal{S}_y\\
				& \mathbf{u} \ra \widehat{\varphi} (\mathbf{u}) := y \,\varphi(\mathbf{u}) 
\end{align*} 
that parameterizes a set $B_y:=\skk{\mathbf{x}\,:\, \mathbf{x}/\|\mathbf{x}\| \in B, \, \|\mathbf{x}\|=y}$ on the sphere $\mathcal{S}_y:=\skk{\mathbf{x} \in \mathbb{R}^d\,:\, \|\mathbf{x}\| = y}$.
Note that both functions $\varphi,\, \widehat{\varphi}$ are differentiable a.e.  

Now, assuming that the pdf of $\mathbf{X}$ satisfies \eqref{eq:factorization}, the joint distribution of the random variable $\|\mathbf{X}\|$ and  the random vector $\mathbf{X}/ \|\mathbf{X}\|$ can be expressed as:
$$
 \p{\|\mathbf{X}\| \leq t,\, \frac{\mathbf{X}}{\|\mathbf{X}\|} \in B} 
= \!\!\!\!\!\!  \int\limits_{\|\mathbf{x}\| \leq t,\,\frac{\mathbf{x}}{\|\mathbf{x}\|}\in B } \!\!\!\!\!\!  h(\|\mathbf{x}\|) g(\mathbf{x}/ \|\mathbf{x}\|) \mathrm{d}\mathbf{x} 
= \!\!\!\!\!\!\int\limits_{y \in [0,t], \,\mathbf{u} \in \varphi^{-1}(B)} \!\!\!\!\!\!\!\!\! h(y) g(\varphi(\mathbf{u})) J(y,\mathbf{u}) \mathrm{d}y \mathrm{d}\mathbf{u}
$$
where $J(y,u)$ denotes the Jacobian for the change of variables
\begin{equation}\label{eq:change_of_var}
    \mathbf{x} = y \varphi(\mathbf{u}).
\end{equation}
Let us compute $J(y,u)$, which is the determinant of the $d\times d$-matrix $\displaystyle \left(\frac{\partial \mathbf{x}}{ \partial(y, \mathbf{u})}\right)$.
\begin{equation*}
\begin{split}
    J(y,\mathbf{u}) &:= 
\left|
\begin{array}{cccc}
    \frac{\partial x_1}{\partial y} & \frac{\partial x_2}{\partial y} &  \dots  & \frac{\partial x_d}{\partial y} \\
    \frac{\partial x_1}{\partial u_1} & \frac{\partial x_2}{\partial u_1} &  \dots  & \frac{\partial x_d}{\partial u_1} \\
    \vdots & \vdots &  \ddots & \vdots \\
    \frac{\partial x_1}{\partial u_{d-1}} & \frac{\partial x_2}{\partial u_{d-1}}  &  \dots  & \frac{\partial x_d}{\partial u_{d-1}} 
\end{array} 
\right| 
=\left|
\begin{array}{cccc}
    \varphi_1 & \varphi_2 &  \dots  & \varphi_d \\
     y \frac{\varphi_1}{\partial u_1} &  y \frac{\varphi_2}{\partial u_1} &  \dots  &  y \frac{\varphi_d}{\partial u_1} \\
    \vdots & \vdots &  \ddots & \vdots \\
     y \frac{\varphi_1}{\partial u_{d-1}} &  y \frac{\varphi_2}{\partial u_{d-1}}  &  \dots  &  y \frac{\varphi_d}{\partial u_{d-1}} 
\end{array} 
\right|\\
& =   y^{d-1} \times\left|
\begin{array}{cccc}
    \varphi_1 & \varphi_2 &  \dots  & \varphi_d \\
      \frac{\varphi_1}{\partial u_1} &   \frac{\varphi_2}{\partial u_1} &  \dots  &   \frac{\varphi_d}{\partial u_1} \\
    \vdots & \vdots &  \ddots & \vdots \\
      \frac{\varphi_1}{\partial u_{d-1}} &   \frac{\varphi_2}{\partial u_{d-1}}  &  \dots  &   \frac{\varphi_d}{\partial u_{d-1}} 
\end{array} 
\right| 
= \,y^{d-1} \times J(1,\mathbf{u}).
\end{split}
\end{equation*}
Therefore, we can write
\begin{equation}\label{eq:jointdistrnorm_direc}
\begin{split}
&  \p{\|\mathbf{X}\| \leq t, \mathbf{X}/ \|\mathbf{X}\| \in B} 
   = \int\limits_{y \in [0,t], \,\mathbf{u} \in \varphi^{-1}(B)  } h(y)y^{d-1} g(\varphi(\mathbf{u}))  J(1,\mathbf{u}) \,\mathrm{d}y \,\mathrm{d}\mathbf{u} \\ 
&=  \int\limits_{y \in [0,t]}  c_d\, h(y)y^{d-1}\,\mathrm{d}y  \;\;\frac{1}{c_d}\!\!\! \int\limits_{ \mathbf{u} \in \varphi^{-1}(B)  } \!\!\!\!\!\! g(\varphi(\mathbf{u}))  J(1,\mathbf{u})  \,\mathrm{d}\mathbf{u} 
  =\p{\|\mathbf{X}\| \leq t}\p{\frac{\mathbf{X}}{\|\mathbf{X}\| }\in B} 
    \end{split}
\end{equation}
with $\displaystyle c_d:=  \int\limits_{ \mathbf{u} \in \varphi^{-1}(\mathcal{S}_1)  }  g(\varphi(\mathbf{u}))  J(1,\mathbf{u})  \,\mathrm{d}\mathbf{u} =\left( \int\limits_0^{\infty}  \, h(y)y^{d-1}\mathrm{d}y\right)^{-1}$. \\
We deduce that $\|\mathbf{X}\|$ and $\displaystyle\frac{\mathbf{X}}{ \|\mathbf{X}\|}$ are independent.

Let us prove the converse, assuming that $\|\mathbf{X}\|$ and $\mathbf{X}/ \|\mathbf{X}\|$ are independent.

Note that, using the same change of variable \eqref{eq:change_of_var}, we can write
\begin{equation}\label{eq:chgeV}
\int\limits_{\|\mathbf{x}\| \leq t,\,\mathbf{x}/ \|\mathbf{x}\| \in B }  f_{\mathbf{X}}(\mathbf{x}) \mathrm{d}\mathbf{x} = 
\int\limits_{y \in [0,t], \,\mathbf{u} \in \varphi^{-1}(B)  }  f_{\mathbf{X}}(y \varphi(\mathbf{u})) \,y^{d-1} J(1,\mathbf{u}) \,\mathrm{d}y \,\mathrm{d}\mathbf{u}.
\end{equation}
Moreover, we have
\begin{align}\label{eq:pdctChge}
\p{\|\mathbf{X}\| \leq t, \mathbf{X}/ \|\mathbf{X}\| \in B} & =\p{\|\mathbf{X}\| \leq t}\p{\mathbf{X}/ \|\mathbf{X}\| \in B} \nonumber\\
& =  \int\limits_{\|\mathbf{x}\| \leq t } \!\!\! f_{\mathbf{X}}(\mathbf{x}) \,\mathrm{d}\mathbf{x} \; \int\limits_{\mathbf{x}/ \|\mathbf{x}\| \in B } \!\!\!\!\!\!  f_{\mathbf{X}}(\mathbf{x}) \,\mathrm{d}\mathbf{x}  \nonumber\\
     & = \int\limits_{y \in [0,t]; \,\mathbf{u}\in\R^{d-1} }\!\!\! \!\!\!\!\!\! f_{\mathbf{X}}(y \varphi(\mathbf{u})) y^{d-1} J(1,\mathbf{u}) \,\mathrm{d}y \,\mathrm{d}\mathbf{u}   
\!\!\!\!\!\!
\int\limits_{y\in\R^+;\, \mathbf{u} \in \varphi^{-1}(B)} \!\!\!\!\!\! \!\!\!\!\!\! f_{\mathbf{X}}(y \varphi(\mathbf{u})) y^{d-1} J(1,\mathbf{u}) \,\mathrm{d}y\, \mathrm{d}\mathbf{u}.   
\end{align}
Since $\displaystyle \p{\|\mathbf{X}\| \leq t, \mathbf{X}/ \|\mathbf{X}\| \in B}=\!\!\!\!\!\! \!\!\!\!\!\!\int\limits_{\|\mathbf{x}\| \leq t,\,\mathbf{x}/ \|\mathbf{x}\| \in B }  \!\!\!\!\!\! \!\!\!\!\!\! f_{\mathbf{X}}(\mathbf{x}) \mathrm{d}\mathbf{x} $, \,then the right hand side of \eqref{eq:chgeV} equals that of \eqref{eq:pdctChge}.
Consequently, for $t_0 \in (0,\infty)$ and $\varphi(\mathbf{u}_0) \in \mathcal{S}_1$, we can write
\begin{equation*}
    f_{\mathbf{X}}(t_0\, \varphi(\mathbf{u}_0)) \,t_0^{d-1} J(1,\mathbf{u}_0) = \!\!\!\!\!
    \int\limits_{\varphi(\mathbf{u}) \in \mathcal{S}_1} \!\!\!\!\!\!\! f_{\mathbf{X}}(t_0 \varphi(\mathbf{u})) \,t_0^{d-1} J(1,\mathbf{u})  \mathrm{d}\mathbf{u}  \, \times \int\limits_{0}^{\infty}    f_{\mathbf{X}}(y \varphi(\mathbf{u}_0))\, y^{d-1} J(1,\mathbf{u}_0) \mathrm{d}y \ \text{ a.e. }
\end{equation*}
thus, simplifying, we obtain
\begin{equation*}
    f_{\mathbf{X}}(t_0 \varphi(\mathbf{u}_0))  = \!\!\!\!\!
    \int\limits_{\varphi(\mathbf{u}) \in \mathcal{S}_1} \!\!\!\!\! f_{\mathbf{X}}(t_0 \varphi(\mathbf{u}))  J(1,\mathbf{u})  \mathrm{d}\mathbf{u}   \times \int\limits_0^{\infty}    f_{\mathbf{X}}(y \varphi(\mathbf{u}_0)) y^{d-1} \mathrm{d}y =:h_1(t_0) \times h_2( \varphi(\mathbf{u}_0)),
\end{equation*}
proving the existence of the function $h_1$ and $h_2$ satisfying \eqref{eq:factorization}.
\end{proof}
Let us end with a remark.
It appears that assuming that the function related to the direction $\mathbf{X}/ \|\mathbf{X}\|$ (here, $h_2$) is a constant, does not guarantee that the distribution of $\mathbf{X}/ \|\mathbf{X}\|$ is uniform on the unit sphere $\mathcal{S}_1$. Nevertheless, for $L^1$, $L^2$ and $L^\infty$ norms (or their weighted versions), this distribution is uniform on the unit sphere $\mathcal{S}_1$:
\begin{equation*}
    \p{\mathbf{X}/ \|\mathbf{X}\| \in B} = \frac{\mu(B)}{\mu(\mathcal{S}_1)},
\end{equation*}
where $\mu$ is a measure on a sphere. But for $L^p$-norms with $p \neq 1,2,\infty$,  it is not uniform, as a measure on the unit sphere is not proportional to a measure on the unit ball.  

\paragraph{Using integration on sphere: an illustration with the norm's cdf -}

Assume the random vector $\mathbf{X}$ has a  pdf defined, for any $\mathbf{x}\in\R^d$, by 
\begin{equation}\label{def:exampledensity}
f(\mathbf{x}) = h_1(\|\mathbf{x}\|)\,h_2(\mathbf{x}/\|\mathbf{x}\|) \quad\text{with} \quad 
h_1(\|\mathbf{x}\|):=C_d\,(1+ \|\mathbf{x}\|)^{-(\alpha+d)},
\end{equation}
$C_d > 0$ being the normalizing constant (making $f$ a pdf). 

Note that, if $h_2$ is a constant, then we get back the multivariate Pareto pdf considered in the example developed in Section~\ref{ssec:exs-dNormex}.

Changing to polar coordinates \eqref{eq:change_of_var}, setting 
 $r = \|\mathbf{x}\|$ and  $\theta = \frac{\mathbf{x}}{\|\mathbf{x}\|}$, we can write
\begin{equation*}
    \frac{1}{C_d} = \int_{\mathbb{R}^d} \frac{h_2(\mathbf{x}/\|\mathbf{x}\|)}{(1+ \|\mathbf{x}\|)^{\alpha+d}} \,\mathrm{d}\mathbf{x} = c_d \int_{0}^{\infty} \frac{r^{d-1}}{(1+r)^{\alpha+d}}\,\mathrm{d}r  = c_d \,\frac{\Gamma(d) \Gamma(\alpha)}{\Gamma(\alpha + d)}, 
\end{equation*}
where $c_d$ is defined after  \eqref{eq:jointdistrnorm_direc}. 
Then the cdf of  $\|\mathbf{X}\|$ can be expressed as
\begin{equation}\label{eq:cdf-norm}
    F_{\|\mathbf{X}\|}(t)= \p{\|\mathbf{X}\| \leq t} = C_d\, c_d \int_{0}^{t} r^{d-1}\,\frac{1}{(1+r)^{\alpha+d}}\, \mathrm{d}r  = 1-\sum_{k=0}^{d-1} \frac{\Gamma(\alpha +k)}{\Gamma(k+1)\Gamma(\alpha)}\frac{t^{k}}{(1+t)^{\alpha+k}},
\end{equation}
proceeding successively by an integration by parts and using the fact that, for any $m \in \mathbb{N},$ $\alpha > 0$ such that $\alpha > m + 1$, 
\begin{equation*}
\int_0^{\infty} \frac{x^m}{(1+x)^{\alpha}} \mathrm{d}x 
= \frac{\Gamma(m+1) \Gamma(\alpha - (m + 1))}{\Gamma(\alpha)},   
\end{equation*}
to obtain the last equality in \eqref{eq:cdf-norm}.

Note that this result is independent of the norm choice and it will be the same when applying  Normex method to any $\mathbf{X}$ having its pdf satisfying \eqref{def:exampledensity}.

The result \eqref{eq:cdf-norm} will be useful for applying both versions of multi-normex, $d$-Normex and MRV-Normex, in the various examples considered in the paper.

\section{Geometric notion of multivariate quantiles}\label{A:secGeomQuant}

In this section, we give an overview of the notion of multivariate quantiles (see e.g. \citet{chaudhuri:1996,DharAl:2014}) and QQ-plots.

\subsubsection*{One dimensional case}

In the one dimensional case, the quantile $q_{\alpha}(X)$ of level $\alpha \in (0,1)$ for a random variable $X \in \mathbb{R}$ can be defined as follows
\begin{equation}\label{def:quantile}
  q_{\alpha}(X) := \arg \min\limits_{ q \in \mathbb{R}} \skk{ \e\left( |X-q| - |X|\right) - q(2\alpha - 1)}.
\end{equation}
The proof of \eqref{def:quantile} is easy to obtain by the following representation for the expectation:
\begin{equation*}
  \e X = \int_{0}^{\infty} \bar{F}(t) \mathrm{d}t - \int_{-\infty}^{0} F(t) \mathrm{d}t.
\end{equation*}
If the first moment of $X$ exists, the term $|X|$ can be dropped out in Definition~\ref{def:quantile}, hence we have
\begin{equation*}
  q_{\alpha}(X) := \arg \min\limits_{ q \in \mathbb{R}} \skk{ \e |X-q|  - q(2\alpha - 1)}, \quad \a\in (0,1).
\end{equation*}
Sometimes, it is convenient to parameterize the levels $\alpha \in (0,1)$  by the set $(-1,1)$ with the linear transformation $u= 2 \alpha - 1$, so that the median ($\a=1/2)$ corresponds to $u=0$. In such a case, we have
\begin{equation*}
  q_{u}(X) := \arg \min\limits_{ q \in \mathbb{R}} \skk{ \e |X-q|  - u q}, \ \ u \in (-1,1).
\end{equation*}

\subsubsection*{Multidimensional case}

Let $F$ denote the multivariate cdf of the random vector $\mathbf{X}\in \mathbb{R}^d$, with $d\ge 1$, $\langle\cdot ,\cdot\rangle$ be the Euclidean inner product and $\|\cdot\|$ the Euclidean norm induced by the inner product.

\begin{defi}[Geometrical quantile, see \citet{chaudhuri:1996}]
\label{def:geom quantile}
For a random vector $\mathbf{X}$ with a probability distribution $F$ on $\mathbb{R}^{d}$, the {\bf $\mathbf d$-dimensional spatial quantile}  or  {\bf geometrical quantile} (GQ) $Q_{F}(\mathbf{u})=\left(Q_{F, 1}(\mathbf{u}), \ldots, Q_{F, d}(\mathbf{u})\right)$ is defined as
\begin{equation}\label{def:multi-quantile}
Q_{F}(\mathbf{u})=\arg \min _{\mathbf{Q} \in \mathbb{R}^{d}} \e \{ \| \mathbf{X}-\mathbf{Q}\| - \| \mathbf{X}\|  - \langle\mathbf{u},\mathbf{Q} \rangle\},
\end{equation}
with $ \mathbf{u} \in B^{d}=\left\{ \mathbf{v} \in \mathbb{R}^{d},\|\mathbf{v}\|<1\right\}.$ 
\end{defi}
If $\mathbf{u} = \mathbf{0}, $ then $Q_{F}(\mathbf{0})$ is a spatial median - a point from $\mathbb{R}^d$ that is equidistant from all values of a random vector $\mathbf{X}$.

GQs satisfy the equivariant property under rotation and shifting:
For deterministic vector $\mathbf{a} \in \mathbb{R}^d$ and orthogonal matrix $A$  we can find a GQ of $A\mathbf{X} + \mathbf{a}$ as
\begin{equation*}
  \mathbf{Q}_{A\mathbf{X} + \mathbf{a}}(\mathbf{u})  = A \mathbf{Q}_{\mathbf{X}}(A^{\mathrm{T}}\mathbf{u}) + \mathbf{a},
\end{equation*}
where $\mathbf{Q}_{\mathbf{X}}(\mathbf{u})$  is a GQ of $\mathbf{X}$.

Unfortunately,  GQs  are not equivariant under arbitrary affine transformations. In \citet{Chakrab:2001}, 
the author used a transformation-retransformation approach to construct affine equivariant estimates of multivariate quantiles, but the object became quite complicated and untractable.

Nevertheless, GQs characterize the distribution,
namely, if two random vectors $\mathbf{X}$ and $\mathbf{Y}$ yield the same quantile function $q$, then they have the same distribution (see \citet{Koltchinskii1997}). 
For extreme multivariate quantiles, this characterization does not hold true whenever $\e \|\mathbf{X}\|^2 < \infty$: The GQs converge to $\infty$ as the length $\lambda = \|\mathbf{u}\| \ra 1$  with the rate $(1 - \lambda)^{-1/2}$, more precisely
$$
\|\mathbf{Q}(\lambda \mathbf{u})\|^{2}(1-\lambda) \rightarrow \frac{1}{2}\left(\operatorname{tr} \Sigma-u^{\mathrm{T}} \Sigma u\right) \text { as } \lambda \uparrow 1,
$$
where $\Sigma$ denotes the covariance matrix of $\mathbf{X}$ (see \citet{girard2015b}). Thus the notion of GQ is not relevant for extreme multivariate quantiles of $L^p$-vectors with $p\ge2$, since for different distributions having the same covariance matrix,  GQs equivalently converge to $\infty$. But, in the case $\e \|\mathbf{X}\|^2 = \infty$, the same authors prove in \citet{girard2015b} that, for MRV distributions with parameter $\alpha \in (0,2]$, the asymptotics of GQs depend on the regular function introduced in the definition of MRV~\eqref{def:mrv}, making then the GQs more relevant.

 If all first moments of the random vector $\mathbf{X}$ exist, then Equation~\eqref{def:multi-quantile} can be written as
\begin{equation*}
Q_{F}(\mathbf{u})=\arg \min _{\mathbf{Q} \in \mathbb{R}^{d}} \e \left[ \| \mathbf{X}-\mathbf{Q}\|   - \langle\mathbf{u},\mathbf{Q} \rangle\right].
\end{equation*}
 Considering a sample $\mathbf{X}_1 , \dots, \mathbf{X}_n$ of iid random vectors of $\mathbb{R}^d$ with distribution $F$, $Q_{F}(\mathbf{u})$ can be estimated by the empirical  quantile
\begin{equation}\label{def:emp-mlti-quantile}
\hat{Q}_{F}(\mathbf{u})=\arg \min _{\mathbf{Q} \in \mathbb{R}^{d}} \frac{1}{n}\sum\limits_{i=1}^n \left( \| \mathbf{X_i}-\mathbf{Q}\|   - \langle\mathbf{u},\mathbf{Q} \rangle\right).
\end{equation}
In this case, the notion of {\bf spatial rank} may be introduced for each element $\mathbf{X}_j$ of the sample, namely
\begin{equation*}
 \mathbf{u}_j :=  n^{-1} \sum_{i \neq j}\left\|\mathbf{X}_j-\mathbf{X}_{i}\right\|^{-1}\left(\mathbf{X}_j-\mathbf{X}_{i}\right).
\end{equation*}
It is straightforward to verify that the GQ $\hat{Q}_{F}(\mathbf{u}_j)$ evaluated at the spatial rank $ \mathbf{u}_j$ is equal to the appropriate element $\mathbf{X}_j$ of the sample:
\begin{equation}\label{eq:spatial rank}
\hat{Q}_{F}(\mathbf{u}_j) = \mathbf{X}_j.
\end{equation}
This is a key property in the construction of QQ-plots of multivariate quantiles for small samples.

To compute the GQs, an algorithm has been proposed in \citet{chaudhuri:1996}, based on the  optimization problem \eqref{def:emp-mlti-quantile},  while another algorithm to construct QQ-plots has been developed in \citet{DharAl:2014}.

\subsubsection*{Multivariate QQ-plot}

To construct QQ-plots, we follow the algorithm from \citet{DharAl:2014} but with a small modification, given the fact that we use $10^7$ simulations to evaluate the considered  distribution. Instead of considering a number of points (in the plot) equal to the general volume of compared samples, as suggested in \citet{DharAl:2014}, which in our case would represent a huge amount of points, we decide to fix a set of 'levels' $\mathbf{u} \in B^{d}$ and compare multi-quantiles from different distributions for these levels.
\\[.7ex] More precisely,  we consider two independent $d$-dimensional data sets, namely, $\mathcal{X}=\left\{\mathbf{x}_{1}, \ldots, \mathbf{x}_{n}\right\}$ and $\mathcal{Y}=\left\{\mathbf{y}_{1}, \ldots, \mathbf{y}_{m}\right\}$, where $\mathbf{x}_{i}=$
$\left(x_{i, 1}, \ldots, x_{i, d}\right)$  are the realizations of a random sample with underlying distribution $F$,
and $\mathbf{y}_{j}=\left(y_{j, 1}, \ldots, y_{j, d}\right)$ with distribution $G$. Let $\mathcal{U} = \skk{\mathbf{u}_{1}, \ldots, \mathbf{u}_{n}}$ be a set of vectors from $B^{d}$. We compute
$Q_{\mathcal{X}}\left(\mathbf{u}_{k}\right)$  and $Q_{\mathcal{Y}}\left(\mathbf{u}_{k}\right)$ for $k=1, \ldots, n$, using the algorithm from \citet{chaudhuri:1996}. Then, we can match the two sets of quantiles by setting the correspondence between $Q_{\mathcal{X}}\left(\mathbf{u}_{k}\right)$ and $Q_{\mathcal{Y}}\left(\mathbf{u}_{k}\right)$ for $k=1, \ldots, n$. Finally, we construct the QQ-plots  as a collection of $d$ $2$-dimensional plots, where each plot corresponds to a component of the spatial quantile. In \citet{DharAl:2014}, Theorem 2.2, the authors proved that for all $i = 1,\cdots ,d$, the points in the
$i$-th $2$-dimensional plot will lie close to a straight line with slope $1$ and intercept $0$ if and only if $F = G$.





\subsection*{Supplementary Material. Truncated moments for certain distributions to apply multi-normex }\label{A:moments_comp} 
%

\subsubsection*{Example~\ref{ex:ParetoClayton} with Clayton copula}
Here we provide truncated moments for the example~\ref{ex:ParetoClayton} in the case $d=2, \alpha\theta = 1$ and $\|\cdot\|_{\infty}$ norm in programming language to be able to apply multi-normex methods.  

$\mathbb{E} \left( X_1; \|\mathbf{X}\|_{\infty} \leq y \right) = $
\begin{verbatim*}
((((-2*a-2)*sqrt(y+1)+(2*a^2+a+1)*y+2*a^2+a+2)*exp(-2*a*log(2*sqrt(y+1)-1))
+(exp(-a*log(y+1))*(((2*sqrt(y+1)-1)^(2*a)*((2*a-1)*(y+1)^a+2*a^2-3*a+1)+(-
2*a^2+3*a-1)*(y+1)^a)*exp(a*log(y+1))+(y+1)^a*(2*sqrt(y+1)-1)^(2*a)*((2-2*a
)*sqrt(y+1)+(-2*a^2+a-1)*y-2*a^2+a-1)))/((y+1)^a*(2*sqrt(y+1)-1)^(2*a)))/(2
*a^2-3*a+1)) 
\end{verbatim*}
$\mathbb{E} \left( X_1^2; \|\mathbf{X}\|_{\infty} \leq y \right) = $
\begin{verbatim}
(((2*a-3)*(2*a-1)*((2*sqrt(y+1)-1)^(2*a)*(2*(y+1)^a-(a-2)*(a-1))+(a-2)*(a-1)
*(y+1)^a)*exp(a*log(y+1))-(y+1)^(a+1)*(2*sqrt(y+1)-1)^(2*a)*(12*(a-2)*sqrt(y
+1)+(4*a^4-12*a^3+11*a^2-3*a+6)*y-4*a^4+20*a^3-19*a^2-19*a+30))*exp(2*a*log(
2*sqrt(y+1)-1))-(y+1)^a*(-y*((4*a^4-4*a^3+11*a^2+7*a+6)*y+2*(8*a^3+10*a^2+23
*a+18))+4*sqrt(y+1)*((a+1)*(2*a^2+a+6)*y+8*a^2+11*a+6)+4*a^4-20*a^3-5*a^2-53
*a-24)*(2*sqrt(y+1)-1)^(2*a)*exp(a*log(y+1)))*exp(-a*(2*log(2*sqrt(y+1)-1)+l
og(y+1))))/((a-2)*(a-1)*(2*a-3)*(2*a-1)*(y+1)^a*(2*sqrt(y+1)-1)^(2*a))
\end{verbatim}
$\mathbb{E} \left( X_1 X_2; \|\mathbf{X}\|_{\infty} \leq y \right)=$
\begin{verbatim*}
(((((2*sqrt(y+1)-1)^(2*a)*((4*a^3-10*a^2+5*a)*(y+1)^a+(-4*a^4+16*a^3-17*a^2
-a+6)*sqrt(y+1)+(4*a^5-16*a^4+17*a^3+a^2-6*a)*y+4*a^4-16*a^3+17*a^2+a-6)+(y
+1)^a*((4*a^4-16*a^3+17*a^2+a-6)*sqrt(y+1)+(-4*a^5+16*a^4-17*a^3-a^2+6*a)*y
-4*a^4+16*a^3-17*a^2-a+6))*exp(a*log(y+1))+(y+1)^a*((-4*a^3+10*a^2-6*a)*y^2
+sqrt(y+1)*((-8*a^4+28*a^3-28*a^2+8*a)*y+4*a^4-24*a^3+41*a^2-15*a-6)+(-4*a^
5+16*a^4-21*a^3+5*a^2+6*a)*y-4*a^4+16*a^3-21*a^2+5*a+6)*(2*sqrt(y+1)-1)^(2*
a))*exp(2*a*log(2*sqrt(y+1)-1))+(y+1)^a*((4*a^5-4*a^4-5*a^3-a^2+6*a)*y^2+sq
rt(y+1)*((-8*a^4+4*a^3+12*a^2-8*a)*y-4*a^4-13*a^2+15*a+6)+(4*a^5+a^3-3*a^2-
6*a)*y+4*a^4+4*a^3+3*a^2-10*a-6)*(2*sqrt(y+1)-1)^(2*a)*exp(a*log(y+1)))*exp
(-2*a*log(2*sqrt(y+1)-1)-a*log(y+1)))/(a*(4*a^4-20*a^3+35*a^2-25*a+6)*(y+1)
^a*(2*sqrt(y+1)-1)^(2*a))
\end{verbatim*}

\subsubsection*{Multivariate Pareto-Lomax with \texorpdfstring{$L^\infty$}{Lg}-norm}

In this section we provide an example of Multivariate Pareto-Lomax distribution, defined by $L^\infty$-norm,  and calculations of truncated (by $L^\infty$-norm as well) moments for an arbitrary dimension $d \geq 1$. One can use this calculations to apply multi-normex methods.  Let the density function of the initial random vector is
\begin{equation*}
    f_{\mathbf{X}}(\mathbf{x}) := \frac{C_d}{(1 + \|\mathbf{x}\|_{\infty})^{\alpha + d}},  x_i > 0, i =1, \cdots, d.
\end{equation*}
As usual, parameter $\alpha > 0,$ and for the constant $c_d$ we have
\begin{equation*}
    \frac{1}{C_d} = \iint_{\mathbb{R}^d_+} f_{\mathbf{X}}(\mathbf{x}) \mathrm{d}\mathbf{x} = d \int_{0}^{\infty} \frac{x^{d-1}}{(1+x)^{\alpha+d}} \mathrm{d}x = \frac{\Gamma(d+1) \Gamma(\alpha)}{\Gamma(\alpha + d)} .
\end{equation*}

\begin{equation*}
\begin{split}
  \mathbb{E} X_1 & =  C_d \iint_{\mathbb{R}_{d}^{+}} \frac{x_{1}}{\left(1+\max_{i} x_{i}\right)^{\alpha+d}} d \mathbf{x} \\ 
  & = C_d\left[\int_{x_{1}-\max_{i} x_{i}} \frac{x_{1}}{\left(1+x_{1}\right)^{\alpha+d} } \mathrm{d} \mathbf{x} +(d-1) \iint_{x_{2}=\max_{i}  x_{i}} \frac{x_{1}}{\left(1+x_{2}\right)^{\alpha+d}} \mathrm{d} \mathbf{x}\right] \\
  & = C_d   \left[\int_{0}^{\infty} \frac{x_{1}^{d}}{\left(1+x_{1}\right)^{\alpha+d}} \mathrm{d}x_1 +\left(d -1 \right) \int_{0}^{\infty} \frac{1}{2} \frac{x_{2}^{d}}{\left(1+x_{2}\right)^{\alpha+d}} \mathrm{d} x_{2}\right]
    \\ 
    & = C_d \cdot\left(1+\frac{(d-1)}{2}\right) \int_{0}^{\infty} \frac{x_{1}^{d}}{\left(1+x_{1}\right)^{\alpha+d}} d x_{1} = \frac{d+1}{2(\alpha-1)}.
  \end{split}
\end{equation*}
Using the same arguments 
\begin{equation*}
    \begin{split}
    \mathbb{E} \left( X_1; \|\mathbf{X}\|_{\infty} \leq t \right) & = C_d \cdot\left(1+\frac{(d-1)}{2}\right) \int_{0}^{t} \frac{x_{1}^{d}}{\left(1+x_{1}\right)^{\alpha+d}} d x_{1}  \\
    & = \left( \frac{d+1}{2}\right)\left( \frac{1}{\alpha - 1} -\sum_{k=0}^{d} \frac{\Gamma(\alpha +k-1)}{\Gamma(k+1)\Gamma(\alpha)}\frac{t^{k}}{(1+t)^{\alpha+k - 1}} \right).
    \end{split}
\end{equation*}
\begin{equation*}
\begin{split}
  \mathbb{E} X_1^2 &=  C_d \iint_{\mathbb{R}_{d}^{+}} \frac{x_{1}^2}{\left(1+\max_{i} x_{i}\right)^{\alpha+d}} d \mathbf{x} \\
    & = C_d \cdot\left(1+\frac{(d-1)}{3}\right) \int_{0}^{\infty} \frac{x_{1}^{d+1}}{\left(1+x_{1}\right)^{\alpha+d}} d x_{1} = \frac{(d+1)(d+2)}{3(\alpha-1)(\alpha-2)}.
  \end{split}
\end{equation*}
\begin{equation*}
    \begin{split}
    \mathbb{E} \left( X_1^2; \|\mathbf{X}\|_{\infty} \leq t \right) & = C_d \cdot\left(1+\frac{(d-1)}{3}\right) \int_{0}^{t} \frac{x_{1}^{d+1}}{\left(1+x_{1}\right)^{\alpha+d}} d x_{1}  \\
    & =  \frac{(d+1)(d+2)}{3} \left( \frac{1}{(\alpha - 1)(\alpha-2)} -\sum_{k=0}^{d+1} \frac{\Gamma(\alpha +k-2)}{\Gamma(k+1)\Gamma(\alpha)}\frac{t^{k}}{(1+t)^{\alpha+k - 2}} \right).
    \end{split}
\end{equation*}
\begin{equation*}
\begin{split}
  \mathbb{E}X_1 X_2 & =  C_d \iint_{\mathbb{R}_{d}^{+}} \frac{x_{1} x_{2}}{\left(1+\max_{i} x_{i}\right)^{\alpha+d}} d \mathbf{x} \\ 
  & = C_d\left[ 2 \int_{x_{1}=\max_{i} x_{i}} \frac{x_{1} x_{2}}{\left(1+x_{1}\right)^{\alpha+d} } \mathrm{d} \mathbf{x} +(d-2) \iint_{x_{3}=\max_{i}  x_{i}} \frac{x_{1} x_2}{\left(1+x_{3}\right)^{\alpha+d}} \mathrm{d} \mathbf{x}\right] \\
  & = C_d   \left[\int_{0}^{\infty} \frac{x_{1}^{d+1}}{\left(1+x_{1}\right)^{\alpha+d}} \mathrm{d}x_1 +\left(d -2 \right) \int_{0}^{\infty} \frac{1}{4} \frac{x_{3}^{d+1}}{\left(1+x_{3}\right)^{\alpha+d}} \mathrm{d} x_{3}\right]
    \\ 
    & = C_d \cdot\left(1+\frac{(d-2)}{4}\right) \int_{0}^{\infty} \frac{x^{d+1}}{\left(1+x\right)^{\alpha+d}} \mathrm{d} x = \frac{(d+1)(d+2)}{4(\alpha-1)(\alpha-2)}.
  \end{split}
\end{equation*}
\begin{equation*}
    \begin{split}
    \mathbb{E} \left( X_1 X_2; \|\mathbf{X}\|_{\infty} \leq t \right) & = C_d \cdot\left(1+\frac{(d-2)}{4}\right) \int_{0}^{t} \frac{x_{1}^{d+1}}{\left(1+x_{1}\right)^{\alpha+d}} d x_{1}  \\
    & =  \frac{(d+1)(d+2)}{4} \left( \frac{1}{(\alpha - 1)(\alpha-2)} -\sum_{k=0}^{d+1} \frac{\Gamma(\alpha +k-2)}{\Gamma(k+1)\Gamma(\alpha)}\frac{t^{k}}{(1+t)^{\alpha+k - 2}} \right).
    \end{split}
\end{equation*}

\end{appendix}

\end{document}